\documentclass{amsart}    

\usepackage{amssymb,graphicx}
\usepackage{amsfonts}
\usepackage{amsmath}
\usepackage{amsthm}
\usepackage{fancyhdr}
\usepackage{indentfirst}
\usepackage{comment}
\usepackage{color}
\usepackage{subcaption}
\usepackage{epsfig}
\usepackage{pst-grad} 
\usepackage{pst-plot} 
\usepackage[space]{grffile}
\usepackage{geometry}
\usepackage{fancyhdr}
\newtheorem{thm}{Theorem}[subsection]
\newtheorem{conjecture}{Conjecture}[subsection]
\newtheorem{corollary}{Corollary}[subsection]
\newtheorem{lemma}{Lemma}[subsection]
\newtheorem{definition}{Definition}[subsection]
\usepackage[marginal]{footmisc}

\begin{document}

\title{The Space of Partially Free Boundary Minimal Half Disks}

\author{Shanjiang Chen}
\address{Department of Mathematics, The Chinese University of Hong Kong, Shatin, N.T., Hong Kong}
\email{sjchen@math.cuhk.edu.hk}


\begin{abstract}
This paper forms part of our ongoing works on the existence of complete non-compact free boundary minimal planes in an asymptotically flat three-dimensional Riemannian manifold with boundary.
We set up the degree theory for the space of properly embedded partially free boundary minimal half disks in a three-dimensional half ball. We prove that the space of such surfaces is a Banach manifold. The projection map which projects partially free boundary minimal half disks into their Dirichlet boundaries is a Fredholm map of index zero and has a well-defined mod-2 degree under suitable assumptions. The new major analytic difficulty is to analyze the properties of a second order elliptic operator with mix boundary conditions on a domain with corners.
\end{abstract}
\maketitle

\tableofcontents

\section{Introduction}
In Chodosh-Ketover\cite{ChodoshKetover2018}, they show that in an asymptotically flat three-dimensional Riemannian manifold containing no closed minimal surface, there is a complete non-compact properly embedded minimal plane passing through a fixed point. In the introduction to \cite{ChodoshKetover2018}, they discuss the motivation of the existence theory of complete non-compact minimal surface and the reason why they apply the mod-2 degree theory developed by White\cite{White1987} instead of min-max theory or solving the Plateau problems in coordinate balls, and then taking limits. In \cite{MazetLaurent2018Mpia}, Mazet and Rosenberg improve Chodosh and Ketover's result. They show that in an asymptotically flat three-dimensional Riemannian manifold containing no closed minimal surface, there is a complete non-compact properly embedded minimal plane passing through a fixed point which is also orthogonally to a given unit vector at that fixed point, and a complete non-compact properly embedded minimal plane passing through three fixed points.

We are interested in the existence of complete non-compact free boundary minimal planes in an asymptotically flat three-dimensional Riemannian manifold with boundary. In particular, we have the following conjecture.
\begin{conjecture}
Let $M = (\mathbb R_{+}^{3},g)$ be an asymptotically flat three-dimensional Riemannian manifold with boundary.
Suppose $M$ contains no compact free boundary minimal surface. Then for any $p \in \partial M$, there exists a complete non-compact free boundary minimal plane passing through $p$.
\end{conjecture}

The assumption that $M$ contains no compact free boundary minimal surface is related to the area and curvature estimates for partially free boundary minimal surfaces when we take limits.

Similar to Chodosh-Ketover\cite{ChodoshKetover2018}, we choose to apply the degree theory to the coordinate half balls to prove the existence of partially free boundary minimal half disks passing through a fixed point. Then we prove the existence of a  convergent subsequence as the radius tend to infinity.

In this paper, we set up the degree theory for the space of partially free boundary minimal half disks in a three-dimensional half ball.

Let $N = (\mathbb B_{+}^3, g)$ be a three-dimensional Riemannian manifold with weakly mean convex boundaries consisting of the half sphere $S$ and the disk $\Sigma$ such that $S$ meets $\Sigma$ orthogonally along their common boundary.
Let $D_{+} = \{(x,y) : x^2+y^2 < 1, y>0\}$ be the unit half disk. $\partial D_{+} = \overline{\Gamma_{1}} \cup \overline{\Gamma_{2}}$ where $\Gamma_{1} = \{(x,y) : x^2+y^2=1, y>0\}$ and $\Gamma_{2} = \{(x,0) : x \in (-1,1)\}$.

We use $C_{\perp}^{2,\sigma}(\overline{\Gamma_{1}}, S) = \{\gamma:\Gamma_{1} \to S \, \big| \, \gamma$ is an embedding and $C^{2,\sigma}$ up to the boundary, $\gamma(\Gamma_{1}) \subset S, \partial [\gamma(\Gamma_{1})] = \overline{\gamma(\Gamma_{1})} \cap \partial S, \gamma(\Gamma_{1}) \perp \Sigma \}$ to denote the space of the Dirichlet boundaries.

We consider the space $C_{*}^{2,\sigma}(\overline{D_{+}}, N) = \{f: D_{+} \to N \, \big| \, f$ is $C^{2,\sigma}$ up to the boundaries and corners, $f(\Gamma_{1}) \subset S, f(\Gamma_{2}) \subset \Sigma, f(D_{+}) \subset N$, $f(D_{+})$ is transversal to $\partial N, f|_{\Gamma_{1}} \in C_{\perp}^{2,\sigma}(\overline{\Gamma_{1}}, S)\}$.

Two maps $f_{1}, f_{2} \in C_{*}^{2,\sigma}(\overline{D_{+}}, N)$ are said to be equivalent if $f_{1} = f_{2} \circ u$ for some $C^{2,\sigma}$ diffeomorphism $u: D_{+} \to D_{+}$ such that $u(x) = x$ for all $x \in \Gamma_{1}$. We use $[f]$ to denote the equivalence class of $f$.

Let $f \in C_{*}^{2,\sigma}(\overline{D_{+}}, N)$ and $\vec{X}$ be a variation on $f(D_{+})$ such that $\vec{X} \in TS,T\Sigma$ on $f(\Gamma_{1}), f(\Gamma_{2})$. Then standard calculations give the first variation formula
$$\dfrac{d}{dt}\Big|_{t=0}Area_{t} = -\int_{f(D_{+})}<\vec{H},\vec{X}> + \int_{f(\Gamma_{1})}<\vec{n},\vec{X}> + \int_{f(\Gamma_{2})}<\vec{n},\vec{X}>$$
where $\vec{H}$ is the mean curvature vector of $f(D_{+})$ in $N$, $\vec{n}$ is the outward unit normal vector of $f(D_{+})$ on $f(\Gamma_{1})$ and $f(\Gamma_{2})$. This is why we need $f$ to be $C^{2,\sigma}$ up to the boundaries and the corners.
\begin{definition}
We say $f(D_{+})$ is a partially free boundary minimal half disk if $\vec{H} = 0$,  $\nu|_{f(\Gamma_{2})} \perp \Sigma$. We use $\mathcal M = \{[f]: f \in C_{*}^{2,\sigma}(\overline{D_{+}}, N), \vec{H} = 0, \nu|_{f(\Gamma_{2})} \perp \Sigma\}$ to denote this space.
\end{definition}

In Section 3, we prove the following global structure theorem and the mod-2 degree theorem.
\begin{thm}
  Let $\mathcal M$ be defined as above. Then \\
  (1) $\mathcal{M}$ is a Banach manifold modelled on $C_{\perp}^{2,\sigma}(\overline{\Gamma_{1}}, S)$ with a countable family of coordinate charts.\\
  (2) The map $\Pi: \mathcal M \to C_{\perp}^{2,\sigma}(\overline{\Gamma_{1}}, S)$ defined by $\Pi([f]) = f|_{\Gamma_{1}}$ is a Fredholm map of index zero. \\
  (3) If $D\Pi([f_{0}])$ has $k$-dimensional kernel, then a neighborhood $U$ of $[f_{0}]$ in $\mathcal M$ may be identified with a codimension-$k$ submanifold of $C_{\perp}^{2,\sigma}(\overline{\Gamma_{1}}, S) \times \mathbb R^k$ (so that $\Pi$ corresponds to the obvious projection of $C_{\perp}^{2,\sigma}(\overline{\Gamma_{1}}, S) \times \mathbb R^k$).
  (4) $U$ may be chosen small enough such that there is a $k$-dimensional subspace $V$ of $C_{\perp}^{2,\sigma}(\overline{\Gamma_{1}}, S)$ and $V \times \{0\}$ complements the tangent space to $U$ at each point of $U$.
\end{thm}
\begin{thm}
  Let $\mathcal M^{'}$ be an open subset of $\mathcal M$, $Z$ be a connected open subset of $C_{\perp}^{2,\sigma}(\overline{\Gamma_{1}}, S)$. If $\Pi$ maps $\mathcal M^{'}$ properly into $Z$, then for a generic $\gamma \in Z$,
  $$\#(\Pi^{-1}(\gamma) \cap \mathcal M^{'}) \quad \text{mod 2}$$
  is a constant (1 or 0).
\end{thm}

The new major difficulty is to analyze the Jacobi operator with mixed boundary conditions (the Dirichlet boundary condition on $\Gamma_{1}$ and the Robin boundary condition on $\Gamma_{2}$) on the unit half disk which has two corners.

In Section 2, we extend Grisvard's results\cite{GrisvardPierre2011EPiN} on the Laplace operator on polygonal domains to general linear second order elliptic operators with general mix boundary conditions on curvilinear polygonal domains (see Definition 2.2.2).

In some special cases, we establish nice properties similar to the case when the domain is smooth. This is what we need in Section 3.
\begin{thm}
  Let $J(u) = a_{ij}D_{ij}u+ b_{i}u_{i} + cu$ be a strongly elliptic operator on $D_{+}$ such that $a_{ij}, b_{i}, c \in C^{0,\sigma}(\overline{D_{+}})$, $a_{ij}$ is diagonal at $S_{1} = (-1,0)$ and $S_{2} = (1,0)$. Let $\beta, d \in C^{1,\sigma}(\overline{\Gamma_{2}})$ such that $\beta(S_{1}) = \beta(S_{2}) = 0$.\\
  Let $C_{\Gamma_{1}}^{2,\sigma}(\overline{D_{+}}) = \{u \in C^{2,\sigma}(\overline{D_{+}}) : u|_{\Gamma_{1}} = 0\}$, $C_{0}^{1,\sigma}(\overline{\Gamma_{2}}) = \{g \in C^{1,\sigma}(\overline{\Gamma_{2}}) :  g(S_{1}) = g(S_{2}) = 0\}$, then the map
  $$L: C_{\Gamma_{1}}^{2,\sigma}(\overline{D_{+}}) \to C^{0,\sigma}(\overline{D_{+}}) \times C_{0}^{1,\sigma}(\overline{\Gamma_{2}})$$
  defined by
  $$L(u) = (J(u), [\dfrac{\partial u}{\partial \nu_{2}}+\beta\dfrac{\partial u}{\partial \tau_{2}} + du]|_{\Gamma_{2}})$$
   is a Fredholm map of index zero where $\nu_{2}$ is the outward unit normal vector of $D_{+}$ on $\Gamma_{2}$, $\tau_{2}$ is the unit tangential vector of $\Gamma_{2}$.
\end{thm}

\section*{acknowledgements}
The author would like to thank his advisor Professor Martin Man Chun LI for introducing him to the field of geometric analysis and all the suggestions and encouragements in his Ph. D. study.

\section{The Elliptic Theory on Polygonal Domains}
\subsection{Introduction}
The purpose of this section is to extend Grisvard's results on the Laplace operator to general elliptic operators on domains with corners.

In Section 2.2, we introduce the notations.

In Section 2.3, we state Grisvard's existence theorem for the Laplace operator with the Dirichlet or the oblique boundary conditions with constant coefficients.

In Section 2.4, we extend Grisvard's existence theorem to the Laplace operator with general boundary conditions.

In Section 2.5, we prove apriori estimates for general elliptic operators with general boundary conditions.

In Section 2.6, we find the index of general elliptic operators by the continuity method.

In Section 2.7, we consider some special cases and prove the existence of regular solutions using the results in Section 2.6.

In Section 2.8 and Section 2.9, we prove inverse trace theorems and deformation theorems in order to extend Grisvard's apriori estimates for homogeneous boundary values to non-homogeneous boundary values in a curvilinear polygonal domain. (When the domain is smooth, this part is trivial.)

\subsection{Notations}
\begin{definition}
Let $\Omega$ be a bounded open subset in $\mathbb R^2$. We say $\Omega$ is a polygonal domain if $\partial \Omega = \Gamma$ is a polygon.\\
For $j = 1,2,...,I_{0}$, we have the following notations:\\
(1) $\Gamma_{j}$: $\Gamma$ is the union of a finite number of linear segments $\overline{\Gamma_{j}}$ (we assume that $\Gamma_{j}$ is open).\\
(2) $S_{j}$: $S_{j} = \overline{\Gamma_{j}} \cap \overline{\Gamma_{j+1}}$. \\
(3) $\omega_{j}$: $\omega_{j}$ is the interior angle of $\Gamma$ at $S_{j}$. \\
(4) $\mathcal D$, $\mathcal N$: We divide $\Gamma_{j}$ into two parts. We consider the Dirichlet boundary condition on $\Gamma_{j}$ for $j \in \mathcal D$ and the Robin boundary condition on $\Gamma_{j}$ for $j \in \mathcal N$.\\
(5) $\nu_{j}$, $\tau_{j}$: $\nu_{j}$ is the outer unit normal vector of $\Omega$ on $\Gamma_{j}$, $\tau_{j}$ is the unit tangent vector of $\Gamma_{j}$. \\
(6) $\beta_{j}$,$\phi_{j}$: $\beta_{j} \in C^{1,\sigma}(\overline{\Gamma_{j}})$.\\
 $\phi_{j}(S_{k}) = \begin{cases}
  \dfrac{\pi}{2}                    & j \in \mathcal D, k \in \{j,j+1\},\\
   \arctan \beta_{j}(S_{k})        & j \in \mathcal N, k \in \{j,j+1\}.
\end{cases}$\\
We will use $\phi_{j}$ instead of $\phi_{j}(S_{k})$ if $\beta_{j}$ is a constant.\\
(7) $\mu_{j}$: $\mu_{j} = \begin{cases}
  \tau_{j}                    & j \in \mathcal D,\\
  \nu_{j}+\beta_{j}\tau_{j}   & j \in \mathcal N.
\end{cases}$\\
(8) $C^{2,\sigma}(\overline{\Omega})$: For any $\sigma \in (0,1)$, a function $f \in C^{2,\sigma}(\overline{\Omega})$ if and only if $f,Df,D^{2}f$ are continuous and bounded in $\overline{\Omega}$ and $D^{2}f$ are uniformly H\"{o}lder continuous with exponent $\sigma$.\\
(9) $C_{c}^{k,\sigma}(\Omega)$, $C_{c}^{k,\sigma}(\Gamma_{j})$, $C_{c}^{k,\sigma}(\mathbb R^n)$:\\
$C_{c}^{k,\sigma}(\Omega) = \{u \in C^{k,\sigma}(\Omega): \text{supp}(u) \subset \Omega\}$,\\
$C_{c}^{k,\sigma}(\Gamma_{j}) = \{g \in C^{k,\sigma}(\Gamma_{j}): \text{supp}(g) \subset \Gamma_{j}\}$,\\
$C_{c}^{k,\sigma}(\mathbb R^n) = \{u \in C^{k,\sigma}(\mathbb R^n) : \text{supp}(u) \text{ is compact}\}$.
\end{definition}

We have a slightly different definition for curvilinear polygonal domains.
\begin{definition}
  Let $\Omega$ be a bounded open subset in $\mathbb R^2$. We say $\Omega$ is a $C^{k,\sigma}$ curvilinear polygonal domain if for any $x \in \Gamma = \partial \Omega$, there is a positive real number $\epsilon$ satisfying (a) or (b):\\
  (a) there is a mapping $\Psi = (\Psi_{1}, \Psi_{2})$ from $B_{\epsilon}(x)$ to $\mathbb R^2$ such that $\Psi$ is injective, $\Psi$ and $\Psi^{-1}$ (defined on $\Psi(B_{\epsilon}(x))$) is $C^{k,\sigma}$ and $\Omega \cap B_{\epsilon}(x) = \{p \in \Omega: \Psi_{2}(p) > 0\}$.\\
  (b) $\Omega \cap B_{\epsilon}(x) = \{x+(r\cos \theta, r\sin \theta):0 < r < \epsilon, \theta_{0} < \theta < \theta_{1}\}$ for some positive $\epsilon$ and $0<\theta_{1}-\theta_{0}<2\pi$.
\end{definition}

Thus curvilinear polygonal domains will include domains with one or two corners. For a curvilinear polygonal domain, we use the same notations as in Definition 2.2.1.

\subsection{Grisvard's results on the Laplace operator on polygonal domains}
Let $\Omega$ be a polygonal domain and $\beta_{j}$ be real numbers. Grisvard considers the homogeneous boundary value problem
$$(*)\begin{cases}
     \triangle u = f      & \text{in $\Omega$},\\
     u = 0         & \text{on $\Gamma_{j}$, $j \in \mathcal D$}, \\
     \dfrac{\partial u}{\partial \nu_{j}} + \beta_{j}\dfrac{\partial u}{\partial \tau_{j}}=0 & \text{on $\Gamma_{j}$, $j \in \mathcal N$},
  \end{cases}$$
the non-homogeneous boundary value problem
$$(**)\begin{cases}
     \triangle u = f      & \text{in $\Omega$},\\
     u = g_{j}         & \text{on $\Gamma_{j}$, $j \in \mathcal D$}, \\
     \dfrac{\partial u}{\partial \nu_{j}} + \beta_{j}\dfrac{\partial u}{\partial \tau_{j}}=g_{j} & \text{on $\Gamma_{j}$, $j \in \mathcal N$},
  \end{cases}$$
and he introduces the following singular solutions to the Laplace operator near the corners which are the real or imaginary part of differentiable complex-valued functions satisfying the given boundary conditions.

Let $\epsilon > 0$ be small enough such that $\Omega \cap B_{\epsilon}(S_{j})$ is a sector for any $j$ and $\eta_{j} \in C_{c}^{\infty}(\mathbb R^2)$ be a cut-off function such that $\eta_{j} = 1$ near $S_{j}$ and $\text{supp}(\eta_{j}) \subset B_{\epsilon}(S_{j})$.

We have the following singular solutions:\\
$\lambda_{j,m} = \dfrac{\phi_{j}-\phi_{j+1}+m\pi}{\omega_{j}}$ \quad \quad \quad for $j=1,2,...,I_{0}$ and $m \in \mathbb Z$,\\
$\mathcal S_{j,m}(r_{j}e^{i\theta_{j}}) = \begin{cases}
  r_{j}^{-\lambda_{j,m}}\cos (\lambda_{j,m}\theta_{j}+\phi_{j+1})\eta(r_{j}e^{i\theta_{j}}) \\
   \quad \quad \text{ if $\lambda_{j,m}$ is not an integer},\\
  r_{j}^{-\lambda_{j,m}}[\ln r_{j} \cos(\lambda_{j,m}\theta_{j}+\phi_{j+1}) +\theta_{j}\sin(\lambda_{j,m}\theta_{j}+\phi_{j+1})]\eta(r_{j}e^{i\theta_{j}})\\
   \quad \quad \text{     if $\lambda_{j,m}$ is an integer},
\end{cases}$\\
where $r_{j}$ and $\theta_{j}$ are the polar coordinates with the origin at $S_{j}$.
\begin{thm}[Grisvard\cite{GrisvardPierre2011EPiN}, Theorem 4.4.4.13]
Assume $\dfrac{\phi_{j} - \phi_{j+1} + \frac{2\omega_{j}}{q}}{\pi}$ is not an integer for any $j$. Then for each $f \in L_{p}(\Omega)$ ($\int_{\Omega}f = 0$ if $\mathcal D$ is empty) there exist real numbers $c_{j,m}$ and a (possibly non-unique) $u$ such that $u$ is a solution to problem ($*$) and
$$u - \sum_{-\frac{2}{q} < \lambda_{j,m} < 0, \lambda_{j,m} \neq -1}c_{j,m} \mathcal S_{j,m} \in W_{p}^{2}(\Omega)$$
where $\dfrac{1}{p} + \dfrac{1}{q} = 1$.
\end{thm}

The point of Theorem 2.3.1 is that the solutions to problem ($*$) may contain singular parts near the corners.

Next Grisvard proves the existence of solutions to problem ($**$) with higher orders and the non-homogeneous boundary value conditions.
\begin{thm}[Grisvard\cite{GrisvardPierre2011EPiN}, Theorem 5.1.3.5]
Assume that $\dfrac{\phi_{j} - \phi_{j+1} + \frac{2\omega_{j}}{q}}{\pi}$ is not an integer for any $j$. Then for each $f \in W_{p}^{k}(\Omega)$, $g_{j} \in W_{p}^{k+2-\frac{1}{p}}(\Gamma_{j})$ if $j \in \mathcal D$, $g_{j} \in W_{p}^{k+1-\frac{1}{p}}(\Gamma_{j})$ if $j \in \mathcal N$ and $g_{j}(S_{j}) = g_{j+1}(S_{j})$ if $j,j+1 \in \mathcal D$, there exist real numbers $c_{j,m}$ and a (possibly non-unique) $u$ such that $u$ is a solution to problem ($**$) and
$$u - \sum_{-k-\frac{2}{q} < \lambda_{j,m} < 0}c_{j,m} \mathcal S_{j,m} \in W_{p}^{k+2}(\Omega)$$
where $\dfrac{1}{p} + \dfrac{1}{q} = 1$.
\end{thm}

The difference between Theorem 2.3.1 and Theorem 2.3.2 is that in some cases (for instance $j \in \mathcal D, j+1 \in \mathcal N, \omega_{j} = \dfrac{\pi}{2}, \beta_{j+1} = 0$), if $u \in W_{p}^{k+2}$ is a solution to problem ($**$), then there are some compatibility conditions for $g_{j}$. Theorem 2.3.2 says that even though the given $g_{j}$ are not compatibility at the corners for some special $\omega_{j}$ and $\beta_{j}$, we can also get a similar result as in Theorem 2.3.1, but we need to add $\mathcal S_{j,m}$ where $\lambda_{j,m}$ is a negative integer.

Next Grisvard proves two similar theorems in the H\"{o}lder spaces.
\begin{thm}[Grisvard\cite{GrisvardPierre2011EPiN}, Theorem 6.4.2.5]
If $\mathcal D$ is not empty and at least two of the vectors $\mu_{j}$ are linearly independent. \\
Assume that $0< \sigma <1$ and $\dfrac{\phi_{j} - \phi_{j+1} + (2+\sigma)\omega_{j}}{\pi}$ is not an integer for any $j$. Then for each $f \in C^{0,\sigma}(\overline{\Omega})$, there exist real numbers $c_{j,m}$ and a function $u$ such that $u$ is a solution to problem ($*$) and
$$u - \sum_{-(2+\sigma) < \lambda_{j,m} < 0, \lambda_{j,m} \neq -1}c_{j,m} \mathcal S_{j,m} \in C^{2,\sigma}(\overline{\Omega}).$$
\end{thm}

Remark: The proof depends on Theorem 2.3.1, Theorem 2.3.2 and the Sobolev embedding theorem $W_{p}^{3}(\Omega) \subset C^{2,\sigma}(\overline{\Omega})$ if $p$ is large enough. When passing to a convergent subsequence, we need to guarantee that the solutions to problem ($*$) is unique, this is why Grisvard assume that $\mathcal D$ is not empty and at least two of the vectors $\mu_{j}$ are linearly independent.

\begin{thm}[Grisvard\cite{GrisvardPierre2011EPiN}, Theorem 6.4.2.6]
If $\mathcal D$ is not empty and at least two of the vectors $\mu_{j}$ are linearly independent.\\
Assume that $0< \sigma <1$ and $\dfrac{\phi_{j} - \phi_{j} + (2+\sigma)\omega_{j}}{\pi}$ is not an integer for any $j$.
Then for each $f \in C^{k,\sigma}(\overline{\Omega})$, $g_{j} \in C^{k+2,\sigma}(\overline{\Gamma_{j}})$ if $j \in \mathcal D$, $g_{j} \in C^{k+1,\sigma}(\overline{\Gamma_{j}})$ if $j \in \mathcal N$ and $g_{j}(S_{j}) = g_{j+1}(S_{j})$ if $j,j+1 \in \mathcal D$, there exist real numbers $c_{j,m}$ and a function $u$ such that $u$ is a solution to problem ($**$) and
$$u - \sum_{-(k+2+\sigma) < \lambda_{j,m} < 0}c_{j,m} \mathcal S_{j,m} \in C^{k+2,\sigma}(\overline{\Omega}).$$
\end{thm}

Observation: By Theorem 2.3.4, the map defined by
$$L(u) = (\triangle u, \prod_{j \in \mathcal D}u\big|_{\Gamma_{j}}, \prod_{j \in \mathcal N}[\dfrac{\partial u}{\partial \nu_{j}} + \beta_{j}\dfrac{\partial u}{\partial \tau_{j}}]\big|_{\Gamma_{j}})$$
from $C^{2,\sigma}(\overline{\Omega})$ to $C^{0,\sigma}(\overline{\Omega}) \times \prod_{j \in \mathcal D}C^{2,\sigma}(\overline{\Gamma}_{j}) \times \prod_{j \in \mathcal N}C^{1,\sigma}(\overline{\Gamma}_{j})$ is a Fredholm map. And the index of $L$ depends on the numbers of elements in $\{\lambda_{j,m}: -(2+\sigma) < \lambda_{j,m} < 0\}$.

When the domain is smooth, to extend the results on the Laplace operator to general elliptic operators, people first prove apriori estimates and then apply the continuity method.

But when the domain is not smooth, the index of $L$ will depend on the values of  $\beta_{j}$. This tells us it is impossible to find a constant $C = C(\Omega, \beta_{j})$ which depends on $\beta_{j}$ continuously such that
$$\|u\|_{2,\sigma, \overline{\Omega}} \leq C\{\|\triangle u\|_{\sigma, \overline{\Omega}} + \sum_{j \in \mathcal D}\|u\|_{2,\sigma,\overline{\Gamma_{j}}} +  \sum_{j \in \mathcal N}\|\dfrac{\partial u}{\partial \nu_{j}} + \beta_{j}\dfrac{\partial u}{\partial \tau_{j}}\|_{1,\sigma,\overline{\Gamma_{j}}}    \}.$$

To deal with this and extend Grisvard's results to general elliptic operators with general boundary conditions, we first prove a existence theorem similar to Theorem 2.3.4 for the Laplace operator with general boundary conditions and then consider general elliptic operators.

\subsection{General boundary conditions}
In this section, we assume that $\Omega$ is a curvilinear polygonal domain (Definition 2.2.2), $\beta_{j}(x)$ is a function on $\Gamma_{j}$ and there is a small positive $\epsilon$ such that $\beta_{j}(x) = \beta_{j}(S_{j})$ for $x \in B_{\epsilon}(S_{j}) \cap \overline{\Gamma_{j}}$, $\beta_{j}(x) = \beta_{j}(S_{j+1})$ for $x \in B_{\epsilon}(S_{j+1}) \cap \overline{\Gamma_{j}}$.

We first show that problem ($*$) has a solution $u \in H^{1}(\Omega)$, then we show that $u \in H^{2}(\Omega \setminus V)$ where $V$ is any closed neighborhood of the corners. Thus we can use the results in Section 2.3 to analyze the behaviour of $u$ near the corners.

To prove the existence theorems, we need two lemmas.
\begin{lemma}[Grisvard\cite{GrisvardPierre2011EPiN}, Theorem 1.5.1.1]
  Let $\Omega$ be a bounded open subset in $\mathbb R^{n}$ with a $C^{k-1,1}$ boundary $\Gamma$. Then the map $u \to u|_{\Gamma}$ defined on $C^{\infty}(\overline{\Omega})$ has a unique continuous extension from $W_{p}^{k}(\Omega)$ onto $W_{p}^{k-\frac{1}{p}}(\Gamma)$.
\end{lemma}

\begin{lemma}[Grisvard\cite{GrisvardPierre2011EPiN}, Theorem 4.4.4.1]
  Let $V$ and $W$ be two Hilbert spaces with a continuous injection of $W$ into $V$ and let $A$ be a continuous bilinear form on $V \times W$. Suppose there exists a constant $\alpha > 0$ such that
  $$A(w,w) \geq \alpha\|w\|^{2}_{V}$$
  for any $w \in W$. Then for each continuous linear form $l$ on $W$, there exists a (possibly non-unique) $u \in V$ such that
  $$A(u,w) = l(w)$$
  for any $w \in W$.
\end{lemma}

Next we prove the following existence theorem.
\begin{thm}
  Let $\Omega$ be a $C^{1,1}$ curvilinear polygonal domain, $\beta_{j}(x) \in C^{1}(\overline{\Gamma_{j}})$ and $d_{j} \in C^{0}(\overline{\Gamma_{j}})$ with $d_{j} \geq 0$.\\
  Then for each $f \in L_{p}(\Omega)$, there exists a (possibly non-unique) $u \in H^{1}(\Omega)$ which is a solution to
  $$(***)\begin{cases}
     \triangle u = f      & \text{in $\Omega$},\\
     u = 0         & \text{on $\Gamma_{j}$, $j \in \mathcal D$}, \\
     \dfrac{\partial u}{\partial \nu_{j}} + \beta_{j}\dfrac{\partial u}{\partial \tau_{j}} + \dfrac{1}{2}\cdot\dfrac{\partial \beta_{j}}{\partial \tau_{j}}u + d_{j}u =0 & \text{on $\Gamma_{j}$, $j \in \mathcal N$}.
  \end{cases}$$
\end{thm}
\emph{Proof:}
Case 1: We consider $\int_{\Omega}f = 0$.\\
The proof is similar to Grisvard\cite{GrisvardPierre2011EPiN}, Lemma 4.4.4.2.\\
Let
$$V = \{v \in H^{1}(\Omega) : v = 0 \text{ on } \Gamma_{j} \text{ for } j \in \mathcal D\},$$
$$W = \{ w \in H^{2}(\Omega) : w = 0 \text{ on } \Gamma_{j} \text{ for } j \in \mathcal D, w(S_{j}) = 0 \text{ for any } j\}.$$
Define
$$A(v,w) = \int_{\Omega}\nabla v \cdot \nabla w + \sum_{j \in \mathcal N} \int_{\Gamma_{j}}\left[\beta_{j}\dfrac{\partial v}{\partial \tau_{j}} + \dfrac{1}{2}\cdot\dfrac{\partial \beta_{j}}{\partial \tau_{j}}v + d_{j}v\right] w.$$
By Lemma 2.4.1,
$$v \in H^{\frac{1}{2}}(\Gamma_{j}), \dfrac{\partial v}{\partial \tau_{j}} \in [\tilde{H}^{\frac{1}{2}}(\Gamma_{j})]^{*},$$
$$w \in H^{\frac{3}{2}}(\Gamma_{j}) \cap H_{0}^{1}(\Gamma_{j}) \subset \tilde{H}^{\frac{1}{2}}(\Gamma_{j})$$
where $\tilde{H}^{s}([-1,1]) = \{u \in H^{s}([-1,1]) :$ the extension of $u$ by setting $u=0$ outside $[-1,1]$ is in $H^{s}(\mathbb R)\}$, $H_{0}^{1}(\Gamma_{j})$ is the closure of $C_{c}^{\infty}(\Gamma_{j})$ in $H^{1}(\Gamma_{j})$. \\
Thus $A$ is a continuous bilinear form.\\
Next we show that $A$ is coercive on $W$.\\
 For any $\phi \in C_{c}^{\infty}(\Gamma_{j})$,
 \begin{align*}
 \int_{\Gamma_{j}}\left[\beta_{j}\dfrac{\partial \phi}{\partial \tau_{j}} + \dfrac{1}{2}\cdot\dfrac{\partial \beta_{j}}{\partial \tau_{j}}\phi\right] \phi
 &= \dfrac{1}{2}\beta(S_{j+1})\phi^{2}(S_{j+1}) - \dfrac{1}{2}\beta(S_{j})\phi^{2}(S_{j}) \\
 &= 0.
 \end{align*}
Since $C_{c}^{\infty}(\Gamma_{j})$ is dense in $\tilde{H}^{\frac{1}{2}}(\Gamma_{j})$, which means $$\int_{\Gamma_{j}}\left[\beta_{j}\dfrac{\partial w}{\partial \tau_{j}} + \dfrac{1}{2}\cdot\dfrac{\partial \beta_{j}}{\partial \tau_{j}}w\right] w = 0$$
for any $w \in W$.\\
When $\mathcal D$ is not empty,
\begin{align*}
A(w,w)
&= \int_{\Omega}\nabla w \cdot \nabla w + \sum_{j \in \mathcal N} \int_{\Gamma_{j}} d_{j}w^2\\
&\geq  \int_{\Omega}\nabla w \cdot \nabla w \\
&\geq \alpha \|w\|_{H^{1}}
\end{align*}
for any $w \in W$ by the assumption that $d_{j} \geq 0$ and the Poincar\'{e}'s inequality. \\
When $\mathcal D$ is empty, we may replace $V$ by $V/C$ where $C$ is the space of constant functions in $\Omega$.\\
Then we can apply Lemma 2.4.2 to show that for any $f \in L_{p}(\Omega)$ ($\int_{\Omega}f = 0$ if $\mathcal D$ is empty), there exists a (possibly non-unique) $u \in V$ such that
$$A(u,w) = -\int_{\Omega}fw$$
for any $w \in W$.\\
Then we need to show this $u$ solves problem $(***)$.\\
Let $w \in C_{c}^{\infty}(\Omega)$, then we have
$$\triangle u = f  \text{ in $\Omega$}.$$
By the definition of $V$,
$$u\big|_{\Gamma_{j}} = 0, j \in \mathcal D.$$
Since $\triangle u = f \in L_{p}(\Omega)$, $w \in H^{2}(\Omega)$, $w(S_{j}) = 0$, by the Green's formula (Grisvard\cite{GrisvardPierre2011EPiN}, Theorem 1.5.3.11), we have
$$\sum_{j \in \mathcal N} \int_{\Gamma_{j}}\left[\beta_{j}\dfrac{\partial u}{\partial \tau_{j}} + \dfrac{1}{2}\cdot\dfrac{\partial \beta_{j}}{\partial \tau_{j}}u + d_{j} u\right] w= -\sum_{j \in \mathcal N} \int_{\Gamma_{j}}\dfrac{\partial u}{\partial \nu_{j}}w$$
for any $w \in W$. Thus we have
$$\left[\dfrac{\partial u}{\partial \nu_{j}} + \beta_{j}\dfrac{\partial u}{\partial \tau_{j}} + \dfrac{1}{2}\cdot\dfrac{\partial \beta_{j}}{\partial \tau_{j}}u + d_{j} u\right]\Big|_{\Gamma_{j}} = 0, j\in \mathcal N.$$
Case 2: We consider $\int_{\Omega}f \neq 0$ even though $\mathcal D$ is empty.\\
Let $\eta: \Omega \to [0,1]$ be a smooth cut-off function such that $\overline{\text{supp}(\eta)} \subset \Omega$ and $\int_{\Omega}\eta = 1$.\\
Let $\tilde{f} = f - (\triangle \eta)\int_{\Omega}f$. Then by the results in Case 1, there exists a $u \in H^{1}(\Omega)$ which is a solution to
$$\begin{cases}
     \triangle u = \tilde{f}      & \text{in $\Omega$},\\
     u = 0         & \text{on $\Gamma_{j}$, $j \in \mathcal D$}, \\
     \dfrac{\partial u}{\partial \nu_{j}} + \beta_{j}\dfrac{\partial u}{\partial \tau_{j}} + \dfrac{1}{2}\cdot\dfrac{\partial \beta_{j}}{\partial \tau_{j}}u + d_{j}u =0 & \text{on $\Gamma_{j}$, $j \in \mathcal N$}.
  \end{cases}$$
Then we have
$$\begin{cases}
     \triangle (u+ \eta\int_{\Omega}f) = f      & \text{in $\Omega$},\\
     (u+ \eta\int_{\Omega}f) = 0         & \text{on $\Gamma_{j}$, $j \in \mathcal D$}, \\
     \dfrac{\partial (u+ \eta\int_{\Omega}f)}{\partial \nu_{j}} + \beta_{j}\dfrac{\partial (u+ \eta\int_{\Omega}f)}{\partial \tau_{j}} \\
     + \dfrac{1}{2}\cdot\dfrac{\partial \beta_{j}}{\partial \tau_{j}}(u+ \eta\int_{\Omega}f) + d_{j}(u+ \eta\int_{\Omega}f) =0 & \text{on $\Gamma_{j}$, $j \in \mathcal N$},
  \end{cases}$$
since $\overline{\text{supp}(\eta)} \subset \Omega$.
\hfill q.e.d.\\

Next we prove that the solution $u$ in Theorem 2.4.1 is regular away from the corners with the help of the following regularity lemmas.
\begin{lemma}[Grisvard\cite{GrisvardPierre2011EPiN}, Corollary 2.5.2.2, Prop 2.5.2.4]
Let $M$ be a bounded subset in $\mathbb R^2$ with a $C^{1,1}$ boundary and $a_{0} \in C^{0}(\overline{M})$ with $a_{0} \leq 0$. We have the following results:\\
(1) If $u$ is a solution to
$$\begin{cases}
     \triangle u + a_{0}u = f \in L_{p}(M)      & \text{in $M$},\\
     u =g \in W_{p}^{2-\frac{1}{p}}(\partial M) & \text{on $\partial M$},
  \end{cases}$$
then $u \in W_{p}^{2}(M)$.\\
(2) Assume furthermore that $a_{0} < 0$, $b_{i} \in C^{1,1}(\partial M)$ and $b_{0} \in C^{0,1}(\partial M)$ such that $b_{0}b_{\nu} = b_{0}\sum_{i}b_{i}\nu_{i} \geq 0$ where $\nu_{i}$ is the $i$-th component of the outward unit normal vector of $M$ on $\partial M$.\\
If $u$ is a solution to
$$\begin{cases}
     \triangle u + a_{0}u = f \in L_{p}(M)      & \text{in $M$},\\
     \sum_{i}b_{i}u_{i} + b_{0}u =g \in W_{p}^{1-\frac{1}{p}}(\partial M) & \text{on $\partial M$},
  \end{cases}$$
then $u \in W_{p}^{2}(M)$.
\end{lemma}

In this lemma we don't need any regularity assumption on $u$, but instead we need a little more smoothness on $b_{i}$.
\begin{lemma}[Grisvard\cite{GrisvardPierre2011EPiN}, Lemma 2.4.1.4]
  Let $M$ be a bounded subset in $\mathbb R^2$ with a $C^{1,1}$ boundary and $\beta \in C^{0,1}(\partial M)$.\\
  Assume $p \leq \dfrac{2r}{2-r}$ for $r < 2$ and $u \in W_{r}^{2}(M)$ is a solution to
  $$\begin{cases}
     \triangle u = f \in L_{p}(M)      & \text{in $M$},\\
     u =g \in W_{p}^{2-\frac{1}{p}}(\partial M) & \text{on $\partial M$},
  \end{cases}$$
  or
  $$\begin{cases}
     \triangle u = f \in L_{p}(M)      & \text{in $M$},\\
     \dfrac{\partial u}{\partial \nu} + \beta \dfrac{\partial u}{\partial \tau} =g \in W_{p}^{1-\frac{1}{p}}(\partial M) & \text{on $\partial M$},
  \end{cases}$$
then $u \in W_{p}^{2}(M)$.
\end{lemma}

In this lemma we can improve the regularity of $u$ if we assume that $u \in W_{r}^{2}(M)$ and $f,g$ are smooth enough.

\begin{thm}
  Let $\Omega$ be a $C^{1,1}$ curvilinear polygonal domain, $\beta_{j}(x) \in C^{1,1}(\overline{\Gamma_{j}})$ and $d_{j} \in C^{0,1}(\overline{\Gamma_{j}})$ with $d_{j} \geq 0$ and $\dfrac{1}{2}\cdot\dfrac{\partial \beta_{j}}{\partial \tau_{j}} + d_{j} \geq 0$.\\
  Then for each $f \in L_{p}(\Omega)$, there exists a (possibly non-unique) $u \in H^{1}(\Omega)$ which is a solution to
  $$(***)\begin{cases}
     \triangle u = f      & \text{in $\Omega$},\\
     u = 0         & \text{on $\Gamma_{j}$, $j \in \mathcal D$}, \\
     \dfrac{\partial u}{\partial \nu_{j}} + \beta_{j}\dfrac{\partial u}{\partial \tau_{j}} + \dfrac{1}{2}\cdot\dfrac{\partial \beta_{j}}{\partial \tau_{j}}u + d_{j}u =0 & \text{on $\Gamma_{j}$, $j \in \mathcal N$},
  \end{cases}$$
  and $u \in W_{p}^{2}(\Omega \setminus V)$ where $V$ is any closed neighborhood of the corners.
\end{thm}
\emph{Proof:} Let $u$ be a solution to problem ($***$) in Theorem 2.4.1. We will localize the problem and use Lemma 2.4.3 and Lemma 2.4.4 to prove the regularity of $u$ away from the corners.\\
Step 1: By the classical elliptic theory, $u \in W_{p}^{2}(B_{\epsilon})$ for any ball $B_{\epsilon} \subset \Omega$.\\
Step 2: For any $x \in \Gamma_{j}$ for some $j \in \mathcal N$, there is a $\epsilon > 0$ and a $C^{2}$ open subset $M_{1} \subset \Omega$ such that $B_{2\epsilon}(x) \cap \Omega \subset M_{1}$. \\
Since $\beta_{j} \in C^{1,1}(\overline{\Gamma_{j}})$, $d_{j} \in C^{0,1}(\overline{\Gamma_{j}})$ and $\Omega$ is a $C^{2}$ curvilinear polygonal domain, we may assume that $\beta_{j} \in C^{1,1}(\overline{\Omega})$ and $d_{j} \in C^{0,1}(\overline{\Omega})$ by extension.\\
Let $\eta(x)$ be a cut-off function such that $\eta = 1$ on $B_{\epsilon}(x)$, $\eta = 0$ outside $B_{2\epsilon}(x)$.
Thus we have
$$\begin{cases}
     \triangle (\eta u) =  (\triangle \eta)u + \nabla \eta \nabla u + \eta f = \tilde{f}     & \text{in $M_{1}$},\\
     \dfrac{\partial (\eta u)}{\partial \nu} + \beta_{j}\dfrac{\partial(\eta u)}{\partial \tau} = -\eta[\dfrac{1}{2}\cdot\dfrac{\partial \beta_{j}}{\partial \tau}u + d_{j}u]+u[\dfrac{\partial \eta}{\partial \nu} + \beta_{j}\dfrac{\partial \eta}{\partial \tau}] = \tilde{g}  & \text{on $\partial M_{1}$}.
  \end{cases}$$
When $1 < p \leq 2$,
$$\nabla u \in L_{2}(\Omega) \subset L_{p}(\Omega), \tilde{f} \subset L_{p}(M_{1}),$$
$$u \in H^{1}(M_{1}) \subset W_{p}^{1}(M_{1}), u \in W_{p}^{1-\frac{1}{p}}(\partial M_{1}), \tilde{g} \in W_{p}^{1-\frac{1}{p}}(\partial M_{1}).$$
Thus by Lemma 2.4.3,
$$\eta u \in W_{p}^{2}(M_{1})$$
which means
$$u \in W_{p}^{2}(B_{\epsilon}(x) \cap \Omega).$$
When $p > 2$,
$$f \in L_{p}(\Omega) \subset L_{r}(\Omega),$$
$$u \in H^{1}(M_{1}) \subset W_{r}^{1}(M_{1}), u \in W_{r}^{1-\frac{1}{r}}(\partial M_{1}), \tilde{g} \in W_{r}^{1-\frac{1}{r}}(\partial M_{1})$$
for any $1<r<2$.\\
Repeat the previous argument, we have
$$\eta u \in W_{r}^{2}(M_{1})$$
for any $1< r < 2$. Which means
$$u \in W_{r}^{2}(B_{\epsilon}(x) \cap \Omega)$$
for any $1< r < 2$. \\
Let $\tilde{\eta}(x)$ be another cut-off function such that $\tilde{\eta} = 1$ on $B_{\frac{\epsilon}{2}}(x)$, $\eta = 0$ outside $B_{\epsilon}(x)$.
Thus we have
$$\begin{cases}
     \triangle (\tilde{\eta} u) =  (\triangle \tilde{\eta})u + \nabla \tilde{\eta} \nabla u + \tilde{\eta} f = \overline{f}     & \text{in $M_{1}$},\\
     \dfrac{\partial (\tilde{\eta} u)}{\partial \nu} + \beta_{j}\dfrac{\partial(\tilde{\eta} u)}{\partial \tau} = -\tilde{\eta}[\dfrac{1}{2}\cdot\dfrac{\partial \beta_{j}}{\partial \tau}u + d_{j}u]+u[\dfrac{\partial \tilde{\eta}}{\partial \nu} + \beta_{j}\dfrac{\partial \tilde{\eta}}{\partial \tau}] = \overline{g} & \text{on $\partial M_{1}$}.
  \end{cases}$$
Choose $r<2$ arbitrarily close to $2$, then $p \leq \dfrac{2r}{2-r}$. Thus
$$u \in W_{r}^{2}(B_{\epsilon}(x) \cap \Omega) \subset W_{p}^{1}(B_{\epsilon}(x) \cap \Omega), \overline{f} \in L_{p}(M_{1}), \overline{g} \in W_{p}^{1-\frac{1}{p}}(M_{1})$$
by the choice of $\tilde{\eta}$. \\
Now we can apply Lemma 2.4.4 to show
$$\tilde{\eta} u \in W_{p}^{2}(M_{1}).$$
Therefore $u \in W_{p}^{2}(B_{\frac{\epsilon}{2}}(x) \cap \Omega)$.\\
Step 3: For any $x \in \Gamma_{j}$ for some $j \in \mathcal D$, there is a $\epsilon > 0$ and a $C^{2}$ open subset $M_{2} \subset \Omega$ such that $B_{2\epsilon}(x) \cap \Omega \subset M_{2}$. \\
Let $\eta(x)$ be a cut-off function such that $\eta = 1$ on $B_{\epsilon}(x)$, $\eta = 0$ outside $B_{2\epsilon}(x)$.
Thus we have
$$\begin{cases}
     \triangle (\eta u) =  (\triangle \eta)u + \nabla \eta \nabla u + \eta f = \tilde{f}     & \text{in $M_{2}$},\\
     \eta u = 0 & \text{on $\partial M_{2}$}.
  \end{cases}$$
When $1 < p \leq 2$,
$$\nabla u \in L_{2}(\Omega) \subset L_{p}(\Omega), \tilde{f} \subset L_{p}(M_{2}).$$
Thus by Lemma 2.4.3,
$$\eta u \in W_{p}^{2}(M_{2})$$
which means $$u \in W_{p}^{2}(B_{\epsilon}(x) \cap \Omega).$$
When $p > 2$,
$$f \in L_{p}(\Omega) \subset L_{r}(\Omega)$$
for any $1< r < 2$. Repeat the previous argument, we have
$$\eta u \in W_{r}^{2}(M_{2})$$
for any $1< r < 2$. Which means
$$u \in W_{r}^{2}(B_{\epsilon}(x) \cap \Omega)$$
for any $1< r < 2$. \\
Let $\tilde{\eta}(x)$ be another cut-off function such that $\tilde{\eta} = 1$ on $B_{\frac{\epsilon}{2}}(x)$, $\eta = 0$ outside $B_{\epsilon}(x)$.
Thus we have
$$\begin{cases}
     \triangle (\tilde{\eta} u) =  (\triangle \tilde{\eta})u + \nabla \tilde{\eta} \nabla u + \tilde{\eta} f = \overline{f}     & \text{in $M_{2}$},\\
     \tilde{\eta} u = 0 & \text{on $\partial M_{2}$}.
  \end{cases}$$
Choose $r<2$ arbitrarily close to $2$, then $p \leq \dfrac{2r}{2-r}$. Thus
$$\nabla u \in W_{r}^{1}(B_{\epsilon}(x) \cap \Omega) \subset L_{p}(B_{\epsilon}(x) \cap \Omega), \overline{f} \in L_{p}(M_{2})$$
by the choice of $\tilde{\eta}$. Now we can apply Lemma 2.4.4 to show
$$\tilde{\eta} u \in W_{p}^{2}(M_{2}).$$
Therefore $u \in W_{p}^{2}(B_{\frac{\epsilon}{2}}(x) \cap \Omega)$.\\
By a partition of unity, we get the desired regularity results.
\hfill q.e.d.\\

With the help of Theorem 2.4.2, we can extend Theorem 2.3.1 to general boundary conditions and general curvilinear polygonal domains.
\begin{thm}
Let $\Omega$ be a $C^{1,1}$ curvilinear polygonal domain, $\beta_{j}(x) \in C^{1,1}(\overline{\Gamma_{j}})$ and $d_{j} \in C_{c}^{0,1}(\Gamma_{j})$ such that there is a small positive $\epsilon$ with $\beta_{j}(x) = \beta_{j}(S_{j})$ for $x \in B_{\epsilon}(S_{j}) \cap \overline{\Gamma_{j}}$, $\beta_{j}(x) = \beta_{j}(S_{j+1})$ for $x \in B_{\epsilon}(S_{j+1}) \cap \overline{\Gamma_{j}}$, $d_{j} \geq 0$ and $\dfrac{1}{2}\cdot\dfrac{\partial \beta_{j}}{\partial \tau_{j}} + d_{j} \geq 0$.\\
Assume $\dfrac{\phi_{j}(S_{j}) - \phi_{j+1}(S_{j}) + \frac{2\omega_{j}}{q}}{\pi}$ is not an integer for any $j$. Then for each $f \in L_{p}(\Omega)$, there exist real numbers $c_{j,m}$ and a (possibly non-unique) $u \in H^{1}(\Omega)$ such that $u$ is a solution to
$$(***)\begin{cases}
     \triangle u = f      & \text{in $\Omega$},\\
     u = 0         & \text{on $\Gamma_{j}$, $j \in \mathcal D$}, \\
     \dfrac{\partial u}{\partial \nu_{j}} + \beta_{j}\dfrac{\partial u}{\partial \tau_{j}} + \dfrac{1}{2}\cdot\dfrac{\partial \beta_{j}}{\partial \tau_{j}}u + d_{j}u =0 & \text{on $\Gamma_{j}$, $j \in \mathcal N$},
  \end{cases}$$
and
$$u - \sum_{-\frac{2}{q} < \lambda_{j,m} < 0, \lambda_{j,m} \neq -1}c_{j,m} \mathcal S_{j,m} \in W_{p}^{2}(\Omega)$$
where $\dfrac{1}{p} + \dfrac{1}{q} = 1$.
\end{thm}
\emph{Proof:}
Let $u$ be a solution to problem ($***$) in Theorem 2.4.2. Next we only need to consider the behaviour of $u$ near each corner.\\
Let $\epsilon > 0$ be small enough such that $\Omega \cap B_{8\epsilon}(S_{j})$ is a sector and $\beta_{j}, \beta_{j+1}$ are constants, $d_{j} = d_{j+1} = 0$ on $\partial \Omega \cap B_{4\epsilon}(S_{j})$. \\
Let $\eta(x)$ be a cut-off function such that $\eta = 1$ on $B_{\epsilon}(S_{j})$, $\eta = 0$ outside $B_{2\epsilon}(S_{j})$ and $\dfrac{\partial \eta}{\partial \nu_{k}} + \beta_{k}\dfrac{\partial \eta}{\partial \tau_{k}}=0$ on $\Gamma_{k}$ for $k \in \{j,j+1\}$. \\
Let $M \subset \Omega$ be a polygonal domain such that $M \cap B_{4\epsilon}(S_{j}) = \Omega \cap B_{4\epsilon}(S_{j})$.\\
Thus $\eta u$ satisfies
$$\begin{cases}
     \triangle (\eta u) = (\triangle \eta) u  + \nabla \eta \cdot \nabla u + \eta f \in L_{p}(M)  & \text{in $\Omega$},\\
     \eta u = 0         & \text{on $\partial M \setminus \Gamma_{k}$, $k \in \mathcal N \cap \{j,j+1\}$}, \\
     \dfrac{\partial (\eta u)}{\partial \nu_{k}} + \beta_{k}\dfrac{\partial (\eta u)}{\partial \tau_{k}} =0 & \text{on $\partial M \cap \Gamma_{k}$, $k \in \mathcal N \cap \{j,j+1\}$},
  \end{cases}$$
since $\eta = 1$ in $B_{\epsilon}(S_{j})$ and $u \in W_{p}^{2}(\Omega \setminus B_{\epsilon}(S_{j}))$.\\
Then by Grisvard\cite{GrisvardPierre2011EPiN}, Theorem 4.4.4.13, there exists a $v \in H^{1}(M)$ and real numbers $c_{j,m}$ such that
$$\begin{cases}
     v = \triangle (\eta u) \in L_{p}(M)  & \text{in $\Omega$},\\
     v = 0         & \text{on $\partial M \setminus \Gamma_{k}$, $k \in \mathcal N \cap \{j,j+1\}$}, \\
     \dfrac{\partial v}{\partial \nu_{k}} + \beta_{k}\dfrac{\partial v}{\partial \tau_{k}} =0 & \text{on $\partial M \cap \Gamma_{k}$, $k \in \mathcal N \cap \{j,j+1\}$},
  \end{cases}$$
and $v - \sum_{-\frac{2}{q} < \lambda_{j,m} < 0, \lambda_{j,m} \neq -1}c_{j,m} \mathcal S_{j,m} \in W_{p}^{2}(M \cap B_{2\epsilon}(S_{j})$.\\
Then $\eta u - v$ satisfies
$$\begin{cases}
     \eta u - v \in H^{1}(M), \\
     \triangle (\eta u - v) = 0  & \text{in $\Omega$},\\
     \eta u - v = 0         & \text{on $\partial M \setminus \Gamma_{k}$, $k \in \mathcal N \cap \{j,j+1\}$}, \\
     \dfrac{\partial (\eta u - v)}{\partial \nu_{k}} + \beta_{k}\dfrac{\partial (\eta u - v)}{\partial \tau_{k}} =0 & \text{on $\partial M \cap \Gamma_{k}$, $k \in \mathcal N \cap \{j,j+1\}$}.
  \end{cases}$$
By Theorem 4.4.2.4 in Grisvard\cite{GrisvardPierre2011EPiN}, $\eta u - v$ has the following expansion near $S_{j}$,
$$\eta u - v = \sum_{\lambda_{j,m} < 0, \lambda_{j,m} \neq -1}\dfrac{c_{m}}{\sqrt{\omega_{j}}}\dfrac{r_{j}^{-\lambda_{j,m}}}{\lambda_{j,m}}\cos(\lambda_{j,m}\theta_{j}+\Phi_{j+1}) + c.$$
We have
$$\eta u - v - \sum_{-\frac{2}{q} < \lambda_{j,m} < 0, \lambda_{j,m} \neq -1}c_{m} \mathcal S_{j,m} \in W_{p}^{2}(M \cap B_{2\epsilon}(S_{j})),$$
$$\eta u - \sum_{-\frac{2}{q} < \lambda_{j,m} < 0, \lambda_{j,m} \neq -1}(c_{m} + c_{j,m}) \mathcal S_{j,m} \in W_{p}^{2}(M \cap B_{2\epsilon}(S_{j})).$$
Which means
$$u - \sum_{-\frac{2}{q} < \lambda_{j,m} < 0, \lambda_{j,m} \neq -1}(c_{m} + c_{j,m}) \mathcal S_{j,m} \in W_{p}^{2}(\Omega \cap B_{\epsilon}(S_{j})).$$
\hfill q.e.d.\\

Remark: For a given $\beta_{j} \in C^{1,1}(\overline{\Gamma_{j}})$ such that there is a small positive $\epsilon$ with $\beta_{j}(x) = \beta_{j}(S_{j})$ for $x \in B_{\epsilon}(S_{j}) \cap \overline{\Gamma_{j}}$, $\beta_{j}(x) = \beta_{j}(S_{j+1})$ for $x \in B_{\epsilon}(S_{j+1}) \cap \overline{\Gamma_{j}}$, such $d_{j}$ always exists since we can choose $d_{j} \geq 0$ such that the support of $d_{j}$ contains the support of $\dfrac{\partial \beta_{j}}{\partial \tau_{j}}$ and $d_{j} \geq -\dfrac{1}{2}\cdot\dfrac{\partial \beta_{j}}{\partial \tau_{j}}$.

With the following classical regularity lemma and the same method as Theorem 2.3.1 yields Theorem 2.3.2, we can prove an existence theorem similar to Theorem 2.3.2.
\begin{lemma}[Grisvard\cite{GrisvardPierre2011EPiN}, Theorem 2.5.1.1]
  Let $M$ be a bounded subset in $\mathbb R^2$ with a $C^{k+1,1}$ boundary. Let $b_{i} \in C^{k,1}(\overline{M})$ such that $\sum_{i}b_{i}\nu_{i} \neq 0$ where $\nu_{i}$ is the $i$-th component of the outward unit normal vector of $\partial M$. If $u \in W_{p}^{2}(M)$ is a solution to
$$\begin{cases}
     \triangle u = f \in W_{p}^{k}(M)      & \text{in $M$},\\
     \sum_{i}b_{i}u_{i} =g \in W_{p}^{k+1-\frac{1}{p}}(\partial M) & \text{on $\partial M$},
  \end{cases}$$
or
$$\begin{cases}
     \triangle u = f \in W_{p}^{k}(M)      & \text{in $M$},\\
     u =g \in W_{p}^{k+2-\frac{1}{p}}(\partial M) & \text{on $\partial M$},
  \end{cases}$$
then $u \in W_{p}^{k+2}(M)$.
\end{lemma}

\begin{thm}
Let $k \geq 0$ be an integer and $\Omega$ be a $C^{k+1,1}$ curvilinear polygonal domain, $\beta_{j}(x) \in C^{k,1}(\overline{\Gamma_{j}}) \cap C^{1,1}(\overline{\Gamma_{j}})$, $d_{j} \in C_{c}^{k,1}(\Gamma_{j})$ such that there is a small positive $\epsilon$ with $\beta_{j}(x) = \beta_{j}(S_{j})$ for $x \in B_{\epsilon}(S_{j}) \cap \overline{\Gamma_{j}}$, $\beta_{j}(x) = \beta_{j}(S_{j+1})$ for $x \in B_{\epsilon}(S_{j+1}) \cap \overline{\Gamma_{j}}$, $d_{j} \geq 0$ and $\dfrac{1}{2}\cdot\dfrac{\partial \beta_{j}}{\partial \tau_{j}} + d_{j} \geq 0$.\\
Assume $\dfrac{\phi_{j}(S_{j}) - \phi_{j+1}(S_{j}) + \frac{2\omega_{j}}{q}}{\pi}$ is not an integer for any $j$. Then for each $f \in W_{p}^{k}(\Omega)$, $g_{j} \in W_{p}^{k+2-\frac{1}{p}}(\Gamma_{j})$ if $j \in \mathcal D$, $g_{j} \in W_{p}^{k+1-\frac{1}{p}}(\Gamma_{j})$ if $j \in \mathcal N$ and $g_{j}(S_{j}) = g_{j+1}(S_{j})$ if $j,j+1 \in \mathcal D$, there exist real numbers $c_{j,m}$ and a (possibly non-unique) $u \in H^{1}(\Omega)$ such that $u$ is a solution to
$$(***)\begin{cases}
     \triangle u = f      & \text{in $\Omega$},\\
     u = g_{j}         & \text{on $\Gamma_{j}$, $j \in \mathcal D$}, \\
     \dfrac{\partial u}{\partial \nu_{j}} + \beta_{j}\dfrac{\partial u}{\partial \tau_{j}} + \dfrac{1}{2}\cdot\dfrac{\partial \beta_{j}}{\partial \tau_{j}}u + d_{j}u =g_{j} & \text{on $\Gamma_{j}$, $j \in \mathcal N$},
  \end{cases}$$
and
$$u - \sum_{-k-\frac{2}{q} < \lambda_{j,m} < 0}c_{j,m} \mathcal S_{j,m} \in W_{p}^{k+2}(\Omega)$$
where $\dfrac{1}{p} + \dfrac{1}{q} = 1$.
\end{thm}

Next we derive the existence theorems in the H\"{o}lder spaces with the help of Theorem 2.4.4 and the following lemma which connects the results in the Sobolev spaces to the H\"{o}lder spaces.
\begin{lemma}[Grisvard\cite{GrisvardPierre2011EPiN}, Theorem 6.3.2.1]
  Let $M$ be a bounded open domain in $\mathbb R^{2}$ with a $C^{2,1}$ boundary and $b \in C^{1,1}(\partial M)$.\\
  Then for any $0<\sigma<1$, if $u \in W_{p}^{2}(M)$ with $p > 2$ is a solution to
  $$\begin{cases}
     \triangle u = f \in C^{0,\sigma}(\overline{M})      & \text{in $M$},\\
     u = g \in C^{2,\sigma}(\partial M)         & \text{on $\partial M$},
  \end{cases}$$
  or
  $$\begin{cases}
     \triangle u = f  \in C^{0,\sigma}(\overline{M})    & \text{in $M$},\\
     \dfrac{\partial u}{\partial \nu} + b\dfrac{\partial u}{\partial \tau} =g \in C^{1,\sigma}(\partial M) & \text{on $\partial M$},
  \end{cases}$$
then $u \in C^{2,\sigma}(\overline{M})$.
\end{lemma}

\begin{thm}
Let $\Omega$ be a $C^{2,1}$ curvilinear polygonal domain, $\beta_{j}(x) \in C^{2,1}(\overline{\Gamma_{j}})$, $d_{j} \in C_{c}^{1,1}(\Gamma_{j})$ such that there is a small positive $\epsilon$ with $\beta_{j}(x) = \beta_{j}(S_{j})$ for $x \in B_{\epsilon}(S_{j}) \cap \overline{\Gamma_{j}}$, $\beta_{j}(x) = \beta_{j}(S_{j+1})$ for $x \in B_{\epsilon}(S_{j+1}) \cap \overline{\Gamma_{j}}$, $d_{j} \geq 0$ and $\dfrac{1}{2}\cdot\dfrac{\partial \beta_{j}}{\partial \tau_{j}} + d_{j} \geq 0$.\\
Assume $0< \sigma <1$ and $\dfrac{\phi_{j}(S_{j}) - \phi_{j+1}(S_{j}) + (2+\sigma)\omega_{j}}{\pi}$ is not an integer for any $j$. Then for each $f \in C^{0,\sigma}(\overline{\Omega})$, there exist real numbers $c_{j,m}$ and a function $u$ such that $u$ is a solution to
$$(***)\begin{cases}
     \triangle u = f      & \text{in $\Omega$},\\
     u = 0         & \text{on $\Gamma_{j}$, $j \in \mathcal D$}, \\
     \dfrac{\partial u}{\partial \nu_{j}} + \beta_{j}\dfrac{\partial u}{\partial \tau_{j}} + \dfrac{1}{2}\cdot\dfrac{\partial \beta_{j}}{\partial \tau_{j}}u + d_{j}u = 0 & \text{on $\Gamma_{j}$, $j \in \mathcal N$},
  \end{cases}$$
and
$$u - \sum_{-(2+\sigma) < \lambda_{j,m} < 0, \lambda_{j,m} \neq -1}c_{j,m} \mathcal S_{j,m} \in C^{2,\sigma}(\overline{\Omega}).$$
\end{thm}
\emph{Proof:}
Let $p$ be large enough such that $\dfrac{\phi_{j}(S_{j}) - \phi_{j+1}(S_{j}) + \frac{2\omega_{j}}{q}}{\pi}$ is not an integer where $\dfrac{1}{p} + \dfrac{1}{q} = 1$ and $W_{p}^{2}(\Omega) \subset C^{1,\sigma}(\overline{\Omega})$. \\
Then we can apply Theorem 2.4.3 to find a $u \in H^{1}(\Omega)$ which solves problem ($***$) and
$$u - \sum_{-\frac{2}{q} < \lambda_{j,m} < 0, \lambda_{j,m} \neq -1}c_{j,m} \mathcal S_{j,m} \in W_{p}^{2}(\Omega).$$
We prove the result by localizing the problem. \\
Let $\eta: [0,\infty) \to [0,1]$ be a smooth cut-off function such that $\eta = 1$ on $[0,1]$ and $\eta = 0$ on $[2,\infty)$. Define $\eta_{\epsilon,y}(x) = \eta(\frac{|x-y|}{\epsilon})$ for any $\epsilon > 0$ and $x,y \in \mathbb R^{2}$.\\
Step 1: For any $y \in \Omega$, there is a $\epsilon > 0$ such that $\overline{B_{4\epsilon}(y)} \subset \Omega$. Then $\eta_{\epsilon,y}u \in W_{p}^{2}(B_{4\epsilon}(y))$ satisfies
$$\begin{cases}
     \triangle (\eta_{\epsilon,y} u) = \triangle (\eta_{\epsilon,y})u + \nabla \eta_{\epsilon,y} \cdot \nabla u + \eta_{\epsilon,y}f \in C^{0,\sigma}(\overline{B_{2\epsilon}(y)})      & \text{in $B_{4\epsilon}(y)$},\\
     \eta_{\epsilon,y}u = 0          & \text{on $\partial B_{4\epsilon}(y)$},
  \end{cases}$$
since $u \in W_{p}^{2}(B_{4\epsilon}(y)) \subset C^{1,\sigma}(\overline{B_{4\epsilon}(y)})$. \\
Thus by Lemma 2.4.6,
$$\eta_{\epsilon,y} u \in C^{2,\sigma}(\overline{B_{4\epsilon}(y)})$$
which means
$$u \in C^{2,\sigma}(\overline{B_{\epsilon}(y)}).$$
Step 2: For any $y \in \Gamma_{j}$ for some $j \in \mathcal D$, there is a $\epsilon > 0$ and a $C^{2,1}$ open subset $M_{1} \subset \Omega$ such that $B_{2\epsilon}(x) \cap \Omega \subset M_{1}$. \\
Similar to Step 1, we have
$$u \in C^{2,\sigma}(\overline{B_{\epsilon}(y) \cap \Omega}).$$
Step 3: For any $y \in \Gamma_{j}$ for some $j \in \mathcal N$, there is a $\epsilon > 0$ and a $C^{2,1}$ open subset $M_{2} \subset \Omega$ such that $B_{2\epsilon}(x) \cap \Omega \subset M_{2}$. By extending $\beta_{j}, d_{j}$ to $C^{2,1}(\overline{\Omega}), C^{1,1}(\overline{\Omega})$, we have
$$\begin{cases}
     \triangle (\eta_{\epsilon,y} u) = \triangle (\eta_{\epsilon,y})u + \nabla \eta_{\epsilon,y} \cdot \nabla u + \eta_{\epsilon,y}f \in C^{0,\sigma}(\overline{M_{2}})      & \text{in $M_{2}$},\\
     \dfrac{\partial (\eta_{\epsilon,y} u)}{\partial \nu} + \beta_{j}\dfrac{\partial (\eta_{\epsilon,y} u)}{\partial \tau} \\
     =  -\eta_{\epsilon,y}[\dfrac{1}{2}\cdot\dfrac{\partial \beta_{j}}{\partial \tau}u + d_{j}u]+u[\dfrac{\partial \eta_{\epsilon,y}}{\partial \nu} + \beta_{j}\dfrac{\partial \eta_{\epsilon,y}}{\partial \tau}]  \in C^{1,\sigma}(\partial M_{2})        & \text{on $\partial M_{2}$},
  \end{cases}$$
since $u \in W_{p}^{2}(M_{2}) \subset C^{1,\sigma}(\overline{M_{2}})$ and $\beta_{j} \in C^{2,1}(\overline{M_{2}})$. \\
Thus by Lemma 2.4.6,
$$\eta_{\epsilon,y} u \in C^{2,\sigma}(\overline{M_{2}})$$
which means
$$u \in C^{2,\sigma}(\overline{B_{\epsilon}(y) \cap \Omega}).$$
Step 4: Let $\epsilon > 0$ be small enough such that $\Omega \cap B_{8\epsilon}(S_{j})$ is a sector and $\beta_{j}, \beta_{j+1}$ are constants, $d_{j} = d_{j+1} = 0$ on $\partial \Omega \cap B_{4\epsilon}(S_{j})$. \\
Let $\eta(x)$ be a cut-off function such that $\eta = 1$ on $B_{\epsilon}(S_{j})$, $\eta = 0$ outside $B_{2\epsilon}(S_{j})$ and $\dfrac{\partial \eta}{\partial \nu_{k}} + \beta_{k}\dfrac{\partial \eta}{\partial \tau_{k}}=0$ on $\Gamma_{k}$ for $k \in \{j,j+1\}$. \\
Let $M \subset \Omega$ be a polygonal domain with at least four edges such that $M \cap B_{4\epsilon}(S_{j}) = \Omega \cap B_{4\epsilon}(S_{j})$.\\
Thus $\eta u$ satisfies
$$(***)\begin{cases}
     \triangle (\eta u) = (\triangle \eta) u  + \nabla \eta \cdot \nabla u + \eta f \in L_{p}(M)  & \text{in $M$},\\
     \eta u = 0         & \text{on $\partial M \setminus \Gamma_{k}$, $k \in \mathcal N \cap \{j,j+1\}$}, \\
     \dfrac{\partial (\eta u)}{\partial \nu_{k}} + \beta_{k}\dfrac{\partial (\eta u)}{\partial \tau_{k}} =0 & \text{on $\partial M \cap \Gamma_{k}$, $k \in \mathcal N \cap \{j,j+1\}$},
  \end{cases}$$
since $\eta = 1$ in $B_{\epsilon}(S_{j})$ and $u \in W_{p}^{2}(\Omega \setminus B_{\epsilon}(S_{j}))$. \\
Since $M$ has at least four edges and $\eta u = 0$ on at least two edges of $M$, problem $(***)$ has at most one solution in $W_{p}^{2}(M)$. \\
This is the case Grisvard considers in \cite{GrisvardPierre2011EPiN} Theorem 6.4.2.5. Thus there exist real numbers $c_{j,m}$ such that
$$\eta u - \sum_{-(2+\sigma) < \lambda_{j,m} < 0,\lambda_{j,m} \neq -1}c_{j,m} \mathcal S_{j,m} \in C^{2,\sigma}(\overline{\Omega \cap B_{\epsilon}(S_{j})}).$$
Which means
$$u - \sum_{-(2+\sigma) < \lambda_{j,m} < 0,\lambda_{j,m} \neq -1}c_{j,m} \mathcal S_{j,m} \in C^{2,\sigma}(\overline{\Omega \cap B_{\epsilon}(S_{j})}).$$
\hfill q.e.d.\\

As Theorem 2.3.3 yields Theorem 2.3.4, we have the following extension of Theorem 2.4.5.
\begin{thm}
Let $k \geq 0$ be an integer and $\Omega$ be a $C^{k+2,1}$ curvilinear polygonal domain, $\beta_{j}(x) \in C^{k+2,1}(\overline{\Gamma_{j}})$, $d_{j} \in C_{c}^{k+1,1}(\Gamma_{j})$ such that there is a small positive $\epsilon$ with $\beta_{j}(x) = \beta_{j}(S_{j})$ for $x \in B_{\epsilon}(S_{j}) \cap \overline{\Gamma_{j}}$, $\beta_{j}(x) = \beta_{j}(S_{j+1})$ for $x \in B_{\epsilon}(S_{j+1}) \cap \overline{\Gamma_{j}}$, $d_{j} \geq 0$ and $\dfrac{1}{2}\cdot\dfrac{\partial \beta_{j}}{\partial \tau_{j}} + d_{j} \geq 0$.\\
Assume that $0< \sigma <1$ and $\dfrac{\phi_{j}(S_{j}) - \phi_{j+1}(S_{j}) + (k+2+\sigma)\omega_{j}}{\pi}$ is not an integer for any $j$.
Then for each $f \in C^{k,\sigma}(\overline{\Omega})$, $g_{j} \in C^{k+2,\sigma}(\overline{\Gamma_{j}})$ if $j \in \mathcal D$, $g_{j} \in C^{k+1,\sigma}(\overline{\Gamma_{j}})$ if $j \in \mathcal N$ and $g_{j}(S_{j}) = g_{j+1}(S_{j})$ if $j,j+1 \in \mathcal D$, there exist real numbers $c_{j,m}$ and a function $u$ such that $u$ is a solution to
$$\begin{cases}
     \triangle u = f      & \text{in $\Omega$},\\
     u = g_{j}         & \text{on $\Gamma_{j}$, $j \in \mathcal D$}, \\
     \dfrac{\partial u}{\partial \nu_{j}} + \beta_{j}\dfrac{\partial u}{\partial \tau_{j}} + \dfrac{1}{2}\cdot\dfrac{\partial \beta_{j}}{\partial \tau_{j}}u + d_{j}u = g_{j} & \text{on $\Gamma_{j}$, $j \in \mathcal N$},
  \end{cases}$$
and
$$u - \sum_{-(k+2+\sigma) < \lambda_{j,m} < 0}c_{j,m} \mathcal S_{j,m} \in C^{k+2,\sigma}(\overline{\Omega}).$$
\end{thm}

\subsection{Apriori estimates}
In this section, we extend Grisvard's apriori estimates to general elliptic operators with general boundary conditions.

We first state Grisvard's apriori estimates.
\begin{thm}[Grisvard\cite{GrisvardPierre2011EPiN}, Theorem 6.4.2.4.]
Let $\Omega$ be a polygonal domain and $\beta_{j}$ be real numbers.\\
Assume $\dfrac{\phi_{j} - \phi_{j+1} + (2+\sigma)\omega_{j}}{\pi}$ is not an integer for any $j$. Then there is a constant $C$ such that for each $u \in C^{2,\sigma}(\overline{\Omega})$ with the boundary conditions \\
  $$\begin{cases}
     u = 0                                & \text{on $\Gamma_{j}$, $j \in \mathcal D$}, \\
     \dfrac{\partial u}{\partial \nu_{j}} + \beta_{j}\dfrac{\partial u}{\partial \tau_{j}}=0 & \text{on $\Gamma_{j}$, $j \in \mathcal N$},
  \end{cases}$$
  we have
  $$\|u\|_{2,\sigma, \overline{\Omega}} \leq C\{\|\triangle u\|_{\sigma, \overline{\Omega}} + \|u\|_{1,\sigma, \overline{\Omega}}     \}.$$
\end{thm}
\subsubsection{$C^{0}$ estimates}
\begin{thm}
Let $\Omega$ be a $C^{2,1}$ curvilinear polygonal domain.
Let $J(u) = a_{ij}D_{ij}u + b_{i}D_{i}u + cu$ be a strongly elliptic operator where $a_{ij} = a_{ji}, b_{i}, c \in C^{0,\sigma}(\overline{\Omega})$, $c \leq -1$ and there is a positive $\alpha$ such that $\sum_{i,j=1}^{2}a_{ij}(x)\xi_{i}\xi_{j} \geq \alpha|\xi|^2$ for all $x \in \overline{\Omega}$ and $\xi \in \mathbb R^2$. Let $\beta_{j}, d_{j} \in C^{1,\sigma}(\overline{\Gamma_{j}})$ such that $d_{j} \geq 1$. Then we have
$$\|u\|_{0,\overline{\Omega}} \leq \|Ju\|_{0,\overline{\Omega}} + \sum_{j \in \mathcal D}\|u\|_{0,\overline{\Gamma}_{j}} +  \sum_{j \in \mathcal N}\|\dfrac{\partial u}{\partial \nu_{j}} + \beta_{j}\dfrac{\partial u}{\partial \tau_{j}} + d_{j}u \|_{0,\overline{\Gamma}_{j}}    $$
for all $u \in C^{2,\sigma}(\overline{\Omega})$ such that $u(S_{j}) = 0$ if $j,j+1 \in \mathcal N$.
\end{thm}
\emph{Proof:}
Let $|u(x_{0})| = \max_{x \in \overline{\Omega}} |u|$.\\
Case 1: $u(x_{0}) > 0$. \\
If $x_{0} \in \Omega$, then $a_{ij}D_{ij}u(x_{0}) \leq 0$. By the assumption that $c \leq -1$, we have
$$(Lu)(x_{0}) \leq c(x_{0})u(x_{0}) \leq -u(x_{0}),$$
which is
$$u(x_{0}) \leq |(Ju)(x_{0})|.$$
If $x_{0} \in \Gamma_{j}$ for some $j \in \mathcal D$, then
$$u(x_{0}) \leq \sum_{j \in \mathcal D}\|u\|_{0,\overline{\Gamma}_{j}}.$$
If $x_{0} \in \Gamma_{j}$ for some $j \in \mathcal N$, then $\dfrac{\partial u}{\partial \nu_{j}} + \beta_{j}\dfrac{\partial u}{\partial \tau_{j}} \geq 0$ at $x_{0}$. Thus
$$\dfrac{\partial u}{\partial \nu_{j}}(x_{0}) + \beta_{j}\dfrac{\partial u}{\partial \tau_{j}}(x_{0}) + d_{j}u(x_{0}) \geq 0 + d_{j}(x_{0})u(x_{0})\geq u(x_{0})$$
which is
$$u(x_{0}) \leq \sum_{j \in \mathcal N}\|\dfrac{\partial u}{\partial \nu_{j}} + \beta_{j}\dfrac{\partial u}{\partial \tau_{j}} + d_{j}u \|_{0,\overline{\Gamma}_{j}}.$$
If $x_{0} = S_{j}$ such that one of $j,j+1 \in \mathcal D$, then
$$|u(x_{0})| \leq \sum_{j \in \mathcal D}\|u\|_{0,\overline{\Gamma}_{j}}.$$
If $x_{0} = S_{j}$ such that both $j,j+1 \in \mathcal N$, then by the assumption $$u(S_{j}) = 0.$$
Case 2: $u(x_{0}) < 0$. The proof is similar to Case 1.
\hfill q.e.d.\\

\begin{corollary}
Let $\Omega$ be a $C^{2,1}$ curvilinear polygonal domain.
  Let $J(u) = a_{ij}D_{ij}u + b_{i}D_{i}u + cu$ be a strongly elliptic operator where $a_{ij} = a_{ji}, b_{i}, c \in C^{0,\sigma}(\overline{\Omega})$ and there is a positive $\alpha$ such that $\sum_{i,j=1}^{2}a_{ij}(x)\xi_{i}\xi_{j} \geq \alpha|\xi|^2$ for all $x \in \overline{\Omega}$ and $\xi \in \mathbb R^2$. Let $\beta_{j}, d_{j} \in C^{1,\sigma}(\overline{\Gamma_{j}})$. \\
Then for any $\mu \leq -\text{sup}|c|-1$, $d \geq \text{sup}|d_{j}|+1$,
$$\|u\|_{0,\overline{\Omega}} \leq \|Ju+\mu u\|_{0,\overline{\Omega}} + \sum_{j \in \mathcal D}\|u\|_{0,\overline{\Gamma}_{j}} +  \sum_{j \in \mathcal N}\|\dfrac{\partial u}{\partial \nu_{j}} + \beta_{j}\dfrac{\partial u}{\partial \tau_{j}} + d_{j}u + d u \|_{0,\overline{\Gamma}_{j}}     $$
for all $u \in C^{2,\sigma}(\overline{\Omega})$ such that $u(S_{j}) = 0$ if $j,j+1 \in \mathcal N$.
\end{corollary}

\subsubsection{$C^{2,\sigma}$ estimates}

Based on Theorem 2.5.1 and the inverse trace theorems we derived in Section 2.8, we prove apriori estimates for the Laplace operator with non-homogeneous boundary value conditions.
\begin{thm}
Let $\Omega$ be a $C^{2,1}$ curvilinear polygonal domain, $\beta_{j}(x) \in C^{1,\sigma}(\overline{\Gamma_{j}})$ such that there is a small positive $\epsilon$ with $\beta_{j}(x) = \beta_{j}(S_{j})$ for $x \in B_{\epsilon}(S_{j}) \cap \overline{\Gamma_{j}}$, $\beta_{j}(x) = \beta_{j}(S_{j+1})$ for $x \in B_{\epsilon}(S_{j+1}) \cap \overline{\Gamma_{j}}$.\\
Assume $0 < \sigma < 1$ and $\dfrac{\phi_{j}(S_{j}) - \phi_{j+1}(S_{j}) + (2+\sigma)\omega_{j}}{\pi}$ is not an integer for any $j$. Then there is a constant $C = C(\Omega,\beta_{j})$ such that for each $u \in C^{2,\sigma}(\overline{\Omega})$, we have
  $$\|u\|_{2,\sigma, \overline{\Omega}} \leq C\{\|\triangle u\|_{\sigma, \overline{\Omega}} + \|u\|_{1,\sigma, \overline{\Omega}} + \sum_{j \in \mathcal D}\|u\|_{2,\sigma,\overline{\Gamma_{j}}} +  \sum_{j \in \mathcal N}\|\dfrac{\partial u}{\partial \nu_{j}} + \beta_{j}\dfrac{\partial u}{\partial \tau_{j}}\|_{1,\sigma,\overline{\Gamma_{j}}}    \}.$$
\end{thm}
\emph{Proof:} The proof relies on Theorem 2.5.1 and the inverse trace theorems we derived in Section 2.8. \\
By the classical elliptic theory, for any $x \in \Gamma_{j}$, there is a $\epsilon > 0$ and a constant $C$ such that
\begin{align*}
\|u\|_{2,\sigma, \overline{\Omega \cap B_{\epsilon}(x)}}
&\leq C\{\|\triangle u\|_{\sigma, \overline{\Omega \cap B_{2\epsilon}(x)}} + \|u\|_{1,\sigma, \overline{\Omega \cap B_{2\epsilon}(x)}}\\
&\quad + \sum_{j \in \mathcal D}\|u\|_{2,\sigma,\overline{\Gamma_{j} \cap B_{2\epsilon}(x)}} +  \sum_{j \in \mathcal N}\|\dfrac{\partial u}{\partial \nu_{j}} + \beta_{j}\dfrac{\partial u}{\partial \tau_{j}}\|_{1,\sigma,\overline{\Gamma_{j} \cap B_{2\epsilon}(x)}}    \}\\
&\leq C\{\|\triangle u\|_{\sigma, \overline{\Omega}} + \|u\|_{1,\sigma, \overline{\Omega}} + \sum_{j \in \mathcal D}\|u\|_{2,\sigma,\overline{\Gamma_{j}}}\\
 &\quad +  \sum_{j \in \mathcal N}\|\dfrac{\partial u}{\partial \nu_{j}} + \beta_{j}\dfrac{\partial u}{\partial \tau_{j}}\|_{1,\sigma,\overline{\Gamma_{j}}}    \}.
\end{align*}
For any $x \in \Omega$, there is a $\epsilon > 0$ and a constant $C$ such that
\begin{align*}
\|u\|_{2,\sigma, \overline{\Omega \cap B_{\epsilon}(x)}}
&\leq C\{\|\triangle u\|_{\sigma, \overline{\Omega \cap B_{2\epsilon}(x)}} + \|u\|_{1,\sigma, \overline{\Omega \cap B_{2\epsilon}(x)}} \}\\
&\leq C\{\|\triangle u\|_{\sigma, \overline{\Omega}} + \|u\|_{1,\sigma, \overline{\Omega}} \}.
\end{align*}
Thus we only need to deal with the behaviour of $u$ near the corners.\\
Let $\eta_{j}$ be smooth cut-off functions such that:\\
(1) $\eta_{j} = 1$ near $S_{j}$, $\text{supp}(\eta_{j}) \cap \Gamma_{l} = \emptyset$ for any $l \notin \{j,j+1\}$.\\
(2) $\dfrac{\partial \eta_{j}}{\partial \nu_{k}} + \beta_{k}\dfrac{\partial \eta_{j}}{\partial \tau_{k}}=0$ on $\Gamma_{k}$ for $k \in \{j,j+1\}$. This is possible since $\beta_{j}$ is a constant near each corner and $\eta_{j} = 1$ near $S_{j}$.\\
Let $M \subset \Omega$ be a polygonal domain with exactly four corners and one of them is $S_{j}$, two edges of $\partial M$ lie on $\Gamma_{j}$ and $\Gamma_{j+1}$. We consider the Dirichlet boundary conditions on the other two edges and we may choose the other two edges of $M$ to satisfy that $\dfrac{\phi_{j} - \phi_{j+1} + (2+\sigma)\omega_{j}}{\pi}$ is not an integer on $M$.\\
Let $w \in C^{2,\sigma}(\overline{M})$ be the function we construct in Section 2.8 with respect to the boundary values of $\eta_{j}u$ on $M$, then by Theorem 2.5.1,
$$\|\eta_{j}u-w\|_{2,\sigma, \overline{M}} \leq C\{\|\triangle (\eta_{j}u-w)\|_{\sigma, \overline{M}} + \|\eta_{j}u-w\|_{1,\sigma, \overline{M}}     \}.$$
Thus we have
\begin{align*}
  \|\eta_{j}u\|_{2,\sigma, \overline{M}} &\leq C\{\|w\|_{2,\sigma, \overline{M}} + \|\triangle u\|_{\sigma, \overline{\Omega}} + \|u\|_{1,\sigma, \overline{\Omega}}     \} \\
  &\leq C\{\|\triangle u\|_{\sigma, \overline{\Omega}} + \|u\|_{1,\sigma, \overline{\Omega}} + \sum_{j \in \mathcal D}\|u\|_{2,\sigma,\overline{\Gamma_{j}}}\\
  &\quad +  \sum_{j \in \mathcal N}\|\dfrac{\partial u}{\partial \nu_{j}} + \beta_{j}\dfrac{\partial u}{\partial \tau_{j}}\|_{1,\sigma,\overline{\Gamma_{j}}}    \}.
\end{align*}
Thus by a partition of unity, we have the desired estimate.
\hfill q.e.d.\\

Next we extend Theorem 2.5.3 to the setting that $\beta_{j}$ is not a constant near the corner by the method of freezing coefficients.
\begin{thm}
Let $\Omega$ be a $C^{2,1}$ curvilinear polygonal domain, $\beta_{j}(x) \in C^{1,\sigma}(\overline{\Gamma_{j}})$.\\
Assume $0 < \sigma < 1$ and $\dfrac{\phi_{j}(S_{j}) - \phi_{j+1}(S_{j}) + (2+\sigma)\omega_{j}}{\pi}$ is not an integer for any $j$. Then there is a constant $C = C(\Omega,\beta_{j})$ such that for each $u \in C^{2,\sigma}(\overline{\Omega})$, we have
  $$\|u\|_{2,\sigma, \overline{\Omega}} \leq C\{\|\triangle u\|_{\sigma, \overline{\Omega}} + \|u\|_{1,\sigma, \overline{\Omega}} + \sum_{j \in \mathcal D}\|u\|_{2,\sigma,\overline{\Gamma_{j}}} +  \sum_{j \in \mathcal N}\|\dfrac{\partial u}{\partial \nu_{j}} + \beta_{j}\dfrac{\partial u}{\partial \tau_{j}}\|_{1,\sigma,\overline{\Gamma_{j}}}    \}.$$
\end{thm}
\emph{Proof:}
The proof depends on Theorem 2.5.3 and the method of freezing coefficients on $\beta_{j}$ near each corner.
\hfill q.e.d.\\

Next we derive apriori estimates for general elliptic operators based on the previous theorems.

Let $J(u) = a_{ij}D_{ij}u + b_{i}D_{i}u + cu$ be a strongly elliptic operator. Let $A_{j}$ be the operator obtained from $J$ by freezing the coefficients of its principle part at $S_{j}$,
  $$A_{j}(f) = \sum_{k,l=1}^{2}a_{kl}(S_{j})D_{kl}f.$$
  And we consider the matrix $\mathcal B_{j}$ such that
  $$\mathcal B_{j} \mathcal A_{j} \mathcal B_{j} =  I$$
  where $\mathcal A_{j}$ is the symmetric matrix whose entries are $a_{kl}(S_{j})$.

  By direct computations, we have
  $$(A_{j}u)(\mathcal B_{j}^{-1}x) = (\triangle (u(\mathcal B_{j}^{-1}x)))(x).$$

We consider a linear coordinate change $\tilde{x} = \mathcal \mathcal B_{j}x$. $\tilde{\Omega} = \mathcal B_{j}(\Omega)$ will still be a curvilinear polygonal domain. We will use $\tilde{\Omega}, \tilde{\phi_{j}}$ and $\tilde{\omega_{j}}$ to denote the corresponding notations.
\begin{thm}
Let $\Omega$ be a $C^{2,1}$ curvilinear polygonal domain, $\beta_{j}(x) \in C^{1,\sigma}(\overline{\Gamma_{j}})$.
Let $J(u) = a_{ij}D_{ij}u + b_{i}D_{i}u + cu$ be a strongly elliptic operator where $a_{ij} = a_{ji}, b_{i}, c \in C^{0,\sigma}(\overline{\Omega})$ and there is a positive $\alpha$ such that $\sum_{i,j=1}^{2}a_{ij}(x)\xi_{i}\xi_{j} \geq \alpha|\xi|^2$ for all $x \in \overline{\Omega}$ and $\xi \in \mathbb R^2$. \\
Assume $0 < \sigma < 1$ and $\dfrac{\tilde{\phi}_{j}(S_{j}) - \tilde{\phi}_{j+1}(S_{j}) + (2+\sigma)\tilde{\omega}_{j}}{\pi}$ is not an integer for any $j$. Then there is a constant $C$ such that for each $u \in C^{2,\sigma}(\overline{\Omega})$, we have
  $$\|u\|_{2,\sigma, \overline{\Omega}} \leq C\{\|J(u)\|_{\sigma, \overline{\Omega}} + \|u\|_{1,\sigma, \overline{\Omega}} + \sum_{j \in \mathcal D}\|u\|_{2,\sigma,\overline{\Gamma_{j}}} +  \sum_{j \in \mathcal N}\|\dfrac{\partial u}{\partial \nu_{j}} + \beta_{j}\dfrac{\partial u}{\partial \tau_{j}}\|_{1,\sigma,\overline{\Gamma_{j}}}    \}.$$
\end{thm}
\emph{Proof:}
The proof depends on the previous estimates and the method of freezing coefficients on $a_{ij}$ near each corner. \\
Step 1: We consider a partition of unity $\{\eta_{j}\}_{j=0,1,...,I_{0}}$ on $\overline{\Omega}$ with $\eta_{j} \in C_{c}^{\infty}(\mathbb R^2)$ such that for $j=1,2,...,I_{0}$, $\eta_{j} = 1$ near $S_{j}$ and $\text{supp}(\eta_{j}) \cap \Gamma_{k} = \emptyset$ if $k \notin \{j,j+1\}$.\\
Thus for $j=0$, $\text{supp}(\eta_{0}u) \cap \overline{\Omega}$ is a smooth domain in $\mathbb R^2$. By the classical elliptic theory, there is a $C_{1}$ such that
$$\|\eta_{0}u\|_{2,\sigma, \text{supp}(\eta_{0}u) \cap \overline{\Omega}} \leq C_{1}\{\|J(\eta_{0}u)\|_{\sigma, \text{supp}(\eta_{0}u) \cap \overline{\Omega}} + \|\eta_{0}u\|_{1,\sigma, \text{supp}(\eta_{0}u) \cap \overline{\Omega}}     \}$$
which is
$$\|\eta_{0}u\|_{2,\sigma, \overline{\Omega}} \leq C_{2}\{\|J(u)\|_{\sigma, \overline{\Omega}} + \|u\|_{1,\sigma, \overline{\Omega}}     \}.$$
Step 2: For $j=1,2,...,I_{0}$, let $v_{j} = \eta_{j}u$. Then we have
$$\triangle v_{j} + \sum_{k,l=1}^{2}b_{kl}D_{kl}v_{j} + b_{i}D_{i}v_{j} + cv_{j} = g$$
where $g = (L(\eta_{j}u)) \circ \mathcal B_{j}^{-1}$, $b_{kl}, b_{i}, c \in C^{0,\sigma}(\overline{\tilde{\Omega}})$ and $b_{k,l}(\mathcal B_{j}(S_{j})) = 0$.\\
By the assumption that $\dfrac{\tilde{\phi}_{j} - \tilde{\phi}_{j+1} + (2+\sigma)\tilde{\omega_{j}}}{\pi}$ is not an integer, we can apply Theorem 2.5.4 to get
\begin{align*}
\|v_{j}\|_{2,\sigma, \overline{\tilde{\Omega}}}
&\leq C\{\|g-(\sum_{k,l=1}^{2}b_{kl}D_{kl}v_{j} + b_{i}D_{i}v_{j} + cv_{j})\|_{\sigma, \overline{\tilde{\Omega}}} + \|v_{j}\|_{1,\sigma, \overline{\tilde{\Omega}}}\\
 &\quad + \sum_{k \in \mathcal D}\|v_{j}\|_{2,\sigma,\overline{\tilde{\Gamma}_{k}}} +  \sum_{k \in \mathcal N}\|\dfrac{\partial v_{j}}{\partial \tilde{\nu}_{k}} + \tilde{\beta}_{k}\dfrac{\partial v_{j}}{\partial \tilde{\tau}_{k}}\|_{1,\sigma,\overline{\tilde{\Gamma}_{k}}}     \}
\end{align*}
where $C$ does not depend on the choose of $\eta_{j}$. \\
So we may choose $\text{supp}(\eta_{j})$ small enough such that $b_{kl} \leq \frac{1}{8C}$ for any $x \in \mathcal B(\text{supp}(\eta_{j}))$, then we have
\begin{align*}
\|v_{j}\|_{2,\sigma, \overline{\tilde{\Omega}}}
&\leq C\{\|g\|_{\sigma, \overline{\tilde{\Omega}}} + \|v_{j}\|_{1,\sigma, \overline{\tilde{\Omega}}} + \sum_{k \in \mathcal D}\|v_{j}\|_{2,\sigma,\overline{\tilde{\Gamma}_{k}}} +  \sum_{k \in \mathcal N}\|\dfrac{\partial v_{j}}{\partial \tilde{\nu}_{k}} + \tilde{\beta}_{k}\dfrac{\partial v_{j}}{\partial \tilde{\tau}_{k}}\|_{1,\sigma,\overline{\tilde{\Gamma}_{k}}}  \} \\
&\quad + 4 C \cdot \dfrac{1}{8C}  \|v_{j}\|_{2,\sigma, \overline{\tilde{\Omega}}}, \\
\|v_{j}\|_{2,\sigma, \overline{\tilde{\Omega}}}
&\leq C^{'}\{\|g\|_{\sigma, \overline{\tilde{\Omega}}} + \|v_{j}\|_{1,\sigma, \overline{\tilde{\Omega}}} + \sum_{k \in \mathcal D}\|v_{j}\|_{2,\sigma,\overline{\tilde{\Gamma}_{k}}} +  \sum_{k \in \mathcal N}\|\dfrac{\partial v_{j}}{\partial \tilde{\nu}_{k}} + \tilde{\beta}_{k}\dfrac{\partial v_{j}}{\partial \tilde{\tau}_{k}}\|_{1,\sigma,\overline{\tilde{\Gamma}_{k}}}  \}.
\end{align*}
It follows that
$$\|\eta_{j}u\|_{2,\sigma, \overline{\Omega}} \leq C\{\|J(u)\|_{\sigma, \overline{\Omega}} +  \|u\|_{1,\sigma, \overline{\Omega}} + \sum_{j \in \mathcal D}\|u\|_{2,\sigma,\overline{\Gamma_{j}}} +  \sum_{j \in \mathcal N}\|\dfrac{\partial u}{\partial \nu_{j}} + \beta_{j}\dfrac{\partial u}{\partial \tau_{j}}\|_{1,\sigma,\overline{\Gamma_{j}}}    \}$$
where $C = C(\eta_{j}, \mathcal B_{j}, a_{ij}, b_{i}, c, \sigma)$. \\
Step 3: By Step 1 and Step 2, we have
\begin{align*}
  \|u\|_{2,\sigma, \overline{\Omega}} &= \|\sum_{j=0}^{N}\eta_{j}u\|_{2,\sigma, \overline{\Omega}} \\
         &\leq \|\eta_{0}u\|_{2,\sigma, \overline{\Omega}} + \sum_{j=1}^{N}\|\eta_{j}u\|_{2,\sigma, \overline{\Omega}} \\
         &\leq C\{\|J(u)\|_{\sigma, \overline{\Omega}} +  \|u\|_{1,\sigma, \overline{\Omega}} + \sum_{j \in \mathcal D}\|u\|_{2,\sigma,\overline{\Gamma_{j}}} \\
         &\quad+  \sum_{j \in \mathcal N}\|\dfrac{\partial u}{\partial \nu_{j}} + \beta_{j}\dfrac{\partial u}{\partial \tau_{j}}\|_{1,\sigma,\overline{\Gamma_{j}}}    \}.
\end{align*}
$\hfill$ q.e.d.\\

\subsubsection{Improved $C^{2,\sigma}$ estimates}
Combine Theorem 2.5.5 and Corollary 2.5.1, we have the following estimate.
\begin{thm}
Let $\Omega$ be a $C^{2,1}$ curvilinear polygonal domain, $\beta_{j}(x), d_{j} \in C^{1,\sigma}(\overline{\Gamma_{j}})$.\\
Assume $0 < \sigma < 1$ and $\dfrac{\tilde{\phi}_{j}(S_{j}) - \tilde{\phi}_{j+1}(S_{j}) + (2+\sigma)\tilde{\omega}_{j}}{\pi}$ is not an integer for any $j$. Then for any $\mu \leq -\text{sup}|c|-1$, $d \geq \text{sup}|d_{j}|+1$, there is a constant $C$ such that
$$\|u\|_{2,\sigma, \overline{\Omega}} \leq C\{\|J(u) + \mu u\|_{\sigma, \overline{\Omega}} + \sum_{j \in \mathcal D}\|u\|_{2,\sigma,\overline{\Gamma_{j}}} +  \sum_{j \in \mathcal N}\|\dfrac{\partial u}{\partial \nu_{j}} + \beta_{j}\dfrac{\partial u}{\partial \tau_{j}} + d_{j}u + d u\|_{1,\sigma,\overline{\Gamma_{j}}}    \}$$
for all $u \in C^{2,\sigma}(\overline{\Omega})$ such that $u(S_{j}) = 0$ if $j,j+1 \in \mathcal N$.
\end{thm}
\emph{Proof:}
The proof depends on Theorem 2.5.5, Corollary 2.5.1 and the inequality that for any $\epsilon > 0$, there is a constant $C$ such that for any $u \in C^{2,\sigma}(\overline{\Omega})$,
$$\|u\|_{1,\sigma, \overline{\Omega}} \leq \epsilon\|u\|_{2,\sigma, \overline{\Omega}} + C\|u\|_{0, \overline{\Omega}}.$$
\hfill q.e.d.

\subsection{Index of general elliptic operators on curvilinear polygonal domains}
From Theorem 2.4.6, we can derive an index theorem for the Laplace operator.
\begin{thm}
Let $\Omega$ be a $C^{2,1}$ curvilinear polygonal domain, $\beta_{j}(x) \in C^{2,1}(\overline{\Gamma_{j}})$, $d_{j} \in C_{c}^{1,1}(\Gamma_{j})$ such that there is a small positive $\epsilon$ with $\beta_{j}(x) = \beta_{j}(S_{j})$ for $x \in B_{\epsilon}(S_{j}) \cap \overline{\Gamma_{j}}$, $\beta_{j}(x) = \beta_{j}(S_{j+1})$ for $x \in B_{\epsilon}(S_{j+1}) \cap \overline{\Gamma_{j}}$, $d_{j} \geq 0$ and $\dfrac{1}{2}\cdot\dfrac{\partial \beta_{j}}{\partial \tau_{j}} + d_{j} \geq 0$.\\
Assume that $0< \sigma <1$ and $\dfrac{\phi_{j}(S_{j}) - \phi_{j+1}(S_{j}) + (2+\sigma)\omega_{j}}{\pi}$ is not an integer for any $j$. Then the operator
  $$L(u) = (\triangle u, \prod_{j \in \mathcal D}u\big|_{\Gamma_{j}}, \prod_{j \in \mathcal N}[\dfrac{\partial u}{\partial \nu_{j}} + \beta_{j}\dfrac{\partial u}{\partial \tau_{j}} + \dfrac{1}{2}\cdot\dfrac{\partial \beta_{j}}{\partial \tau_{j}}u + d_{j}u]\big|_{\Gamma_{j}})$$
  from $C^{2,\sigma}(\overline{\Omega})$ to $C^{0,\sigma}(\overline{\Omega}) \times \prod_{j \in \mathcal D}C^{2,\sigma}(\overline{\Gamma}_{j}) \times \prod_{j \in \mathcal N}C^{1,\sigma}(\overline{\Gamma}_{j})$ is a Fredholm map.
\end{thm}
\emph{Proof:} The space
$$\{(g_{j}, g_{j+1}) \in C^{2,\sigma}(\overline{\Gamma_{j}}) \times C^{2,\sigma}(\overline{\Gamma_{j+1}}) : g_{j}(S_{j}) = g_{j+1}(S_{j})\}$$
is of codimension-1 in $C^{2,\sigma}(\overline{\Gamma_{j}}) \times C^{2,\sigma}(\overline{\Gamma_{j+1}})$.\\
For any
$$(f,\prod_{j \in \mathcal D}g_{j}, \prod_{j \in \mathcal N}g_{j}) \in C^{0,\sigma}(\overline{\Omega}) \times \prod_{j \in \mathcal D}C^{2,\sigma}(\overline{\Gamma}_{j}) \times \prod_{j \in \mathcal N}C^{1,\sigma}(\overline{\Gamma}_{j})$$
such that $g_{j}(S_{j}) = g_{j+1}(S_{j})$ when $j, j+1 \in \mathcal D$, by Theorem 2.4.6, there exist real numbers $c_{j,m}$ and a function $u$ such that
$$L(u) = (f,\prod_{j \in \mathcal D}g_{j}, \prod_{j \in \mathcal N}g_{j}),$$
$$u-\sum_{-(2+\sigma)<\lambda_{j,m}<0}c_{j,m}\mathcal S_{j,m} \in C^{2,\sigma}(\overline{\Omega}).$$
Therefore we have
$$L(u-\sum_{-(2+\sigma)<\lambda_{j,m}<0}c_{j,m}\mathcal S_{j,m}) + \sum_{-(2+\sigma)<\lambda_{j,m}<0}c_{j,m}L(\mathcal S_{j,m}) = (f,\prod_{j \in \mathcal D}g_{j}, \prod_{j \in \mathcal N}g_{j}).$$
Thus $\text{coker}(L)$ is of finite dimension.\\
By Theorem 2.5.3,
\begin{align*}
& \quad \|u\|_{2,\sigma, \overline{\Omega}} \\
&\leq C\{\|\triangle u\|_{\sigma, \overline{\Omega}} + \|u\|_{1,\sigma, \overline{\Omega}} + \sum_{j \in \mathcal D}\|u\|_{2,\sigma,\overline{\Gamma_{j}}} +  \sum_{j \in \mathcal N}\|\dfrac{\partial u}{\partial \nu_{j}} + \beta_{j}\dfrac{\partial u}{\partial \tau_{j}}\|_{1,\sigma,\overline{\Gamma_{j}}}    \}\\
&\leq C^{'}\{\|\triangle u\|_{\sigma, \overline{\Omega}} + \|u\|_{1,\sigma, \overline{\Omega}} + \sum_{j \in \mathcal D}\|u\|_{2,\sigma,\overline{\Gamma_{j}}} \\
&\quad +  \sum_{j \in \mathcal N}\|\dfrac{\partial u}{\partial \nu_{j}} + \beta_{j}\dfrac{\partial u}{\partial \tau_{j}} + \dfrac{1}{2}\cdot\dfrac{\partial \beta_{j}}{\partial \tau_{j}}u + d_{j}u\|_{1,\sigma,\overline{\Gamma_{j}}}    \}.
\end{align*}
Then $\text{ker}(L)$ is also of finite dimension.
$\hfill$ q.e.d.\\

Next we derive the index theorem for general elliptic operators.
\begin{thm}
Let $\Omega$ be a $C^{2,1}$ curvilinear  polygonal domain, $\beta_{j}(x), d_{j} \in C^{1,\sigma}(\overline{\Gamma_{j}})$. Let $J(u) = a_{ij}D_{ij}u + b_{i}D_{i}u + cu$ be a strongly elliptic operator where $a_{ij} = a_{ji}, b_{i}, c \in C^{0,\sigma}(\overline{\Omega})$ and there is a positive $\alpha$ such that $\sum_{i,j=1}^{2}a_{ij}(x)\xi_{i}\xi_{j} \geq \alpha|\xi|^2$ for all $x \in \overline{\Omega}$ and $\xi \in \mathbb R^2$.\\
Assume $0 < \sigma < 1$ and $\dfrac{\tilde{\phi}_{j}(S_{j}) - \tilde{\phi}_{j+1}(S_{j}) + (2+\sigma)\tilde{\omega}_{j}}{\pi}$ is not an integer for any $j$. Then the operator
  $$L(u) = (J(u), \prod_{j \in \mathcal D}u\big|_{\Gamma_{j}}, \prod_{j \in \mathcal N}[\dfrac{\partial u}{\partial \nu_{j}} + \beta_{j}\dfrac{\partial u}{\partial \tau_{j}} + d_{j}u]\big|_{\Gamma_{j}})$$
  from $C^{2,\sigma}(\overline{\Omega})$ to $C^{0,\sigma}(\overline{\Omega}) \times \prod_{j \in \mathcal D}C^{2,\sigma}(\overline{\Gamma}_{j}) \times \prod_{j \in \mathcal N}C^{1,\sigma}(\overline{\Gamma}_{j})$ is a Fredholm map.
\end{thm}
\emph{Proof:} We consider two cases.\\
Case 1: $a_{kl}(S_{j}) = c_{j}\delta_{kl}$ for some positive constant $c_{j}$.\\
Let $\tilde{\beta}_{j}(x) \in C^{2,1}(\overline{\Gamma_{j}})$, $\tilde{d}_{j} \in C_{c}^{1,1}(\Gamma_{j})$ such that there is a small positive $\epsilon$ with $\tilde{\beta}_{j}(x) = \beta_{j}(S_{j})$ for $x \in B_{\epsilon}(S_{j}) \cap \overline{\Gamma_{j}}$, $\tilde{\beta}_{j}(x) = \beta_{j}(S_{j+1})$ for $x \in B_{\epsilon}(S_{j+1}) \cap \overline{\Gamma_{j}}$, $\tilde{d}_{j} \geq 0$ and $\dfrac{1}{2}\cdot\dfrac{\partial \tilde{\beta}_{j}}{\partial \tau_{j}} + \tilde{d}_{j} \geq 0$.\\
Then by Theorem 2.6.1,
$$L_{0}(u) = (\triangle u, \prod_{j \in \mathcal D}u\big|_{\Gamma_{j}}, \prod_{j \in \mathcal N}[\dfrac{\partial u}{\partial \nu_{j}} + \tilde{\beta}_{j}\dfrac{\partial u}{\partial \tau_{j}} + \dfrac{1}{2} \cdot \dfrac{\partial \tilde{\beta}_{j}}{\partial \tau_{j}}u + \tilde{d}_{j}u]\big|_{\Gamma_{j}})$$
is a Fredholm map.\\
Let $M > \max\{\sup|c|+1, \sup|d_{j}|+1\}$, define
\begin{align*}
L_{t}(u)
&= tL(u) + (1-t)L_{0}(u) + (-M u, \prod_{j \in \mathcal D} 0,\prod_{j \in \mathcal N} M u)\\
&= (J_{t}(u), \prod_{j \in \mathcal D}u\big|_{\Gamma_{j}}, \prod_{j \in \mathcal N}Q_{j,t}(u)\big|_{\Gamma_{j}}).
\end{align*}
Then the principal part of $J_{t}$ at each $S_{j}$ is a multiple of the Laplace operator. Thus by Theorem 2.5.6, there is a constant $C(t)$ which depends on $t$ continuously such that
$$\|u\|_{2,\sigma, \overline{\Omega}} \leq C(t)\{\|J_{t}(u)\|_{\sigma, \overline{\Omega}} + \sum_{j \in \mathcal D}\|u\|_{2,\sigma,\overline{\Gamma_{j}}}
 +  \sum_{j \in \mathcal N}\|Q_{j,t}(u)\|_{1,\sigma,\overline{\Gamma_{j}}}    \}$$
for all $u \in C^{2,\sigma}(\overline{\Omega})$ such that $u(S_{j}) = 0$ if $j,j+1 \in \mathcal N$.\\
By a theorem in Kato(1966), the index of $L_{t}$ is invariant.\\
Since $\{u \in C^{2,\sigma}(\overline{\Omega}): u(S_{j}) = 0 \text{ if $j,j+1 \in \mathcal N$}\}$ is of finite codimension in $C^{2,\sigma}(\overline{\Omega})$ and $(-M u, 0, M u)$ is a compact operator, $L$ is also of finite index.\\
Case 2: Similar to the proof in Theorem 2.5.5, there are symmetric $2 \times 2$ matrices $\mathcal B_{j}$ with eigenvalues in $(0,1]$ such that
$$\mathcal B_{j} \mathcal A_{j} \mathcal B_{j} =  c_{j}I$$
where $c_{j}$ are some positive constants.\\
By Theorem 2.9.2, there is a diffeomorphism $\Phi : \mathbb R^2 \to \mathbb R^2$ such that\\
(1) $\Phi(x,y) = S_{j} +  B_{j}\left\{\begin{bmatrix}
    x\\
    y
  \end{bmatrix}-S_{j}\right\}$ on $B_{\epsilon}(S_{j})$.\\
(2) $\Phi(x,y) = (x,y)$ outside $\cup_{j} B_{2\epsilon}(S_{j})$. \\
Thus we reduce our problem to Case 1 in $\Phi(\Omega)$
$\hfill$ q.e.d.\\

\subsection{Regular solutions}
In this section, we consider the solvability of the elliptic problem in some special cases.

\subsubsection{When $\omega_{j}$ is small}
\begin{thm}
  Let $\Omega$ be a $C^{k+2,1}$ curvilinear polygonal domain, $\beta_{j}(x) \in C^{k+2,1}(\overline{\Gamma_{j}})$, $d_{j} \in C_{c}^{k+1,1}(\Gamma_{j})$ such that there is a small positive $\epsilon$ with $\beta_{j}(x) = \beta_{j}(S_{j})$ for $x \in B_{\epsilon}(S_{j}) \cap \overline{\Gamma_{j}}$, $\beta_{j}(x) = \beta_{j}(S_{j+1})$ for $x \in B_{\epsilon}(S_{j+1}) \cap \overline{\Gamma_{j}}$, $d_{j} \geq 0$ and $\dfrac{1}{2}\cdot\dfrac{\partial \beta_{j}}{\partial \tau_{j}} + d_{j} \geq 0$.\\
Assume $0< \sigma <1$, then there is a positive angle $\omega_{0}$ such that if $0 < \omega_{j} < \omega_{0}$, for each $f \in C^{k,\sigma}(\overline{\Omega})$, $g_{j} \in C^{k+2,\sigma}(\overline{\Gamma_{j}})$ if $j \in \mathcal D$, $g_{j} \in C^{k+1,\sigma}(\overline{\Gamma_{j}})$ if $j \in \mathcal N$ and $g_{j}(S_{j}) = g_{j+1}(S_{j})$ if $j,j+1 \in \mathcal D$, there exists a (possibly non-unique) $u \in C^{k+2,\sigma}(\overline{\Omega})$ which solves the following problem
$$(***)\begin{cases}
     \triangle u = f      & \text{in $\Omega$},\\
     u = g_{j}         & \text{on $\Gamma_{j}$, $j \in \mathcal D$}, \\
     \dfrac{\partial u}{\partial \nu_{j}} + \beta_{j}\dfrac{\partial u}{\partial \tau_{j}} + \dfrac{1}{2}\cdot\dfrac{\partial \beta_{j}}{\partial \tau_{j}}u + d_{j}u =g_{j} & \text{on $\Gamma_{j}$, $j \in \mathcal N$}.
  \end{cases}$$
\end{thm}
\emph{Proof:} Since $0<\omega_{j}<\omega_{0}$ is small enough, we know $$\dfrac{\phi_{j}(S_{j}) - \phi_{j+1}(S_{j}) + (2+\sigma)\omega_{j}}{\pi}$$
is not an integer for any $j$. \\
By Theorem 2.4.6, there exists a $u$ which is a solution to problem ($***$) and
$$u - \sum_{-(k+2+\sigma) < \lambda_{j,m} < 0}c_{j,m} \mathcal S_{j,m} \in C^{k+2,\sigma}(\overline{\Omega}).$$
If $\omega_{j}$ is small enough, then
$$\{\lambda_{j,m} : -(k+2+\sigma) < \lambda_{j,m} < 0\} = \emptyset,$$
which means $u \in C^{k+2,\sigma}(\overline{\Omega})$.
\hfill q.e.d.\\

Theorem 2.7.1 is similar to the results in Azzam-Kreyszig\cite{AzzamKreyszig} concerning small angles. Next we derive the existence of regular solutions for general elliptic operators using the same continuity method as in the proofs of the index theorems in Section 2.6.
\begin{thm}
  Let $\Omega$ be a $C^{2,1}$ curvilinear  polygonal domain, $\beta_{j}(x), d_{j} \in C^{1,\sigma}(\overline{\Gamma_{j}})$. Let $J(u) = a_{ij}D_{ij}u + b_{i}D_{i}u + cu$ be a strongly elliptic operator where $a_{ij} = a_{ji}, b_{i}, c \in C^{0,\sigma}(\overline{\Omega})$ and there is a positive $\alpha$ such that $\sum_{i,j=1}^{2}a_{ij}(x)\xi_{i}\xi_{j} \geq \alpha|\xi|^2$ for all $x \in \overline{\Omega}$ and $\xi \in \mathbb R^2$.\\
  Assume $0 < \sigma < 1$, $c \leq -1$, $d_{j} \geq 1$, then there is a positive angle $\omega_{0}$ such that if $0 < \omega_{j} < \omega_{0}$, for each $f \in C^{0,\sigma}(\overline{\Omega})$, $g_{j} \in C^{2,\sigma}(\overline{\Gamma_{j}})$ if $j \in \mathcal D$, $g_{j} \in C^{1,\sigma}(\overline{\Gamma_{j}})$ if $j \in \mathcal N$ and $g_{j}(S_{j}) = g_{j+1}(S_{j})$ if $j,j+1 \in \mathcal D$, there exists a unique $u \in C^{2,\sigma}(\overline{\Omega})$ which solves the following problem
$$\begin{cases}
     J u = f      & \text{in $\Omega$},\\
     u = g_{j}         & \text{on $\Gamma_{j}$, $j \in \mathcal D$}, \\
     \dfrac{\partial u}{\partial \nu_{j}} + \beta_{j}\dfrac{\partial u}{\partial \tau_{j}} + d_{j}u=g_{j} & \text{on $\Gamma_{j}$, $j \in \mathcal N$}, \\
     u(S_{j}) = 0  & \text{if $\{j,j+1\} \subset \mathcal N$}.
  \end{cases}$$
\end{thm}

Remark: We assume that $c \leq -1$, $d_{j} \geq 1$ and $u(S_{j}) = 0$ if $\{j,j+1\} \subset \mathcal N$ to make sure that we have the estimate
$$\|u\|_{2,\sigma, \overline{\Omega}} \leq C\{\|J(u)\|_{\sigma, \overline{\Omega}}+ \sum_{j \in \mathcal D}\|u\|_{2,\sigma,\overline{\Gamma_{j}}} +  \sum_{j \in \mathcal N}\|\dfrac{\partial u}{\partial \nu_{j}} + \beta_{j}\dfrac{\partial u}{\partial \tau_{j}} + d_{j}u\|_{1,\sigma,\overline{\Gamma_{j}}}    \}.$$

\subsubsection{When $\omega_{j} = \dfrac{\pi}{2}$}
We consider the case when $\omega_{j} = \dfrac{\pi}{2}$, one of $\{j,j+1\}$ is in $\mathcal D$, the other one is in $\mathcal N$ and $\beta_{j} = 0$. For any real numbers $d_{j}$, we say $g_{j} \in C^{2,\sigma}(\overline{\Gamma_{j}})$ if $j \in \mathcal D$ and $g_{j} \in C^{1,\sigma}(\overline{\Gamma_{j}})$ if $j \in \mathcal N$ are compatible at $S_{j}$ if
$$\begin{cases}
  g_{j+1}(S_{j}) = \dfrac{\partial g_{j}}{\partial \tau_{j}}(S_{j}) + d_{j+1}g_{j}(S_{j})  &\text{ \quad if $j \in \mathcal D$, $j+1 \in \mathcal N$},\\
  g_{j}(S_{j}) = -\dfrac{\partial g_{j+1}}{\partial \tau_{j+1}}(S_{j}) + d_{j}g_{j+1}(S_{j})  &\text{ \quad if $j \in \mathcal N$, $j+1 \in \mathcal D$}.
\end{cases}$$

As a direct corollary of Theorem 2.4.5, we have the following theorem.
\begin{thm}
  Let $\Omega$ be a $C^{2,1}$ curvilinear polygonal domain. Assume $\omega_{j} = \dfrac{\pi}{2}$, one of $\{j,j+1\}$ is in $\mathcal D$ and the other one is in $\mathcal N$.\\
Then for $0 < \sigma < 1$ and each $f \in C^{0,\sigma}(\overline{\Omega})$, there exists a function $u \in C^{2,\sigma}(\overline{\Omega})$ which is a solution to
$$(*)\begin{cases}
     \triangle u = f      & \text{in $\Omega$},\\
     u = 0         & \text{on $\Gamma_{j}$, $j \in \mathcal D$}, \\
     \dfrac{\partial u}{\partial \nu_{j}} = 0 & \text{on $\Gamma_{j}$, $j \in \mathcal N$}.
  \end{cases}$$
\end{thm}
\emph{Proof:} The proof is similar to Theorem 2.7.1. By the assumptions, we have $\dfrac{\phi_{j} - \phi_{j+1} + (2+\sigma)\omega_{j}}{\pi} = \dfrac{\pm \frac{\pi}{2} + (2+\sigma)\frac{\pi}{2}}{\pi}$ is not an integer. Apply Theorem 2.4.5, there exist real numbers $c_{j,m}$ and a function $u$ which solves problem ($*$) such that
$$u - \sum_{-(2+\sigma) < \lambda_{j,m} < 0, \lambda_{j,m} \neq -1}c_{j,m} \mathcal S_{j,m} \in C^{2,\sigma}(\overline{\Omega}).$$
Since
$$\lambda_{j,m} = \dfrac{\phi_{j}-\phi_{j+1}+m\pi}{\omega_{j}} = \dfrac{\pm \frac{\pi}{2}+m\pi}{\frac{\pi}{2}} = \pm 1 + 2m,$$
we have
$$\{\lambda_{j,m} : -(2+\sigma) < \lambda_{j,m} < 0, \lambda_{j,m} \neq -1\} = \emptyset,$$
which means $u \in C^{2,\sigma}(\overline{\Omega})$.
\hfill q.e.d \\

Remark: In this case, suppose $\omega_{j} = \dfrac{\pi}{2}$, $j \in \mathcal D$, $j+1 \in \mathcal N$ and $u \in C^{2,\sigma}(\overline{\Omega})$, then there are compatibility conditions of $u$ at the corners. The singular solution $\mathcal S_{j,m}$ when $\lambda_{j,m} = -1$ is to balance the incompatibility of the given boundary values.

With the above remark, we have the following existence theorem of regular solutions if the given boundary values are compatible at the corners.

\begin{thm}
Let $\Omega$ be a $C^{2,1}$ curvilinear polygonal domain. Assume $\omega_{j} = \dfrac{\pi}{2}$, one of $\{j,j+1\}$ is in $\mathcal D$ and the other one is in $\mathcal N$.\\
Then for $0 < \sigma < 1$ and each $f \in C^{0,\sigma}(\overline{\Omega})$, $g_{j} \in C^{2,\sigma}(\overline{\Gamma_{j}})$ if $j \in \mathcal D$, $g_{j} \in C^{1,\sigma}(\overline{\Gamma_{j}})$ if $j \in \mathcal N$, $g_{j}$ and $g_{j+1}$ are compatible at $S_{j}$, there exists a function $u \in C^{2,\sigma}(\overline{\Omega})$ which is a solution to
$$(**) \begin{cases}
     \triangle u = f      & \text{in $\Omega$},\\
     u = g_{j}         & \text{on $\Gamma_{j}$, $j \in \mathcal D$}, \\
     \dfrac{\partial u}{\partial \nu_{j}} = g_{j} & \text{on $\Gamma_{j}$, $j \in \mathcal N$}.
  \end{cases}$$
\end{thm}
\emph{Proof:} Since $g_{j}$ and $g_{j+1}$ are compatible at $S_{j}$, by the inverse trace theorems in Section 2.8, there is a $w \in C^{2,\sigma}(\overline{\Omega})$ such that
$$\begin{cases}
     w = g_{j}         & \text{on $\Gamma_{j}$, $j \in \mathcal D$}, \\
     \dfrac{\partial w}{\partial \nu_{j}} = g_{j} & \text{on $\Gamma_{j}$, $j \in \mathcal N$}.
  \end{cases}$$
By Theorem 2.7.3, there exists a $v \in C^{2,\sigma}(\overline{\Omega})$ such that
 $$\begin{cases}
     \triangle v = f- \triangle w     & \text{in $\Omega$},\\
     v = 0         & \text{on $\Gamma_{j}$, $j \in \mathcal D$}, \\
     \dfrac{\partial v}{\partial \nu_{j}} = 0 & \text{on $\Gamma_{j}$, $j \in \mathcal N$}.
  \end{cases}$$
Thus $u = w + v$ is what we need.
\hfill q.e.d.\\

Next we have the index of the Laplace operator with the Dirichlet boundary condition and the Neumann boundary condition.
\begin{thm}
  Let $\Omega$ be a $C^{2,1}$ curvilinear polygonal domain. Assume $\omega_{j} = \dfrac{\pi}{2}$, one of $\{j,j+1\}$ is in $\mathcal D$ and the other one is in $\mathcal N$.\\
  Then the operator
  $$L(u) = (\triangle u, \prod_{j \in \mathcal D}u\big|_{\Gamma_{j}}, \prod_{j \in \mathcal N}\dfrac{\partial u}{\partial \nu_{j}}\big|_{\Gamma_{j}})$$
  from $C^{2,\sigma}(\overline{\Omega})$ to $C^{0,\sigma}(\overline{\Omega}) \times \prod_{j \in \mathcal D}C^{2,\sigma}(\overline{\Gamma}_{j}) \times \prod_{j \in \mathcal N}C^{1,\sigma}(\overline{\Gamma}_{j})$ is a Fredholm map of index equal to $-I_{0}$.
\end{thm}
\emph{Proof:}
By Theorem 2.7.4, for each $f \in C^{0,\sigma}(\overline{\Omega})$, $g_{j} \in C^{2,\sigma}(\overline{\Gamma_{j}})$ if $j \in \mathcal D$, $g_{j} \in C^{1,\sigma}(\overline{\Gamma_{j}})$ if $j \in \mathcal N$, $g_{j}$ and $g_{j+1}$ are compatible at $S_{j}$, there exists a function $u \in C^{2,\sigma}(\overline{\Omega})$ which is a solution to
$$\begin{cases}
     \triangle u = f      & \text{in $\Omega$},\\
     u = g_{j}         & \text{on $\Gamma_{j}$, $j \in \mathcal D$}, \\
     \dfrac{\partial u}{\partial \nu_{j}} = g_{j} & \text{on $\Gamma_{j}$, $j \in \mathcal N$}.
  \end{cases}$$
Thus $\dim(\text{coker}(L)) = I_{0}$ since for any $u \in C^{2,\sigma}(\overline{\Omega})$ the boundary value of $u$ must be compatible at the corners and $g_{j} \in C^{2,\sigma}(\overline{\Gamma_{j}})$ if $j \in \mathcal D$, $g_{j} \in C^{1,\sigma}(\overline{\Gamma_{j}})$ if $j \in \mathcal N$ such that $g_{j}$ and $g_{j+1}$ is compatible at $S_{j}$ is of codimension-$I_{0}$ in $\prod_{j \in \mathcal D}C^{2,\sigma}(\overline{\Gamma}_{j}) \times \prod_{j \in \mathcal N}C^{1,\sigma}(\overline{\Gamma}_{j})$.\\
By the estimate in Theorem 2.5.5, we have
$$\|u\|_{2,\sigma, \overline{\Omega}} \leq C\{\|\triangle u\|_{\sigma, \overline{\Omega}} + \|u\|_{1,\sigma, \overline{\Omega}} + \sum_{j \in \mathcal D}\|u\|_{2,\sigma,\overline{\Gamma_{j}}} +  \sum_{j \in \mathcal N}\|\dfrac{\partial u}{\partial \nu_{j}}\|_{1,\sigma,\overline{\Gamma_{j}}}    \},$$
which means $\text{ker}(L)$ is of finite dimension.\\
Thus $\text{index}(L) = \dim(\text{ker}(L)) - \dim(\text{coker}(L)) \geq 0-I_{0} = -I_{0}$.\\
On the other hand, define
$$\tilde{L}(u) = (\triangle u - u, \prod_{j \in \mathcal D}u\big|_{\Gamma_{j}}, \prod_{j \in \mathcal N}[\dfrac{\partial u}{\partial \nu_{j}} + u]\big|_{\Gamma_{j}}).$$
We have $\text{index}(\tilde{L}) = \text{index}(L)$.\\
By the estimate in Theorem 2.5.6, for any $u \in C^{2,\sigma}(\overline{\Omega})$,
$$\|u\|_{2,\sigma, \overline{\Omega}} \leq C\{\|\triangle u - u\|_{\sigma, \overline{\Omega}} + \sum_{j \in \mathcal D}\|u\|_{2,\sigma,\overline{\Gamma_{j}}} +  \sum_{j \in \mathcal N}\|\dfrac{\partial u}{\partial \nu_{j}} + u\|_{1,\sigma,\overline{\Gamma_{j}}}    \}$$
which means $\tilde{L}$ is one-to-one. Since the boundary values of $u$ must be compatible at the corners, $\dim(\text{coker}(\tilde{L})) \geq I_{0}$.\\
Then we have $\text{index}(\tilde{L}) = \dim(\text{ker}(\tilde{L})) - \dim(\text{coker}(\tilde{L})) \leq 0-I_{0} = -I_{0}$. And
$$-I_{0} \geq \text{index}(\tilde{L}) = \text{index}(L) \geq -I_{0}.$$
We get $\text{index}(L) = -I_{0}$.
\hfill q.e.d.\\

Finally we have the following existence theorem for general elliptic operators.
\begin{thm}
   Let $\Omega$ be a $C^{2,1}$ curvilinear  polygonal domain, $\beta_{j}(x), d_{j} \in C^{1,\sigma}(\overline{\Gamma_{j}})$. Let $J(u) = a_{ij}D_{ij}u + b_{i}D_{i}u + cu$ be a strongly elliptic operator where $a_{ij} = a_{ji}, b_{i}, c \in C^{0,\sigma}(\overline{\Omega})$ and there is a positive $\alpha$ such that $\sum_{i,j=1}^{2}a_{ij}(x)\xi_{i}\xi_{j} \geq \alpha|\xi|^2$ for all $x \in \overline{\Omega}$ and $\xi \in \mathbb R^2$.\\
  Assume $0 < \sigma < 1$, $\omega_{j} = \dfrac{\pi}{2}$, one of $\{j,j+1\}$ is in $\mathcal D$ and the other one is in $\mathcal N$, $c \leq -1$, $d_{j} \geq 1$, $\beta_{j}(S_{j}) = \beta_{j}(S_{j+1}) = 0$, $a_{kl}(S_{j}) = \delta_{kl}(S_{j})$.\\
   Then for each $f \in C^{0,\sigma}(\overline{\Omega})$, $g_{j} \in C^{2,\sigma}(\overline{\Gamma_{j}})$ if $j \in \mathcal D$, $g_{j} \in C^{1,\sigma}(\overline{\Gamma_{j}})$ if $j \in \mathcal N$ such that $g_{j}$ and $g_{j+1}$ are compatible at $S_{j}$, there exists a unique $u \in C^{2,\sigma}(\overline{\Omega})$ which is a solution to
$$\begin{cases}
     J u = f      & \text{in $\Omega$},\\
     u = g_{j}         & \text{on $\Gamma_{j}$, $j \in \mathcal D$}, \\
     \dfrac{\partial u}{\partial \nu_{j}} + \beta_{j}\dfrac{\partial u}{\partial \tau_{j}} + d_{j}u=g_{j} & \text{on $\Gamma_{j}$, $j \in \mathcal N$}.
  \end{cases}$$
\end{thm}
\emph{Proof:} Define
$$L_{0}(u) = (\triangle u - u, \prod_{j \in \mathcal D}u\big|_{\Gamma_{j}}, \prod_{j \in \mathcal N}[\dfrac{\partial u}{\partial \nu_{j}} + u]\big|_{\Gamma_{j}}).$$
Thus by Theorem 2.7.5, $\text{index}(L_{0}) = -I_{0}$.\\
Define
$$L_{1} = (J(u), \prod_{j \in \mathcal D}u\big|_{\Gamma_{j}}, \prod_{j \in \mathcal N}[\dfrac{\partial u}{\partial \nu_{j}} + \beta_{j}\dfrac{\partial u}{\partial \tau_{j}} + d_{j}u]\big|_{\Gamma_{j}})$$
and
\begin{align*}
  L_{t}(u)
  &= (1-t)L_{0}(u) + tL_{1}(u) \\
  &= (J_{t}(u), \prod_{j \in \mathcal D}u\big|_{\Gamma_{j}}, \prod_{j \in \mathcal N}Q_{j,t}(u)\big|_{\Gamma_{j}}) \text{ \quad for $t \in [0,1]$.}
\end{align*}
By the estimate in Theorem 2.5.6, there is a $C_{t}$ which depends on $t$ continuously such that for any $u \in C^{2,\sigma}(\overline{\Omega})$,
$$\|u\|_{2,\sigma, \overline{\Omega}} \leq C\{\|J_{t}(u)\|_{\sigma, \overline{\Omega}}+ \sum_{j \in \mathcal D}\|u\|_{2,\sigma,\overline{\Gamma_{j}}} +  \sum_{j \in \mathcal N}\|Q_{j,t}(u)\|_{1,\sigma,\overline{\Gamma_{j}}}    \}.$$
Then we have $\text{index}(L_{1}) = \text{index}(L_{0}) = -I_{0}$.\\
From the previous estimate, $L_{1}$ is one-to-one. Then $\text{image}(L_{1})$ is of codimension-$I_{0}$ in $$C^{0,\sigma}(\overline{\Omega}) \times \prod_{j \in \mathcal D}C^{2,\sigma}(\overline{\Gamma}_{j}) \times \prod_{j \in \mathcal N}C^{1,\sigma}(\overline{\Gamma}_{j}).$$
We denote
$$C_{*} = \{(f, \prod_{j \in \mathcal D}g_{j}, \prod_{j \in \mathcal N} g_{j})\} \subset C^{0,\sigma}(\overline{\Omega}) \times \prod_{j \in \mathcal D}C^{2,\sigma}(\overline{\Gamma}_{j}) \times \prod_{j \in \mathcal N}C^{1,\sigma}(\overline{\Gamma}_{j})$$
such that $g_{j}$ and $g_{j+1}$ are compatible at $S_{j}$. \\
Then $C_{*}$ is also of codimension-$I_{0}$ in
$$C^{0,\sigma}(\overline{\Omega}) \times \prod_{j \in \mathcal D}C^{2,\sigma}(\overline{\Gamma}_{j}) \times \prod_{j \in \mathcal N}C^{1,\sigma}(\overline{\Gamma}_{j}).$$
Since $\text{image}(L_{1}) \subset C_{*}$ and both are closed, we have $$\text{image}(L_{1}) = C_{*}.$$
\hfill q.e.d.\\

Remark: As in Theorem 2.7.1, we assume that $c \leq -1$, $d_{j} \geq 1$ to make sure we have the estimate
$$\|u\|_{2,\sigma, \overline{\Omega}} \leq C\{\|J(u)\|_{\sigma, \overline{\Omega}}+ \sum_{j \in \mathcal D}\|u\|_{2,\sigma,\overline{\Gamma_{j}}} +  \sum_{j \in \mathcal N}\|\dfrac{\partial u}{\partial \nu_{j}} + \beta_{j}\dfrac{\partial u}{\partial \tau_{j}} + d_{j}u\|_{1,\sigma,\overline{\Gamma_{j}}}    \}.$$

The reason we assume that $a_{kl}(S_{j}) = \delta_{kl}(S_{j})$ is that when we apply the continuity method, we cannot guarantee that $\dfrac{\tilde{\phi}_{j}(S_{j}) - \tilde{\phi}_{j+1}(S_{j})+ (2+\sigma)\tilde{\omega}_{j}}{\pi}$ is not an integer for any $t$ if $a_{kl}(S_{j}) \neq \delta_{kl}(S_{j})$.

\subsection{Inverse trace theorems}
To extend Grisvard's apriori estimates for homogeneous boundary values to non-homogeneous boundary values, we need inverse trace theorems in a curvilinear polygonal domain. We first prove a sequence of inverse trace theorems on the first quadrant $\mathbb R_{+} \times \mathbb R_{+}$, the upper half space $\mathbb R_{+}^{2}$, the complement of the first quadrant $\mathbb R^2 \setminus \mathbb R_{+} \times \mathbb R_{+}$ concerning the mixed boundary conditions. The method we use here is a modification of Theorem 5.8 in \cite{Necas}.

First we state an important extension theorem.
\begin{thm}[Grisvard\cite{GrisvardPierre2011EPiN}, Theorem 6.2.5]
  Let $\Omega$ be any bounded open domain in $\mathbb R^2$ with a polygonal boundary, there is a continuous linear operator $P: C^{m,\sigma}(\overline{\Omega}) \to C^{m,\sigma}(\mathbb R^2)$, for every $\sigma \in [0,1]$, such that for any $u \in C^{m,\sigma}(\overline{\Omega})$
  $$P(u)\big|_{\overline{\Omega}} = u.$$
  In particular, the result still hold if $\Omega$ is replaced by $\mathbb R_{+} \times \mathbb R_{+}$, $\mathbb R_{+}^{2}$ or $\mathbb R^2 \setminus (\mathbb R_{+} \times \mathbb R_{+})$ and $u \in $ $C_{c}^{m,\sigma}(\overline{\mathbb R_{+} \times \mathbb R_{+}})$, $C_{c}^{m,\sigma}(\overline{\mathbb R_{+}^{2}})$ or $C_{c}^{m,\sigma}(\overline{\mathbb R^2 \setminus (\mathbb R_{+} \times \mathbb R_{+})})$.
\end{thm}

Let $\phi(t) \in C_{c}^{\infty}(\mathbb R, [0,1])$ be a cut-off function such that \\
(1) $\phi(t) = 1$ for $t \in [-\frac{1}{4},\frac{1}{4}]$ and $\phi$ vanishes outside $[-1,1]$. \\
(2) $\int_{-\infty}^{\infty} \phi(t)dt = 1$.

Let $\{c_{j}\}$ satisfies $\sum_{j=1}^{3}(-j)^{l}c_{j} = 1$ for any $l=0,1,2$.
We will use them frequently in this section.

Let $C_{c}^{1,\sigma}(\mathbb R) = \{g \in C^{1,\sigma}(\mathbb R) \, \big| \, \text{supp}(g)$ is compact$\}$. We have the following theorem.
\begin{lemma}
There is a bounded linear operator
$$T_{X}:C_{c}^{1,\sigma}(\mathbb R) \to C_{c}^{2,\sigma}(\mathbb R^2)$$
such that for any $g \in C_{c}^{1,\sigma}(\mathbb R)$,
$$[T_{X}(g)](x,0)=0,$$
$$\left[-\dfrac{\partial (T_{X}(g))}{\partial y}\right]_{(x,0)}=g(x).$$
Similarly, there is a bounded linear operator
$$T_{Y}:C_{c}^{1,\sigma}(\mathbb R) \to C_{c}^{2,\sigma}(\mathbb R^2)$$
such that for any $h \in C_{c}^{1,\sigma}(\mathbb R)$,
$$[T_{Y}(h)](0,y)=0,$$
$$\left[-\dfrac{\partial (T_{Y}(h))}{\partial x}\right]_{(0,y)}=h(y).$$
\end{lemma}
\emph{Proof:}
Define
$$u(x,y) = \int_{\mathbb R} \phi(\frac{s-x}{-y})g(s)ds,$$
$$T_{X}(g)(x,y) = \phi(y)u(x,y).$$
We claim that $T_{X}$ is what we need.\\
Since $\phi \in C_{c}^{\infty}(\mathbb R)$ and $g \in C_{c}^{1,\sigma}(\mathbb R)$, $u(x,y)$ is smooth when $y \neq 0$. \\
(1) We show that $T_{X}(g) \in C_{c}^{2,\sigma}(\mathbb R^2)$ and $\|T_{X}\|$ is bounded.
\begin{align*}
  u(x,y)     &= \int_{\mathbb R} \phi(\frac{s-x}{-y})g(s)ds \\
             &= \int_{\mathbb R} -y\phi(t)g(x-ty)dt,
\end{align*}
\begin{align*}
  u_{x}(x,y) &= \int_{\mathbb R} \frac{1}{y}\phi^{'}(\frac{s-x}{-y})g(s)ds \\
             &= \int_{\mathbb R} \frac{1}{y}\phi^{'}(t)g(x-ty)(-y)dt \\
             &= -\int_{\mathbb R} \phi^{'}(t)g(x-ty)dt,
\end{align*}
\begin{align*}
  u_{y}(x,y) &= \int_{\mathbb R} \left[\frac{s-x}{y^2}\right]\phi^{'}(\frac{s-x}{-y})g(s)ds \\
             &= \int_{\mathbb R} \left[-\frac{t}{y}\right]\phi^{'}(t)g(x-ty)(-y)dt \\
             &= \int_{\mathbb R} t\phi^{'}(t)g(x-ty)dt,
\end{align*}
\begin{align*}
  u_{xx}(x,y)&= -\int_{\mathbb R} \phi^{'}(t)g^{'}(x-ty)dt, \\
  u_{xy}(x,y)&= \int_{\mathbb R} t\phi^{'}(t)g^{'}(x-ty)dt, \\
  u_{yy}(x,y)&= -\int_{\mathbb R} t^2\phi^{'}(t)g^{'}(x-ty)dt.
\end{align*}
From the above formulas, $u,Du,D^{2}u$ are all continuous up to $y=0$ and $D^{2}u$ is locally $C^{0,\sigma}$. Thus
$$T_{X}(g) \in C^{2,\sigma}(\mathbb R^2).$$
By the choice of $\phi(y)$, $\text{supp}(\phi) \subset [-1,1]$, we have $\text{supp}(T_{X}(g)) \subset \{(x,y):|y|\leq 1\}$.\\
For $|y|\leq 1$, $g \in C_{c}^{1,\sigma}(\mathbb R)$, $u(x,y) = \int_{\mathbb R} -y\phi(t)g(x-ty)dt$, then
$$\text{supp}(T_{X}(g)) \text{ is compact}.$$
Since $\text{supp}(\phi) \subset [-1,1]$,
\begin{align*}
|T_{X}(g)| &\leq \max_{|y|\leq 1} |\phi(y)\int_{\mathbb R} -y\phi(t)g(x-ty)dt|\\
           &\leq \int_{-1}^{1} |g(x-ty)|dt \\
           &\leq 2\|g\|_{0},
\end{align*}
\begin{align*}
|D[T_{X}(g)]| &\leq \max_{|y|\leq 1}|D\phi||u|+|\phi||Du|\\
           &\leq C(\phi)\|g\|_{0}+C(\phi)\max_{|y|\leq 1} |Du|\\
           &\leq C(\phi)\|g\|_{0}+C(\phi)\|g\|_{0}\\
           &\leq C(\phi)\|g\|_{0},
\end{align*}
\begin{align*}
\|D^{2}[T_{X}(g)]\|_{0,\sigma} &\leq C(\phi)\|u\|_{2,\sigma,\{|y|\leq 1\}}\\
               &\leq C(\phi)\|g\|_{1,\sigma}.
\end{align*}
Thus, we have
$$\|T_{X}(g)\|_{2,\sigma} \leq C(\phi)\|g\|_{1,\sigma}.$$
(2) We show that $T_{X}(g)$ satisfies the conditions on $\{(x,0)\}$,
\begin{align*}
  [T_{X}(g)](x,0) &= \phi(0)u(x,0) \\
                  &= \int_{\mathbb R} -0\phi(t)g(x-ty)dt\\
                  &= 0.
\end{align*}
Since $\phi(y) = 1$ near $y=0$,
\begin{align*}
  \left[-\dfrac{\partial (T_{X}(g))}{\partial y}\right]_{(x,0)}
  &= \left[-\dfrac{\partial u}{\partial y}\right]_{(x,0)} \\
  &= \left[-\int_{\mathbb R} t\phi^{'}(t)g(x-ty)dt\right]_{(x,0)} \\
  &= g(x)\int_{\mathbb R} -t\phi^{'}(t)dt \\
  &= g(x)\int_{\mathbb R} \phi(t)dt\\
  &= g(x).
\end{align*}
The construction of $T_{Y}$ is similar to $T_{X}$.
$\hfill$ q.e.d.\\

\subsubsection{On $\mathbb R_{+} \times \mathbb R_{+}$}
Let $C_{c}^{k,\sigma}(\overline{\mathbb R_{+} \times \mathbb R_{+}}) = \{g \in C^{k,\sigma}(\overline{\mathbb R_{+} \times \mathbb R_{+}}) \, \big| \, \text{supp}(g)$ is compact$\}$. We have the following theorems.
\begin{thm}
  There is a constant $C$ such that for any $u \in C_{c}^{2,\sigma}(\overline{\mathbb R_{+} \times \mathbb R_{+}})$, there exists a $w \in C_{c}^{2,\sigma}(\overline{\mathbb R_{+} \times \mathbb R_{+}})$ with
  $$u(x,0) = w(x,0) \text{ for $x \geq 0$},$$
  $$u(0,y) = w(0,y) \text{ for $y \geq 0$},$$
  $$\|w\|_{2,\sigma,\overline{\mathbb R_{+} \times \mathbb R_{+}}} \leq C[\|u\|_{2,\sigma,\overline{\{(x,0)\}}_{x\geq 0}} + \|u\|_{2,\sigma,\overline{\{(0,y)\}}_{y\geq 0}}].$$
\end{thm}
\emph{Proof:}
Define
$$w(x,y) = u(0,y)\phi(x)+u(x,0)\phi(y)-u(0,0)\phi(x)\phi(y).$$
Then we have
$$w(x,0) = u(0,0)\phi(x)+u(x,0)-u(0,0)\phi(x) = u(x,0),$$
$$w(0,y) = u(0,y)+u(0,0)\phi(y)-u(0,0)\phi(y) = u(0,y),$$
and $w \in C_{c}^{2,\sigma}(\overline{\mathbb R_{+} \times \mathbb R_{+}})$ such that
$$\|w\|_{2,\sigma,\overline{\mathbb R_{+} \times \mathbb R_{+}}} \leq C(\phi)[ \|u\|_{2,\sigma,\overline{\{(x,0)\}}_{x\geq 0}} + \|u\|_{2,\sigma,\overline{\{(0,y)\}}_{y\geq 0}}].$$
$\hfill$ q.e.d.\\

Remark: We use $u(0,y)\phi(x)$ instead of $u(0,y)$ because $u(0,y)$ is not compactly supported.

\begin{thm}
For any real number $\beta$, there is a constant $C = C(\beta)$ such that for any $u \in C_{c}^{2,\sigma}(\overline{\mathbb R_{+} \times \mathbb R_{+}})$, there exists a $w \in C_{c}^{2,\sigma}(\overline{\mathbb R_{+} \times \mathbb R_{+}})$ with
  $$\left[-\dfrac{\partial w}{\partial y} + \beta\dfrac{\partial w}{\partial x}\right]_{(x,0)} = g(x) \text{ for $x \geq 0$,}$$
  $$w(0,y) = h(y) \text{ for $y \geq 0$},$$
  $$\|w\|_{2,\sigma,\overline{\mathbb R_{+} \times \mathbb R_{+}}} \leq C[\|g\|_{1,\sigma,\overline{\mathbb R_{+}}} + \|h\|_{2,\sigma,\overline{\mathbb R_{+}}}]$$
  where
  $$g(x) = \left[-\dfrac{\partial u}{\partial y} + \beta\dfrac{\partial u}{\partial x}\right]_{(x,0)},$$
  $$h(y) = u(0,y).$$
\end{thm}
\emph{Proof:}
Let $\tilde{u} = u - h(y)\phi(x)$, then $\tilde{u}$ satisfies
\begin{align*}
\left[-\dfrac{\partial \tilde{u}}{\partial y} + \beta\dfrac{\partial \tilde{u}}{\partial x}\right]_{(x,0)} &= g(x)+h^{'}(0)\phi(x)-\beta h(0)\phi^{'}(x) &\text{ for $x \geq 0$}, \\
\tilde{u}(0,y) &= 0 &\text{ for $y \geq 0$}.
\end{align*}
Thus we reduce our problem to the case
\begin{align*}
\left[-\dfrac{\partial u}{\partial y} + \beta\dfrac{\partial u}{\partial x}\right]_{(x,0)} &= g(x) &\text{ for $x \geq 0$}, \\
u(0,y) &= 0 &\text{ for $y \geq 0$}.
\end{align*}
Case 1: $\beta = 0$. Then $g(0) = -\partial_{y}u(0,0) = 0$. Define\\
$$\widetilde{g}(x) = \begin{cases}
  g(x)-g^{'}(0)\phi(x)x             &\text{if $x \geq 0$}, \\
  0                                             &\text{if $x < 0$}.
\end{cases}$$\\
Then $\widetilde{g} \in C_{c}^{1,\sigma}(\mathbb R)$. Define\\
$$w(x,y) = T_{X}(\widetilde{g})(x,y)-\sum_{j=1}^{3}c_{j}T_{X}(\widetilde{g})(-jx,y) - g^{'}(0)\phi(x)\phi(y)xy.$$
Now we check if $w(x,y)$ satisfies the conditions we need. \\
Clearly $w \in C_{c}^{2,\sigma}(\overline{\mathbb R_{+} \times \mathbb R_{+}})$ and
\begin{align*}
  \left[-\dfrac{\partial w}{\partial y} + \beta\dfrac{\partial w}{\partial x}\right]_{(x,0)}
         &= \left[-\dfrac{\partial w}{\partial y} + 0\right]_{(x,0)}\\
         &= \tilde{g}(x) - \sum_{j=1}^{3}c_{j}(-j)\tilde{g}(-jx)+g^{'}(0)\phi(x)x \\
         &= g(x)-g(0)\phi(x)-g^{'}(0)\phi(x)x + 0 + g^{'}(0)\phi(x)x \\
         &= g(x) \text{    for $x\geq0$}, \\
  w(0,y) &= T_{X}(\widetilde{g})(0,y)-\sum_{j=1}^{3}c_{j}T_{X}(\widetilde{g})(0,y) - 0 \\
         &= T_{X}(\widetilde{g})(0,y)-T_{X}(\widetilde{g})(0,y) \\
         &= 0
\end{align*}
and we have the estimate of $w$,
$$\|w\|_{2,\sigma,\overline{\mathbb R_{+} \times \mathbb R_{+}}}
  \leq C(\|T_{X}\|,c_{j},\phi)\|g\|_{1,\sigma,\overline{\mathbb R_{+} }}.$$
Case 2: $\beta \neq 0$. Define\\
$\widetilde{g}(x) = \begin{cases}
  g(x)- \phi^{'}(x)[g(0)x+\dfrac{g^{'}(0)}{2}x^2] - \phi(x)[g(0)+g^{'}(0)x]           &\text{if $x \geq 0$}, \\
  0                                             &\text{if $x < 0$}.
\end{cases}$\\
Then $\widetilde{g} \in C_{c}^{1,\sigma}(\mathbb R)$. Define\\
$$w(x,y) = T_{X}(\widetilde{g})(x,y)-\sum_{j=1}^{3}c_{j}T_{X}(\widetilde{g})(-jx,y) + \dfrac{\phi(x)\phi(y)}{\beta}[g(0)x+\dfrac{g^{'}(0)}{2}x^2].$$
Similar to Case 1, we check if $w(x,y)$ satisfies the conditions we need. \\
Clearly $w \in C_{c}^{2,\sigma}(\overline{\mathbb R_{+} \times \mathbb R_{+}})$ and
\begin{align*}
  \left[-\dfrac{\partial w}{\partial y} + \beta\dfrac{\partial w}{\partial x}\right]_{(x,0)}
         &= \tilde{g}(x) - \sum_{j=1}^{3}c_{j}(-j)\tilde{g}(-jx) \\
         &\quad + \phi^{'}(x)[g(0)x+\dfrac{g^{'}(0)}{2}x^2] + \phi(x)[g(0)+g^{'}(0)x]   \\
         &= \tilde{g}(x) - 0 \\
         &\quad + \phi^{'}(x)[g(0)x+\dfrac{g^{'}(0)}{2}x^2] + \phi(x)[g(0)+g^{'}(0)x]   \\
         &=g(x) \text{  for $x\geq0$}, \\
  w(0,y) &= T_{X}(\widetilde{g})(0,y)-\sum_{j=1}^{3}c_{j}T_{X}(\widetilde{g})(0,y) + 0\\
         &= T_{X}(\widetilde{g})(0,y)- T_{X}(\widetilde{g})(0,y) \\
         &= 0
\end{align*}
and we have the estimate of $w$,
$$\|w\|_{2,\sigma,\overline{\mathbb R_{+} \times \mathbb R_{+}}}
  \leq C(T_{X}, \phi, \frac{1}{\beta}, c_{j} )\|g\|_{1,\sigma,\overline{\mathbb R_{+}}}.$$
$\hfill$  q.e.d.\\

\begin{thm}
For any real numbers $\alpha$ and $\beta$, there is a constant \\
$C = C(\alpha, \beta)$ such that for any $u \in C_{c}^{2,\sigma}(\overline{\mathbb R_{+} \times \mathbb R_{+}})$, there exists a $w \in C_{c}^{2,\sigma}(\overline{\mathbb R_{+} \times \mathbb R_{+}})$ with
  $$\left[-\dfrac{\partial w}{\partial y} + \beta\dfrac{\partial w}{\partial x}\right]_{(x,0)} = g(x) \text{ for $x \geq 0$},$$
  $$\left[-\dfrac{\partial w}{\partial x} + \alpha\dfrac{\partial w}{\partial y}\right]_{(0,y)} = h(y) \text{ for $y \geq 0$},$$
  $$\|w\|_{2,\sigma,\overline{\mathbb R_{+} \times \mathbb R_{+}}} \leq C[\|g\|_{1,\sigma,\overline{\mathbb R_{+}}} + \|h\|_{1,\sigma,\overline{\mathbb R_{+}}}]$$
  where
  $$g(x) = \left[-\dfrac{\partial u}{\partial y} + \beta\dfrac{\partial u}{\partial x}\right]_{(x,0)},$$
  $$h(y) = \left[-\dfrac{\partial u}{\partial x} + \alpha\dfrac{\partial u}{\partial y}\right]_{(0,y)}.$$
\end{thm}
 \emph{Proof:} Define
$\widetilde{h}(y) = \begin{cases}
  h(y)                      &\text{if $y \geq 0$}, \\
  \sum_{j=1}^{3}c_{j}h(-jy) &\text{if $y < 0$}.
\end{cases}$\\
Let $v = T_{Y}[\widetilde{h}(y)]$ and $u_{1} = u-v$, then $u_{1} \in C_{c}^{2,\sigma}(\overline{\mathbb R_{+} \times \mathbb R_{+}})$ satisfies
$$\left[-\dfrac{\partial u_{1}}{\partial y} + \beta\dfrac{\partial u_{1}}{\partial x}\right]_{(x,0)} = g_{1}(x) \text{ for $x \geq 0$},$$
$$\left[-\dfrac{\partial u_{1}}{\partial x} + \alpha\dfrac{\partial u_{1}}{\partial y}\right]_{(0,y)} = 0 \quad \,\,\, \text{ for $y \geq 0$},$$
where $g_{1}(x) = g(x) - [-\partial_{y}v+\beta\partial_{x}v]_{(x,0)}$. Thus we reduce our problem to the case
$$\left[-\dfrac{\partial u}{\partial y} + \beta\dfrac{\partial u}{\partial x}\right]_{(x,0)} = g(x) \text{ for $x \geq 0$},$$
$$\left[-\dfrac{\partial u}{\partial x} + \alpha\dfrac{\partial u}{\partial y}\right]_{(0,y)} = 0 \quad \,\,\, \text{ for $y \geq 0$},$$
Case 1: $\alpha\beta = 1$.\\
As $\alpha\beta = 1$, we have
$$g(0) = (-\partial_{y}u+\beta\partial_{x}u)(0,0) = -\beta(-\partial_{x}u+\alpha\partial_{y}u)(0,0) = 0.$$
Define
$$\widetilde{g}(x) = \begin{cases}
  g(x)- \dfrac{g^{'}(0)}{2}\phi^{'}(x)x^2 - g^{'}(0)\phi(x)x         &\text{if $x \geq 0$}, \\
  0                                             &\text{if $x < 0$}.
\end{cases}$$
Then $\widetilde{g} \in C_{c}^{1,\sigma}(\mathbb R)$. Define
$$w(x,y) = T_{X}(\widetilde{g})(x,y) - \sum_{j=1}^{3}c_{j}T_{X}(\widetilde{g})(-jx,y) + \dfrac{g^{'}(0)}{2\beta}\phi(x)\phi(y)x^2.$$
Case 2.1: $\alpha\beta \neq 1$, $\beta \neq 0$. Define
$$\widetilde{g}(x) = \begin{cases}
  g(x)- g(0)\phi(x) - \dfrac{g^{'}(0)}{2}\phi^{'}(x)x^2 - g^{'}(0)\phi(x)x            &\text{if $x \geq 0$}, \\
  0                                             &\text{if $x < 0$}.
\end{cases}$$
Then $\widetilde{g} \in C_{c}^{1,\sigma}(\mathbb R)$. Define
\begin{align*}
w(x,y) &= T_{X}(\widetilde{g})(x,y) - \sum_{j=1}^{3}c_{j}T_{X}(\widetilde{g})(-jx,y) + \dfrac{g^{'}(0)}{2\beta}\phi(x)\phi(y)x^2 \\
       &\quad -\dfrac{\alpha g(0)}{1-\alpha\beta}\phi(x)\phi(y)x -  \dfrac{g(0)}{1-\alpha\beta}\phi(x)\phi(y)y.
\end{align*}
Case 2.2.1: $\alpha\beta \neq 1$, $\beta = 0$, $\alpha \neq 0$. Define
$$\widetilde{g}(x) = \begin{cases}
  g(x)- g(0)\phi(x) - g^{'}(0)\phi(x)x            &\text{if $x \geq 0$}, \\
  0                                             &\text{if $x < 0$}.
\end{cases}$$
Then $\widetilde{g} \in C_{c}^{1,\sigma}(\mathbb R)$. Define
\begin{align*}
w(x,y) &= T_{X}(\widetilde{g})(x,y) - \sum_{j=1}^{3}c_{j}T_{X}(\widetilde{g})(-jx,y)\\
       &\quad -\alpha g(0)\phi(x)\phi(y)x -  g(0)\phi(x)\phi(y)y\\
       &\quad - g^{'}(0)\phi(x)\phi(y)xy - \dfrac{g^{'}(0)}{2\alpha}\phi(x)\phi(y)y^2.
\end{align*}
Case 2.2.2: $\beta = 0$, $\alpha = 0$.\\
As $\alpha = \beta = 0$, we have
$$(\partial_{x}u)(0,y) = 0,-(\partial_{y}u)(x,0) = g(x).$$
So $g^{'}(0) = -(\partial_{x}\partial_{y}u)(0,0) = 0$. Define
$$\widetilde{g}(x) = \begin{cases}
  g(x)- g(0)\phi(x)           &\text{if $x \geq 0$}, \\
  0                                             &\text{if $x < 0$}.
\end{cases}$$
Then $\widetilde{g} \in C_{c}^{1,\sigma}(\mathbb R)$. Define
\begin{align*}
w(x,y) &= T_{X}(\widetilde{g})(x,y) - \sum_{j=1}^{3}c_{j}T_{X}(\widetilde{g})(-jx,y)\\
       &\quad -  g(0)\phi(x)\phi(y)y.
\end{align*}
We can check that $w(x,y)$ we defined above satisfies the respect boundary conditions by direct computations and the norm of $w$ satisfies
 $$ \|w+T_{Y}[\widetilde{h}(y)]\|_{2,\sigma,\overline{\mathbb R_{+} \times \mathbb R_{+}}}
  \leq C(c_{j}, T_{X}, T_{Y}, \phi, \alpha, \beta )[\|h\|_{1,\sigma,\overline{\mathbb R_{+}}} + \|g\|_{1,\sigma,\overline{\mathbb R_{+}}}].$$
$\hfill$  q.e.d.\\

\subsubsection{On $\mathbb R_{+}^{2}$}
Let $C_{c}^{k,\sigma}(\overline{\mathbb R_{+}^{2}}) = \{g \in C^{k,\sigma}(\overline{\mathbb R_{+}^{2}}) \, \big| \, \text{supp}(g)$ is compact$\}$. We have the following theorems.
\begin{thm}
There is a constant $C$ such that for any $u \in C_{c}^{2,\sigma}(\overline{\mathbb R_{+}^{2}})$, there exists a $w \in C_{c}^{2,\sigma}(\overline{\mathbb R_{+}^{2}})$ with
  $$w(x,0) = u(x,0) \text{ for $x \geq 0$},$$
  $$w(x,0) = u(x,0) \text{ for $x \leq 0$},$$
  $$\|w\|_{2,\sigma,\overline{\mathbb R_{+}^{2}}} \leq C[\|u\|_{2,\sigma,\overline{\{(x,0)\}}_{x\geq 0}} + \|u\|_{2,\sigma,\overline{\{(x,0)\}}_{x\leq 0}}].$$
\end{thm}
\emph{Proof:}
Define
$$w(x,y) = u(x,0)\phi(y).$$
Then $w$ is what we need.
$\hfill$ q.e.d.\\

\begin{thm}
For any real number $\beta$, there is a constant $C = C(\beta)$ such that for any $u \in C_{c}^{2,\sigma}(\overline{\mathbb R_{+}^{2}})$, there exists a $w \in C_{c}^{2,\sigma}(\overline{\mathbb R_{+}^{2}})$ with
  $$\left[-\dfrac{\partial w}{\partial y} + \beta\dfrac{\partial w}{\partial x}\right]_{(x,0)} = g(x) \text{ for $x \geq 0$},$$
  $$w(x,0) = h(x) \text{ for $x \leq 0$},$$
  $$\|w\|_{2,\sigma,\overline{\mathbb R_{+}^{2}}} \leq C[\|g\|_{1,\sigma,\overline{\{(x,0)\}}_{x\geq 0}} + \|h\|_{2,\sigma,\overline{\{(x,0)\}}_{x\leq 0}}]$$
  where
  $$g(x) = \left[-\dfrac{\partial u}{\partial y} + \beta\dfrac{\partial u}{\partial x}\right]_{(x,0)} \text{ for $x \geq 0$},$$
  $$h(x) = u(x,0) \text{ for $x \leq 0$}.$$
\end{thm}
\emph{Proof:}
Let $c_{j}$ satisfy $\sum_{j=1}^{3}(-j)^{l}c_{j} = 1$ for any $l=0,1,2$. Define
$$\widetilde{g}(x) = \begin{cases}
  g(x)                      &\text{if $x \geq 0$}, \\
  \sum_{j=1}^{3}c_{j}g(-jx) &\text{if $x < 0$}.
\end{cases}$$
$$\widetilde{h}(x) = \begin{cases}
  \sum_{j=1}^{3}c_{j}h(-jx) &\text{if $x > 0$}, \\
  h(x)                      &\text{if $x \leq 0$}. \\
\end{cases}$$
and
$$w(x,y) = \tilde{h}(x)\phi(y) + T_{X}(\tilde{g}-\beta\partial_{x}\tilde{h})(x,y).$$
Then we have
\begin{align*}
\left[-\dfrac{\partial w}{\partial y} + \beta\dfrac{\partial w}{\partial x}\right]_{(x,0)} &= (0+g(x)-\beta\partial_{x}\tilde{h}(x))+\beta(\partial_{x}\tilde{h}(x)+0) \\
&= g(x) \text{$\quad$ for $x \geq 0$}, \\
w(x,0) &= h(x)+0 = h(x) \text{$\quad$ for $x \leq 0$}
\end{align*}
and
$$\|w\|_{2,\sigma,\overline{\mathbb R_{+}^{2}}} \leq C(T_{X},c_{j},\beta, \phi)[\|g\|_{1,\sigma,\overline{\{(x,0)\}}_{x\geq 0}} + \|h\|_{2,\sigma,\overline{\{(x,0)\}}_{x\leq 0}}].$$
$\hfill$ q.e.d\\

\begin{thm}
For any real numbers $\alpha$ and $\beta$, there is a constant \\
$C = C(\alpha, \beta)$ such that for any $u \in C_{c}^{2,\sigma}(\overline{\mathbb R_{+}^{2}})$, there exists a $w \in C_{c}^{2,\sigma}(\overline{\mathbb R_{+}^{2}})$ with
  $$\left[-\dfrac{\partial w}{\partial y} + \beta\dfrac{\partial w}{\partial x}\right]_{(x,0)} = g(x) \text{ for $x \geq 0$},$$
  $$\left[-\dfrac{\partial w}{\partial y} + \alpha\dfrac{\partial w}{\partial x}\right]_{(x,0)} = h(x) \text{ for $x \leq 0$},$$
  $$\|w\|_{2,\sigma,\overline{\mathbb R_{+}^{2}}} \leq C(\alpha,\beta)[\|g\|_{1,\sigma,\overline{\{(x,0)\}}_{x\geq 0}} + \|h\|_{1,\sigma,\overline{\{(x,0)\}}_{x\leq 0}}]$$
  where
  $$g(x) = \left[-\dfrac{\partial u}{\partial y} + \beta\dfrac{\partial u}{\partial x}\right]_{(x,0)} \text{$\quad$ for $x \geq 0$},$$
  $$h(x) = \left[-\dfrac{\partial u}{\partial y} + \alpha\dfrac{\partial u}{\partial x}\right]_{(x,0)} \text{$\quad$ for $x \leq 0$}.$$
\end{thm}
\emph{Proof:}
Case 1: $\alpha = \beta$. Define
$$w(x,y) = T_{X}[(-\partial_{y}u+\alpha\partial_{x}u)(x,0)]$$
We can check that $w(x,y)$ is what we need.\\
Case 2: $\alpha \neq \beta$. Let $c_{j}$ satisfy $\sum_{j=1}^{3}(-j)^{l}c_{j} = 1$ for any $l=0,1,2$. Define
$$\widetilde{h}(x) = \begin{cases}
  \sum_{j=1}^{3}c_{j}h(-jx) &\text{if $x > 0$}, \\
  h(x)                      &\text{if $x \leq 0$}.
\end{cases}$$
Let $\tilde{u}(x,y) = u - T_{X}(\tilde{h})$, then $\tilde{u}$ satisfies
$$\left[-\dfrac{\partial \tilde{u}}{\partial y} + \beta\dfrac{\partial \tilde{u}}{\partial x}\right]_{(x,0)} = g_{1}(x) \text{ for $x \geq 0$},$$
$$\left[-\dfrac{\partial \tilde{u}}{\partial y} + \alpha\dfrac{\partial \tilde{u}}{\partial x}\right]_{(x,0)} = 0 \quad \,\,\,\, \text{ for $x \leq 0$}$$
where $g_{1}(x) = g(x) - \tilde{h}(x)$ for $x \geq 0$. So we reduce our problem to the case
$$\left[-\dfrac{\partial u}{\partial y} + \beta\dfrac{\partial u}{\partial x}\right]_{(x,0)} = g(x) \text{ for $x \geq 0$},$$
$$\left[-\dfrac{\partial u}{\partial y} + \alpha\dfrac{\partial u}{\partial x}\right]_{(x,0)} = 0 \quad \,\,\, \text{ for $x \leq 0$}.$$
Define
$$\widetilde{g}(x) = \begin{cases}
  g(x)-\phi(x)[g(0)+g^{'}(0)x] - \dfrac{\beta \phi^{'}(x)}{\beta - \alpha}[g(0)x+\dfrac{g^{'}(0)}{2}x^2]                      &\text{if $x \geq 0$}, \\
  -\dfrac{\alpha \phi^{'}(x)}{\beta - \alpha}[g(0)x+\dfrac{g^{'}(0)}{2}x^2] &\text{if $x < 0$}.
\end{cases}$$
We have $\widetilde{g} \in C_{c}^{1,\sigma}(\mathbb R)$ since $\phi = 1$ near $x=0$. Define
$$w(x,y) = T_{X}(\tilde{g}) + \dfrac{\phi(x)\phi(y)}{\beta - \alpha}[g(0)x + \dfrac{g^{'}(0)}{2}x^2] + \dfrac{\alpha\phi(x)\phi(y)}{\beta - \alpha}[g(0) + g^{'}(0)x]y.$$
Then $w$ is what we need.
$\hfill$ q.e.d.\\

\subsubsection{On $\mathbb R^2 \setminus \mathbb R_{+} \times \mathbb R_{+}$}
Let $C_{c}^{k,\sigma}(\overline{\mathbb R^2 \setminus \mathbb R_{+} \times \mathbb R_{+}}) = \{g \in C^{k,\sigma}(\overline{\mathbb R^2 \setminus \mathbb R_{+} \times \mathbb R_{+}}) \, \big| \, \text{supp}(g)$ is compact$\}$ and $P$ be the extension operator in Theorem 2.8.1.

For any $u \in C_{c}^{2,\sigma}(\overline{\mathbb R^2 \setminus \mathbb R_{+} \times \mathbb R_{+}})$, let $\tilde{u} = P(u)\big|_{\mathbb R_{+} \times \mathbb R_{+}}$. By the results in Section 2.8.1, there exists a $\tilde{w} \in C_{c}^{2,\sigma}(\overline{\mathbb R_{+} \times \mathbb R_{+}})$ such that $\tilde{w}$ satisfies the respect Dirichlet or oblique boundary conditions of $\tilde{u}$ and the norm of $\tilde{w}$ is controlled by the boundary values of $\tilde{u}$. Then we extend $\tilde{w}$ to a function $w \in C_{c}^{2,\sigma}(\mathbb R^2)$. Since $\partial[\mathbb R^2 \setminus \mathbb R_{+} \times \mathbb R_{+}] = \partial[\mathbb R_{+} \times \mathbb R_{+}]$, $w\big|_{\mathbb R^2 \setminus \mathbb R_{+} \times \mathbb R_{+}}$ satisfies the respect Dirichlet or oblique boundary conditions of $u$ and the norm of $w$ is controlled by the respect norm of the boundary values of $u$.

\subsubsection{On sectors}
For any sector $T_{\theta_{0}, \omega_{0}}=\{(r\cos \theta, r\sin \theta): r > 0, \theta_{0} < \theta < \theta_{0}+\omega_{0}\}$ with $\omega_{0} \in (0,2\pi)$, let $C_{c}^{k,\sigma}(\overline{T_{\theta_{0}, \omega_{0}}}) = \{g \in C^{k,\sigma}(\overline{T_{\theta_{0}, \omega_{0}}}) \, \big| \, \text{supp}(g)$ is compact$\}$.

We can deform $T_{\theta_{0}, \omega_{0}}$ to the first quadrant $\mathbb R_{+} \times \mathbb R_{+}$, the upper half space $\mathbb R_{+}^{2}$ or $\mathbb R^2 \setminus \mathbb R_{+} \times \mathbb R_{+}$ under rotations and dilations. And $\pm \partial \nu + \beta^{'}\partial \tau$ will be $-\partial y + \beta \partial x \big|_{(x,0)}$ where $\nu$ is the outer unit normal vector of $T_{\theta_{0}, \omega_{0}}$ on $\partial T_{\theta_{0}, \omega_{0}}$, $\tau$ is the unit tangent vector of $\partial T_{\theta_{0}, \omega_{0}}$. Thus the previous inverse trace theorems still hold for $u \in C_{c}^{2,\sigma}(\overline{T_{\theta_{0}, \omega_{0}}})$ if we consider the Dirichlet or oblique boundary conditions on $\partial T_{\theta_{0}, \omega_{0}}$. And the constant $C$ will be $C(\alpha, \beta, \theta_{0}, \omega_{0})$.

\subsubsection{On polygonal domains}
Based on the previous results, we have the following inverse trace theorems on polygonal domains.
\begin{thm}
  Let $\Omega$ be a polygonal domain in $\mathbb R^2$ and $\omega_{j} \in (0,2\pi)$. For any real numbers $\beta_{j}$, there is a constant $C = C(\Omega, \beta_{j})$ such that for any $u \in C^{2,\sigma}(\overline{\Omega})$, there exists a $w \in C^{2,\sigma}(\overline{\Omega})$ with
  $$\begin{cases}
     w = u         & \text{on $\Gamma_{j}$, $j \in \mathcal D$}, \\
     \dfrac{\partial w}{\partial \nu_{j}} + \beta_{j}\dfrac{\partial w}{\partial \tau_{j}}=\dfrac{\partial u}{\partial \nu_{j}} + \beta_{j}\dfrac{\partial u}{\partial \tau_{j}} & \text{on $\Gamma_{j}$, $j \in \mathcal N$},
  \end{cases}$$
  $$\|w\|_{2,\sigma,\overline{\Omega}} \leq C[\sum_{j \in \mathcal D}\|u\|_{2,\sigma,\overline{\Gamma_{j}}} +  \sum_{j \in \mathcal N}\|\dfrac{\partial u}{\partial \nu_{j}} + \beta_{j}\dfrac{\partial u}{\partial \tau_{j}}\|_{1,\sigma,\overline{\Gamma_{j}}}].$$
\end{thm}
\emph{Proof:}
Let $\eta_{j},\phi_{j}$ be smooth functions such that:\\
(1) $\eta_{j}(x),\phi_{j}(x) \in [0,1]$ for any $x \in \mathbb R^2$, $\sum\eta_{j}(x)+\sum \phi_{j}(x) = 1$ for any $x \in \Gamma$.\\
(2) $\eta_{j} = 1$ near $S_{j}$, $\text{supp}(\eta_{j}) \cap \Gamma_{l} = \emptyset$ for any $l \notin \{j,j+1\}$.\\
(3) $\text{supp}(\phi_{j}) \cap \Gamma_{l} = \emptyset$ for any $l \neq j$.\\
(4) $\dfrac{\partial \eta_{j}}{\partial \nu_{k}} + \beta_{k}\dfrac{\partial \eta_{j}}{\partial \tau_{k}}= \dfrac{\partial \phi_{j}}{\partial \nu_{k}} +
\beta_{k}\dfrac{\partial \phi_{j}}{\partial \tau_{k}} = 0$ if $k \in \mathcal N$.\\
Let $P$ be the extension operator in Theorem 2.8.1. \\
By the results in Section 2.8.4 and a translation, there are functions $w_{j} \in C_{c}^{2,\sigma}(\mathbb R^2)$ such that
$$\begin{cases}
     &w_{j} = P(u) = u        \\
     & \text{\quad \quad \quad on $\Gamma_{k}$, $k \in \{j,j+1\} \cap \mathcal D$}, \\
     &\dfrac{\partial w_{j}}{\partial \nu_{k}} + \beta_{k}\dfrac{\partial w_{j}}{\partial \tau_{k}}=\dfrac{\partial P(u)}{\partial \nu_{k}} + \beta_{k}\dfrac{\partial P(u)}{\partial \tau_{k}}= \dfrac{\partial u}{\partial \nu_{k}} + \beta_{k}\dfrac{\partial u}{\partial \tau_{k}}\\
     &\text{\quad \quad \quad on $\Gamma_{k}$, $k \in \{j,j+1\} \cap \mathcal N$}.
  \end{cases}$$
  By Lemma 2.8.1, there are functions $\tilde{w}_{j} \in C_{c}^{2,\sigma}(\mathbb R^2)$ such that
  $$\begin{cases}
     \tilde{w}_{j} = P(u) = u         & \text{on $\Gamma_{j}$, if $j \in \mathcal D$}, \\
     \dfrac{\partial \tilde{w}_{j}}{\partial \nu_{j}} + \beta_{j}\dfrac{\partial \tilde{w}_{j}}{\partial \tau_{j}} = \dfrac{\partial P(u)}{\partial \nu_{j}} + \beta_{j}\dfrac{\partial P(u)}{\partial \tau_{j}} = \dfrac{\partial u}{\partial \nu_{j}} + \beta_{j}\dfrac{\partial u}{\partial \tau_{j}} & \text{on $\Gamma_{j}$, if $j \in \mathcal N$}.
  \end{cases}$$
Define
$$w = \sum \eta_{j}w_{j} + \sum \phi_{j}\tilde{w}_{j}.$$
We claim $w$ is what we need.\\
If $k \in \mathcal D$, for any $x \in \Gamma_{k}$,
\begin{align*}
  w(x) &= \eta_{k-1}(x)w_{k-1}(x) + \eta_{k}(x)w_{k}(x) + \phi_{k}(x)\tilde{w}_{k}(x)\\
       &= \eta_{k-1}(x)u(x) + \eta_{k}(x)u(x) + \phi_{k}(x)u(x)\\
       &= \left[\sum \eta_{j}(x) + \sum \phi_{j}(x)\right]u(x) \\
       &= u(x).
\end{align*}
If $k \in \mathcal N$, for any $x \in \Gamma_{k}$,
\begin{align*}
  \left[\dfrac{\partial w}{\partial \nu_{j}} + \beta_{j}\dfrac{\partial w}{\partial \tau_{j}}\right](x)
  &= \eta_{k-1}(x)\left[\dfrac{\partial w_{k-1}}{\partial \nu_{k}} + \beta_{k}\dfrac{\partial w_{k-1}}{\partial \tau_{k}}\right](x) + \\
  &\quad \eta_{k}(x)\left[\dfrac{\partial w_{k}}{\partial \nu_{k}} + \beta_{k}\dfrac{\partial w_{k}}{\partial \tau_{k}}\right](x) \\
  &\quad + \phi_{k}(x)\left[\dfrac{\partial \tilde{w}_{k}}{\partial \nu_{k}} + \beta_{k}\dfrac{\partial \tilde{w}_{k}}{\partial \tau_{k}}\right](x)\\
  &= \eta_{k-1}(x)\left[\dfrac{\partial u}{\partial \nu_{k}} + \beta_{k}\dfrac{\partial u}{\partial \tau_{k}}\right](x) \\
  &\quad + \eta_{k}(x)\left[\dfrac{\partial u}{\partial \nu_{k}} + \beta_{k}\dfrac{\partial u}{\partial \tau_{k}}\right](x) \\
    &\quad + \phi_{k}(x)\left[\dfrac{\partial u}{\partial \nu_{k}} + \beta_{k}\dfrac{\partial u}{\partial \tau_{k}}\right](x)\\
  &= \left[\sum \eta_{j}(x) + \sum \phi_{j}(x)\right]\left[\dfrac{\partial u}{\partial \nu_{k}} + \beta_{k}\dfrac{\partial u}{\partial \tau_{k}}\right](x)\\
  &= \left[\dfrac{\partial u}{\partial \nu_{k}} + \beta_{k}\dfrac{\partial u}{\partial \tau_{k}}\right](x).
\end{align*}
And we have
$$\|w\|_{2,\sigma,\overline{\Omega}} \leq C(\eta_{j}, \phi_{j}, P, \beta_{j}, \Omega)[\sum_{j \in \mathcal D}\|u\|_{2,\sigma,\overline{\Gamma_{j}}} +  \sum_{j \in \mathcal N}\|\dfrac{\partial u}{\partial \nu_{j}} + \beta_{j}\dfrac{\partial u}{\partial \tau_{j}}\|_{1,\sigma,\overline{\Gamma_{j}}}].$$
$\hfill$ q.e.d.

\subsection{Deformations}
\begin{thm}
  Let $B$ be a symmetric $2 \times 2$ matrix with eigenvalues $\lambda_{1}, \lambda_{2} \in (0,1]$. Let $\epsilon$ be a small positive number. Then there is a diffeomorphism $\Phi : \mathbb R^2 \to \mathbb R^2$ such that\\
  (1) $\Phi(x,y) = B\begin{bmatrix}
    x\\
    y
  \end{bmatrix}$ on $B_{\epsilon}(0,0)$.\\
  (2) $\Phi(x,y) = (x,y)$ outside $B_{2\epsilon}(0,0)$.
\end{thm}
\emph{Proof:} Let $\eta:[0,\infty) \to [0,1]$ be a non-increasing cut-off function on $\overline{\mathbb R_{+}}$ such that $\eta = 1$ on $[0,\epsilon]$ and $\eta = 0$ on $[2\epsilon, \infty)$.\\
Since $B$ is a symmetric $2 \times 2$ matrix with eigenvalues $\lambda_{1}, \lambda_{2} \in (0,1)$,
$B = \begin{bmatrix}
  \lambda_{1} & 0\\
  0           & \lambda_{2}
\end{bmatrix}$ in the coordinate system spanned by the eigenvectors of $B$. Thus we can assume that $B = \begin{bmatrix}
  \lambda_{1} & 0\\
  0           & \lambda_{2}
\end{bmatrix}$ and we define
$$\Phi(x,y) = \eta(\sqrt{x^2+y^2})\begin{bmatrix}
    \lambda_{1}x\\
    \lambda_{2}y
  \end{bmatrix}
   + [1-\eta(\sqrt{x^2+y^2})]\begin{bmatrix}
    x\\
    y
  \end{bmatrix}.$$
Obviously $\Phi$ satisfies the conditions (1) and (2). We only need to check that $\Phi$ is globally bijective and locally the Jacobian is non-zero.\\
Since $\lambda_{1}, \lambda_{2} \in (0,1]$ and $\eta \in [0,1]$, $\Phi$ is bijective.\\
$|d\Phi|$\\
$=\left| \begin{matrix}
    [(\lambda_{1}-1)\eta+1]+\frac{\eta^{'}}{\sqrt{x^2+y^2}}(\lambda_{1}-1)x^2    &    \frac{(\lambda_{1}-1)xy}{\sqrt{x^2+y^2}}\eta^{'} \\
    \frac{(\lambda_{2}-1)xy}{\sqrt{x^2+y^2}}\eta^{'}   &  [(\lambda_{2}-1)\eta+1]+\frac{\eta^{'}}{\sqrt{x^2+y^2}}(\lambda_{2}-1)y^2
   \end{matrix} \right|$\\
$= [(\lambda_{1}-1)\eta+1][(\lambda_{2}-1)\eta+1] \\
  + \frac{\eta^{'}}{\sqrt{x^2+y^2}}\left\{(\lambda_{1}-1)[(\lambda_{2}-1)\eta+1]x^2 + (\lambda_{2}-1)[(\lambda_{1}-1)\eta+1]y^2  \right\}\\
  \geq \lambda_{1} \lambda_{2} + 0 \hfill \text{ since } ( \eta^{'}, (\lambda_{1}-1), (\lambda_{2}-1) \leq0)$ \\
$> 0$.
$\hfill$ q.e.d.

Remark: The reason we assume that the eigenvalues are in $(0,1]$ is that it will guarantee that $\Phi$ is ``increasing'' to the identity operator.

Next we prove that there is a deformation which only deforms a polygonal domain near its corners.
\begin{thm}
  Let $\Omega$ be a polygonal domain, $\epsilon$ be a small positive number and $B_{j}$ be symmetric $2 \times 2$ matrices with eigenvalues in $(0,1]$. Then there is a diffeomorphism $\Phi : \mathbb R^2 \to \mathbb R^2$ such that\\
  (1) $\Phi(x,y) = S_{j} +  B_{j}\left\{\begin{bmatrix}
    x\\
    y
  \end{bmatrix}-S_{j}\right\}$ on $B_{\epsilon}(S_{j})$.\\
  (2) $\Phi(x,y) = (x,y)$ outside $\cup_{j} B_{2\epsilon}(S_{j})$.
\end{thm}
\emph{Proof:} Similar to the previous theorem, let $\eta:[0,\infty) \to [0,1]$ be a non-increasing cut-off function on $\overline{\mathbb R_{+}}$ such that $\eta = 1$ on $[0,\epsilon]$ and $\eta = 0$ on $[2\epsilon, \infty)$. We define
\begin{align*}
\Phi(x,y)
&=\sum_{j=1}^{I_{0}} \left( S_{j} + \eta(|(x,y)-S_{j}|)B_{j}\left\{\begin{bmatrix}
    x\\
    y
  \end{bmatrix}-S_{j}\right\}\right)\\
&\quad  +\left(1-\sum_{j=1}^{I_{0}} \eta(|(x,y)-S_{j}|)\right)\begin{bmatrix}
    x\\
    y
  \end{bmatrix}.
\end{align*}
Since $\epsilon$ can be arbitrarily small, $\Phi$ satisfies the conditions (1) and (2).\\
The same method as in the previous theorem shows that $\Phi$ is globally bijective and locally the Jacobian is non-zero.
$\hfill$ q.e.d.

\section{The Mod-2 Degree Theory}
\subsection{Introduction}
In this section, we extend the degree theory developed by White\cite{White1987} to the partially free boundary setting. We first consider the embeddings from the unit half disk to $(\mathbb R_{+}^{3}, g)$ and then all theorems still hold for the embeddings from the unit half disk to $(\mathbb B_{+}^{3}, g)$. We prove that the space of properly embedded partially free boundary minimal half disks in $(\mathbb R_{+}^{3}, g)$ is a Banach manifold and the projection map which projects such surfaces into their Dirichlet boundaries has a well-defined mod-2 degree under suitable assumptions.

We prove the main results using similar arguments in White\cite{White1987}. The new major analytic difficulty is to analyze the
properties of the Jacobi operator with mixed boundary conditions on a domain with corners which we solved in Section 2.

\subsection{Notations}
\subsubsection{Domains}
$N_{+}$: Let $N_{+} = (\mathbb R_{+}^3, g)$ be a three-dimensional Riemannian manifold with weakly mean convex boundary and the metric $g$ be $C^{q,\beta}$ up to the boundary.

$D_{+}$: Let $D_{+} = \{(x,y): x^2+y^2 < 1, y>0\}$ be the unit half disk. And $\partial D_{+} = \overline{\Gamma_{1}} \cup \overline{\Gamma_{2}}$ where $\Gamma_{1} = \{(x,y):x^2+y^2=1, y>0\}$ and $\Gamma_{2} = \{(x,0):x \in (-1,1)\}$. The corners are $S_{1} = (-1,0)$, $S_{2} = (1,0)$.

$\nu_{j}, \tau_{j}$: $\nu_{j}$ is the outward unit normal vector of $D_{+}$ on $\Gamma_{j}$, $\tau_{j}$ is the unit tangent vector of $\Gamma_{j}$.

Remark: We use $N_{+}, D_{+}$ to denote the open space and $\overline{N_{+}}, \overline{D_{+}}$ for their closure. We can deform $D_{+}$ to a $C^{2,1}$ curvilinear polygonal domain with angles equal to $\dfrac{\pi}{2}$. For example, the map $\phi: \{(x,y): x \geq -1, y \geq 0\} \to D_{+}$ such that $\phi(x,y) = (x+1-\sqrt{1-y^2},y)$ near $S_{1} = (-1,0)$ satisfies our requirement.

\subsubsection{$\mathbf{C_{\perp}^{2,\sigma}(\overline{\Gamma_{1}}, N_{+})}$, $\mathbf{C_{*}^{2,\sigma}(\overline{D_{+}}, N_{+})}$}
Let $\sigma \in (0,\beta)$ where $\beta$ is given in Section 3.2.1.

$C_{\perp}^{2,\sigma}(\overline{\Gamma_{1}}, N_{+}) = \{\gamma :\Gamma_{1} \to N_{+} \, \big| \, \gamma $ is an embedding and $C^{2,\sigma}$ up to the boundary, $\partial [\gamma(\Gamma_{1})] = \overline{\gamma(\Gamma_{1})} \cap \partial N_{+}$, $\gamma(\Gamma_{1}) \perp \partial N_{+}$ with respect to the metric $g  \}$.

$C_{*}^{2,\sigma}(\overline{D_{+}}, N_{+}) = \{f:D_{+} \to N_{+} \, \big| \, f$ is an embedding and $C^{2,\sigma}$ up to the boundaries and corners, $\partial [f(D_{+})] \cap \partial N_{+} = \overline{f(\Gamma_{2})}$, $\overline{f(D_{+})}$ is transversal to $\partial N_{+}$, $f|_{\Gamma_{1}} \in C_{\perp}^{2,\sigma}(\overline{\Gamma_{1}}, N_{+}) \}$.

Two maps $f_{1}, f_{2} \in C_{*}^{2,\sigma}(\overline{D_{+}}, N_{+})$ are said to be equivalent if $f_{1} = f_{2} \circ u$ for some $C^{2,\sigma}$ diffeomorphism $u: D_{+} \to D_{+}$ such that $u(x) = x$ for all $x \in \Gamma_{1}$. We use $[f]$ to denote the equivalence class of $f$.

\subsubsection{The first variation}
Let $f \in C_{*}^{2,\sigma}(\overline{D_{+}}, N_{+})$ and $\vec{X}: \overline{D_{+}} \to T N_{+}$ such that $\vec{X}$ is $C^{2,\sigma}$ and $\vec{X}(x) \in T_{f(x)}(\partial N_{+})$ for each $x \in \Gamma_{2}$, then it generates a $C^{2,\sigma}$ one-parameter family of maps $F(x,t): D_{+} \times (-\epsilon, \epsilon) \to N_{+}$ such that
$$F(x,0) = f, F_{t}(x,0) = \vec{X}, F(\cdot,t) \in C_{*}^{2,\sigma}(\overline{D_{+}}, N_{+}) \text{ for each $t$}.$$

We have the following first variation formula of the area.
\begin{thm}
We denote $F_{i}= \partial_{x_{i}}F(x,t)$ and $\tilde{g}_{ij} = <F_{i},F_{j}>_{F(x,t)}$ for $i,j=1,2$. Then
$$\dfrac{d}{dt}\Big|_{t=0}\text{Area}(F(\cdot,t)) = -\int_{f(D_{+})}<\vec{H}_{f},\vec{X}> + \int_{f(\Gamma_{1})}<\vec{n},\vec{X}> + \int_{f(\Gamma_{2})}<\vec{n},\vec{X}>$$
where $\vec{H}_{f} = \tilde{g}^{ij}\nabla_{f_{i}}f_{j}$ is the mean curvature vector of $f(D_{+})$ in $N_{+}$, $\tilde{g}^{ij}$ is the inverse matrix of $\tilde{g}_{ij}$, $\nabla$ is the Riemannian connection of $N_{+}$ and $\vec{n}$ is the outward unit normal vector of $f(D_{+})$ on $f(\Gamma_{1})$ and $f(\Gamma_{2})$.
\end{thm}
\emph{Proof:} By definition
$$\text{Area}(F(\cdot,t)) = \int_{D_{+}}\sqrt{\det (\tilde{g}_{ij})}.$$
Then we have\\
$\text{ \,\,\,\,}\dfrac{d}{dt}\text{Area}(F(\cdot,t))$\\
$= \dfrac{d}{dt}\int_{D_{+}}\sqrt{\det (\tilde{g}_{ij})}$\\
$= \int_{D_{+}}\dfrac{\partial_{t}\det (\tilde{g}_{ij})}{2\sqrt{\det (\tilde{g}_{ij})}}$\\
$= \int_{D_{+}}\dfrac{1}{2}\tilde{g}^{ij}\partial_{t}\tilde{g}_{ij}\sqrt{\det (\tilde{g}_{ij})}$\\
$= \int_{F(\cdot,t)(D_{+})}\dfrac{1}{2}\tilde{g}^{ij}\nabla_{F_{t}}<F_{i},F_{j}>$\\
$= \int_{F(\cdot,t)(D_{+})}\tilde{g}^{ij}<\nabla_{F_{t}}F_{i},F_{j}>$ \\
$\text{ \,\,\,\,}$ since $\tilde{g}^{ij}$ is symmetric\\
$= \int_{F(\cdot,t)(D_{+})}\tilde{g}^{ij}<\nabla_{F_{i}}F_{t},F_{j}> $
$\text{ \,\,\,\,}$ since $F$ is $C^{2,\sigma}$\\
$= \int_{F(\cdot,t)(D_{+})}[\tilde{g}^{ij}\nabla_{F_{i}}<F_{t},F_{j}> - \tilde{g}^{ij}<F_{t},\nabla_{F_{i}}F_{j}>]$\\
$= \int_{F(\cdot,t)(D_{+})}[\tilde{g}^{ij}\nabla_{F_{i}}<(F_{t})^{\top},F_{j}> - \tilde{g}^{ij}<F_{t},\nabla_{F_{i}}F_{j}>]$\\
$\text{ \,\,\,\,}$ where $(F_{t})^{\top}$ is the tangential component of $F_{t}$\\
$= \int_{F(\cdot,t)(D_{+})}[\text{div}_{F(\cdot,t)(D_{+})}[(F_{t})^{\top}] -<F_{t},\tilde{g}^{ij}\nabla_{F_{i}}F_{j}>]$\\
$= \int_{\partial[F(\cdot,t)(D_{+})]}<(F_{t})^{\top},\vec{n}> - \int_{F(\cdot,t)(D_{+})}<F_{t},\tilde{g}^{ij}\nabla_{F_{i}}F_{j}>$\\
$\text{ \,\,\,\,}$ by the Stoke's Theorem since $F$ is $C^{2,\sigma}$\\
$= \int_{\partial[F(\cdot,t)(D_{+})]}<F_{t},\vec{n}> - \int_{F(\cdot,t)(D_{+})}<F_{t},\tilde{g}^{ij}\nabla_{F_{i}}F_{j}>.$\\
Thus we have
$$\dfrac{d}{dt}\Big|_{t=0}\text{Area}(F(\cdot,t)) = -\int_{f(D_{+})}<\vec{H}_{f},\vec{X}> + \int_{f(\Gamma_{1})}<\vec{n},\vec{X}> + \int_{f(\Gamma_{2})}<\vec{n},\vec{X}>.$$
\hfill q.e.d.

We say $f \in C_{*}^{2,\sigma}(\overline{D_{+}}, N_{+})$ is minimal if $\vec{H}_{f} = 0$ and has partially free boundary if $<\vec{n}, \vec{v}>|_{f(\Gamma_{2})} = 0$ for any $\vec{v} \in T_{f(\Gamma_{2})}(\partial N_{+})$. We use $\mathcal{M} = \{[f]: f \in C_{*}^{2,\sigma}(\overline{D_{+}},N_{+}), \vec{H}_{f} = 0, <\vec{n}, \vec{v}>|_{f(\Gamma_{2})} = 0$ for any $\vec{v} \in T_{f(\Gamma_{2})}(\partial N_{+}) \}$ to denote the space of partially free boundary minimal half disks.

\subsubsection{The second variation}
Let $[f_{0}] \in \mathcal{M}$. The unit normal vector field of $f_{0}(D_{+})$ in $N_{+}$ would only be $C^{1,\sigma}$. We choose a smooth unit vector field $\vec{p}: \overline{D_{+}} \to \mathbb R^3$ such that $\vec{p}$ complements the image of $Df_{0}$ and $\vec{p}(x) \in T_{f_{0}(x)}(\partial N_{+})$ for each $x \in \Gamma_{2}$.

For any $f$ near to $f_{0}$, since $\vec{H}_{f}$ is normal to $f(D_{+})$, by the choice of $\vec{p}$, $\vec{H}_{f} = 0$ if and only if $<\vec{H}_{f}, \vec{p}> = 0$.

Since $T_{f(\Gamma_{2})}(\partial N_{+}) = \text{span}\{\vec{p}, \dfrac{\partial f}{\partial \tau_{2}}\}$, we have $<\vec{n}, \vec{v}>|_{f(\Gamma_{2})} = 0$ for any $\vec{v} \in T_{f(\Gamma_{2})}(\partial N_{+})$ if and only if $<\vec{n}, \vec{p}> = 0$ on $f(\Gamma_{2})$.

Define $H_{f}(x) = <\vec{H}_{f}(f(x)), \vec{p}>_{f(x)}$ for each $x \in D_{+}$, $\Theta_{f}(x) = <\vec{n}(f(x)), \vec{p}>_{f(x)}$ for each $x \in \Gamma_{2}$.

Thus for $f$ near to $f_{0}$, $[f] \in \mathcal M$ if and only if $H_{f}=0, \Theta_{f} = 0$.
\begin{thm}
With the above notations, the map
$$(H,\Theta): C_{*}^{2,\sigma}(\overline{D_{+}}, N_{+}) \to C^{0,\sigma}(\overline{D_{+}}) \times C^{1,\sigma}(\overline{\Gamma_{2}})$$
is a $C^{q-1,\beta-\sigma}$ map.
\end{thm}

To prove Theorem 3.2.2, we need a lemma in White\cite{White1987}.
\begin{lemma}[White\cite{White1987}, Appendix]
Let $M$ be a bounded subset in $\mathbb R^n$, $E_{1},E_{2}$ be two Euclidean spaces.\\
If $0<\alpha < \beta < 1$ and $g \in C^{j+k,\beta}(E_{1},E_{2})$, then the map
$$G_{g}: u \to g\circ u$$
is a $C^{k,\beta-\alpha}$ map from $C^{j+1,\alpha}(M,E_{1})$ to $C^{j,\alpha}(M,E_{2})$.
\end{lemma}
\emph{Proof of Theorem 3.2.2:}
For any $f \in C_{*}^{2,\sigma}(\overline{D_{+}}, N_{+})$, let $f = (f^{1},f^{2},f^{3})$ and $\{\partial_{1},\partial_{2},\partial_{3}\}$ be a basis of $TN_{+}$. Thus we have
\begin{align*}
&\quad \, \vec{H}_{f}\\
&= \tilde{g}^{ij}\nabla_{f_{i}}f_{j}\\
&= \tilde{g}^{ij}\nabla_{f_{i}^{k}\partial_{k}}f_{j}^{l}\partial_{l}\\
&= \tilde{g}^{ij}f_{ij} + \tilde{g}^{ij}f_{i}^{k}f_{j}^{l}\Gamma_{kl}^{m}(f(x))\partial_{m}
\end{align*}
where $\Gamma_{kl}^{m}$ is the Christoffel symbol of $N_{+}$. Since $g$ is $C^{q,\beta}$ and $\Gamma_{kl}^{m}$ is $C^{q-1,\beta}$, by Lemma 3.2.1,
$$\tilde{g}_{ij}(Df): C^{1,\sigma}(\overline{D_{+}}, TN_{+}) \times C^{1,\sigma}(\overline{D_{+}}, TN_{+}) \to C^{0,\sigma}(\overline{D_{+}})$$
defined by
$$[\tilde{g}_{ij}(Df)](x) = <f_{i},f_{j}>_{f(x)}$$
is a $C^{q,\beta-\sigma}$ map. And also the map
$$\Gamma_{kl}^{m}(f): C^{1,\sigma}(\overline{D_{+}}, N_{+}) \to C^{0,\sigma}(\overline{D_{+}})$$
defined by
$$[\Gamma_{kl}^{m}(f)](x) = \Gamma_{kl}^{m}(f(x))$$
is a $C^{q-1,\beta-\sigma}$ map.\\
Since $f \to Df, D^{2}f$ are analytic, we have
$$H: C_{*}^{2,\sigma}(\overline{D_{+}}, N_{+}) \to C^{0,\sigma}(\overline{D_{+}})$$
defined by
$$H_{f} = <\vec{H}_{f},\vec{p}>$$
is a $C^{q-1,\beta-\sigma}$ map.\\
Similarly, the map
$$\Theta: C_{*}^{2,\sigma}(\overline{D_{+}}, N_{+}) \to C^{1,\sigma}(\overline{\Gamma_{2}})$$
defined by
$$\Theta_{f} = <\vec{n},\vec{p}>\big|_{f(\Gamma_{2})}$$
is also a $C^{q-1,\beta-\sigma}$ map since it only involves the metric $g$ and $Df$.

\hfill q.e.d.\\

Let $C_{\Gamma_{1}}^{2,\sigma}(\overline{D_{+}}) = \{h \in C^{2,\sigma}(\overline{D_{+}}): h|_{\Gamma_{1}}=0\}$ and $C_{0}^{1,\sigma}(\overline{\Gamma_{2}}) = \{g \in C^{1,\sigma}(\overline{\Gamma_{2}}): g(S_{1}) = g(S_{2})=0 \}$.

Then for any $f \in C_{*}^{2,\sigma}(\overline{D_{+}}, N_{+})$ which is near to $f_{0}$, there is a $h \in C_{\Gamma_{1}}^{2,\sigma}(\overline{D_{+}})$ and a $\tilde{f} \in C_{*}^{2,\sigma}(\overline{D_{+}}, N_{+})$ such that
$$[f] = [\tilde{f}], \tilde{f} = \phi(f|_{\Gamma_{1}}) + h\vec{p}$$
where $\phi: C_{\perp}^{2,\sigma}(\overline{\Gamma_{1}},N_{+}) \to C_{*}^{2,\sigma}(\overline{D_{+}}, N_{+})$ is defined as follows:\\
for any $\gamma \in C_{\perp}^{2,\sigma}(\overline{\Gamma_{1}},N_{+})$, $\phi(\gamma)$ is the unique solution to
$$\begin{cases}
  \triangle(\phi(\gamma)-f_{0}) - (\phi(\gamma)-f_{0})=0 &\text{on $D_{+}$},\\
  \phi(\gamma)-f_{0}|_{\Gamma_{1}} = \gamma-f_{0}|_{\Gamma_{1}} &\text{on $\Gamma_{1}$}, \\
  \dfrac{\partial (\phi(\gamma)-f_{0})}{\partial \nu_{2}}(x,0) + (\phi(\gamma)-f_{0}) \\
  = \dfrac{1-x}{2}\cdot\dfrac{\partial (\gamma-f_{0}|_{\Gamma_{1}})}{\partial \tau_{1}}(-1,0) + \dfrac{1+x}{2}\cdot\dfrac{\partial (\gamma-f_{0}|_{\Gamma_{1}})}{-\partial \tau_{1}}(1,0) &\text{on $\Gamma_{2}$}.
\end{cases}$$
This problem is uniquely solvable in $C_{*}^{2,\sigma}(\overline{D_{+}}, N_{+})$ by Theorem 2.7.6 and we have the estimate
\begin{align*}
&\|\phi(\gamma) - f_{0}\|_{2,\sigma, \overline{D_{+}}} \\
\leq &C\{0 + \|\gamma - f_{0}|_{\Gamma_{1}}\|_{2,\sigma,\overline{\Gamma}_{1}} \\
\quad &+  \|\dfrac{1-x}{2}\cdot\dfrac{\partial (\gamma-f_{0}|_{\Gamma_{1}})}{\partial \tau_{1}}(-1,0) + \dfrac{1+x}{2}\cdot\dfrac{\partial (\gamma-f_{0}|_{\Gamma_{1}})}{-\partial \tau_{1}}(1,0)\|_{1,\sigma,\overline{\Gamma_{2}}}    \}\\
\leq &C^{'}\|\gamma - f_{0}|_{\Gamma_{1}}\|_{2,\sigma,\overline{\Gamma}_{1}}.
\end{align*}
by Theorem 2.5.6.

Thus we consider the map
$$\Phi: C_{\perp}^{2,\sigma}(\overline{\Gamma_{1}}, N_{+}) \times C_{\Gamma_{1}}^{2,\sigma}(\overline{D_{+}}) \to C^{0,\sigma}(\overline{D_{+}}) \times C_{0}^{1,\sigma}(\overline{\Gamma_{2}})$$
defined by
$$\Phi(\gamma, h) = (H_{\phi(\gamma)+h\vec{p}}, \Theta_{\phi(\gamma)+h\vec{p}}).$$

The reason $<\vec{n}(f(x)), \vec{p}>_{f(x)} = 0$ at $S_{1}$ and $S_{2}$ for any $f \in C_{*}^{2,\sigma}(\overline{D_{+}}, N_{+})$ is that we already assume that  $f(\Gamma_{1}) \perp \partial N_{+}$.

Next we compute the derivative of $\Phi$.
\begin{thm}
With the above notations, we have
$$[D_{2}\Phi(f_{0}|_{\Gamma_{1}}, 0)](h) = (\tilde{g}^{ij}h_{ij}+Q(h), [\alpha_{2}\dfrac{\partial h}{\partial \nu_{2}} + \beta_{2}\dfrac{\partial h}{\partial \tau_{2}}+Rh]|_{\Gamma_{2}} )$$
where $\tilde{g}^{ij}$ is the inverse matrix of $\tilde{g}_{ij} = <\partial_{i}f_{0}, \partial_{j}f_{0}>$, $Q(h)$ is a linear first order differential operator on $h$, $\nu_{2}$ is the outward unit normal vector of $D_{+}$ on $\Gamma_{2}$, $\tau_{2}$ is the unit tangent vector of $\Gamma_{2}$ and $\alpha_{2}, \beta_{2}, R \in C^{1,\sigma}(\overline{\Gamma_{2}})$ such that $\alpha_{2} > 0$, $\beta_{2}(S_{1}) = \beta(S_{2}) = 0$.
\end{thm}
\emph{Proof:} Since $[f_{0}] \in \mathcal M$, we have $\vec{H}_{f_{0}} = 0$. \\
For any $h \in C_{\Gamma_{1}}^{2,\sigma}(\overline{D_{+}})$, define
$$F(x,t): D_{+} \times (-\epsilon, \epsilon) \to N_{+}$$
by
$$F(x,t) = f_{0}+th\vec{p}.$$
Then $F(\cdot,t) \in C_{*}^{2,\sigma}(\overline{D_{+}}, N_{+})$ since by the assumption that $N_{+} = (\mathbb R_{+}^{3},g)$ and $\vec{p}(x) \in T_{f_{0}(x)}(\partial N_{+})$ for each $x \in \Gamma_{2}$.\\
Let $F = (F^{1},F^{2},F^{3})$ and $\{\partial_{1},\partial_{2},\partial_{3}\}$ be a basis of $TN_{+}$, we have\\
$\text{ \,\,\,\,}\dfrac{d}{dt}\Big|_{t=0} H_{F(\cdot,t)}$\\
$= \dfrac{d}{dt}\Big|_{t=0} <\vec{H}_{F(\cdot,t)},\vec{p}>$\\
$= <\dfrac{d}{dt}\Big|_{t=0}\vec{H}_{F(\cdot,t)},\vec{p}>$ since $\vec{H}_{f_{0}} = 0$\\
$= <\dfrac{d}{dt}\Big|_{t=0} \tilde{g}^{ij}\nabla_{F_{i}}F_{j}, \vec{p}>$ \\
$= <\left[[\partial_{t}\tilde{g}^{ij}]\nabla_{F_{i}}F_{j}+ \tilde{g}^{ij}\nabla_{F_{t}}\nabla_{F_{i}}F_{j}\right], \vec{p}>\Big|_{t=0}$\\
$= <\left[-\tilde{g}^{ik}\tilde{g}^{jl}[\partial_{t}\tilde{g}_{kl}]\nabla_{F_{i}}F_{j}+ \tilde{g}^{ij}\nabla_{F_{t}}\nabla_{F_{i}^{k}\partial_{k}}(F_{j}^{l}\partial_{l})\right], \vec{p}>\Big|_{t=0}$\\
$= <\left[-\tilde{g}^{ik}\tilde{g}^{jl}[\partial_{t}<F_{k},F_{l}>]\nabla_{F_{i}}F_{j}+ \tilde{g}^{ij}\nabla_{F_{t}}\nabla_{F_{i}^{k}\partial_{k}}(F_{j}^{l}\partial_{l})\right], \vec{p}>\Big|_{t=0}$\\
$= <\left[-\tilde{g}^{ik}\tilde{g}^{jl}[<\nabla_{F_{t}}F_{k},F_{l}> + <F_{k},\nabla_{F_{t}}F_{l}>]\nabla_{F_{i}}F_{j}\right.$\\
$\text{ \,\,\,\,} \left.+ \tilde{g}^{ij}\nabla_{F_{t}}(F_{ij}^{l}\partial_{l} + F_{i}^{k}F_{j}^{l}\nabla_{\partial_{k}}\partial_{l})\right], \vec{p}>\Big|_{t=0}$\\
$= <\left[-\tilde{g}^{ik}\tilde{g}^{jl}[<\nabla_{F_{k}}F_{t},F_{l}> + <F_{k},\nabla_{F_{l}}F_{t}>]\nabla_{F_{i}}F_{j}\right.$\\
$\text{ \,\,\,\,} \left.+ \tilde{g}^{ij}(F_{ijt}^{l}\partial_{l} + F_{ij}^{l}\nabla_{F_{t}}\partial_{l} + F_{it}^{k}F_{j}^{l}\nabla_{\partial_{k}}\partial_{l} + F_{i}^{k}F_{jt}^{l}\nabla_{\partial_{k}}\partial_{l} \right.$\\
$\text{ \,\,\,\,} \left.+ F_{i}^{k}F_{j}^{l}\nabla_{F_{t}}\nabla_{\partial_{k}}\partial_{l})\right], \vec{p}>\Big|_{t=0}$\\
$= <\tilde{g}^{ij}F_{ijt}^{l}\partial_{l}$\\
$\text{ \,\,\,\,} + \left[-\tilde{g}^{ik}\tilde{g}^{jl}[<\nabla_{F_{k}}F_{t},F_{l}> + <F_{k},\nabla_{F_{l}}F_{t}>]\nabla_{F_{i}}F_{j}\right.$\\
$\text{ \,\,\,\,} \left.+ \tilde{g}^{ij}(F_{ij}^{l}\nabla_{F_{t}}\partial_{l} + F_{it}^{k}F_{j}^{l}\nabla_{\partial_{k}}\partial_{l} + F_{i}^{k}F_{jt}^{l}\nabla_{\partial_{k}}\partial_{l} \right.$\\
$\text{ \,\,\,\,} \left.+ F_{i}^{k}F_{j}^{l}\nabla_{F_{t}}\nabla_{\partial_{k}}\partial_{l})\right], \vec{p}>\Big|_{t=0}.$\\
Since $F(x,t) = f_{0}(x) + th\vec{p}$, we have
$$\dfrac{d}{dt}\Big|_{t=0} H_{F(\cdot,t)} = \tilde{g}^{ij}h_{ij} + Q(h)$$
where $Q(h)$ is a linear first order differential operator on $h$.\\
Next we compute $\dfrac{d}{dt}\Big|_{t=0} \Theta_{F(\cdot,t)}$.\\
By definition,
$$\vec{n}_{F(\cdot,t)}\big|_{f_{0}(\Gamma_{2})} = \dfrac{-F_{y}+\frac{<F_{x},F_{y}>}{<F_{x},F_{x}>}F_{x}}{\|-F_{y}+\frac{<F_{x},F_{y}>}{<F_{x},F_{x}>}F_{x}\|}$$
where $\dfrac{<F_{x},F_{y}>}{<F_{x},F_{x}>}(S_{1}) = \dfrac{<F_{x},F_{y}>}{<F_{x},F_{x}>}(S_{2}) = 0$ since $f_{0}\big|_{\Gamma_{1}} \perp \partial N_{+}$ and $h\big|_{\Gamma_{1}} = 0$.\\
Since $\Theta_{f_{0}} = 0$, we have
\begin{align*}
\dfrac{d}{dt}\Big|_{t=0} \Theta_{F(\cdot,t)}
&= \dfrac{d}{dt}\Big|_{t=0}<\vec{n}_{F(\cdot,t)}, \vec{p}> \\
&=\dfrac{d}{dt}\Big|_{t=0}<\dfrac{-F_{y}+\frac{<F_{x},F_{y}>}{<F_{x},F_{x}>}F_{x}}{\|-F_{y}+\frac{<F_{x},F_{y}>}{<F_{x},F_{x}>}F_{x}\|}, \vec{p}> \\
&= \dfrac{\dfrac{d}{dt}\Big|_{t=0}<-F_{y}+\frac{<F_{x},F_{y}>}{<F_{x},F_{x}>}F_{x}, \vec{p}>}{\|-F_{y}+\frac{<F_{x},F_{y}>}{<F_{x},F_{x}>}F_{x}\|} \\
&= \left[\dfrac{-<\nabla_{F_{t}}F_{y}+(\nabla_{F_{t}}\frac{<F_{x},F_{y}>}{<F_{x},F_{x}>})F_{x} + \frac{<F_{x},F_{y}>}{<F_{x},F_{x}>}\nabla_{F_{t}}F_{x}, \vec{p}>}{\|-F_{y}+\frac{<F_{x},F_{y}>}{<F_{x},F_{x}>}F_{x}\|}\right. \\
&\quad \left. + \dfrac{<-F_{y}+\frac{<F_{x},F_{y}>}{<F_{x},F_{x}>}F_{x}, \nabla_{F_{t}}\vec{p}>}{\|-F_{y}+\frac{<F_{x},F_{y}>}{<F_{x},F_{x}>}F_{x}\|}\right]\Bigg|_{t=0}.
\end{align*}
Since $\dfrac{<F_{x},F_{y}>}{<F_{x},F_{x}>}(S_{1}) = \dfrac{<F_{x},F_{y}>}{<F_{x},F_{x}>}(S_{2}) = 0$ for each $t$ and $F(\cdot,t) = f_{0} + th\vec{p}$, we have
$$\dfrac{d}{dt}\Big|_{t=0} \Theta_{F(\cdot,t)} = -\alpha_{2}h_{y} + \beta_{2}h_{x}+Rh$$
where $\alpha_{2}, \beta_{2}, R \in C^{1,\sigma}(\overline{\Gamma_{2}})$ such that $\alpha_{2} > 0$, $\beta_{2}(S_{1}) = \beta(S_{2}) = 0$.

\hfill q.e.d.\\

Remark: In these two formulas, there is no Ricci tensor or the second fundamental form which appear in the usual second variation formula of the area. It is because from the above computations, we know that $DH_{f_{0}}$ is a linear second order differential operator, $D\Theta_{f_{0}}$ is a linear first order differential operator and we only care about the highest order terms.

We also have the Green's formula, for all $h_{1}, h_{2} \in C^{2,\sigma}(\overline{D_{+}})$,\\
$\text{ \,\,\,\,} -\int_{D_{+}} h_{1}[DH_{f_{0}}(h_{2})]dA_{f_{0}}
+ \int_{\Gamma_{1}} h_{1}[\alpha_{1}\dfrac{\partial h_{2}}{\partial \nu_{1}} + \beta_{1}\dfrac{\partial h_{2}}{\partial \tau_{1}}]|_{\Gamma_{1}} ds_{f_{0}} + \int_{\Gamma_{2}} h_{1}[D\Theta_{f_{0}}(h_{2})]ds_{f_{0}}$\\
$= -\int_{D_{+}} h_{2}[DH_{f_{0}}(h_{1})]dA_{f_{0}}
+ \int_{\Gamma_{1}} h_{2}[\alpha_{1}\dfrac{\partial h_{1}}{\partial \nu_{1}} + \beta_{1}\dfrac{\partial h_{1}}{\partial \tau_{1}}]|_{\Gamma_{1}} ds_{f_{0}} + \int_{\Gamma_{2}} h_{2}[D\Theta_{f_{0}}(h_{1})]ds_{f_{0}}$\\
where $\alpha_{1}, \beta_{1} \in C^{1,\sigma}(\overline{\Gamma_{1}})$ such that $Df_{0}(\alpha_{1}\nu_{1}+\beta_{1}\tau_{1})$ is the outward unit normal vector of $f_{0}(D_{+})$ on $f_{0}(\Gamma_{1})$, $dA_{f_{0}}$ and $ds_{f_{0}}$ are the pull-back area and length integrals of $f_{0}$.

\subsection{The elliptic theory}
From the elliptic theory we construct in Section 2, we have the following theorem.
\begin{thm}
  Let $J(u) = a_{ij}D_{ij}u+ b_{i}u_{i} + cu$ be a strongly elliptic operator on $D_{+}$ such that $a_{ij} = a_{ji},b_{i},c \in C^{0,\sigma}(\overline{D_{+}})$, $a_{ij}$ is diagonal at $S_{1}$ and $S_{2}$. Let $\beta, d \in C^{1,\sigma}(\overline{\Gamma_{2}})$ such that $\beta(S_{1}) = \beta(S_{2}) = 0$. \\
  then the map
  $$L: C_{\Gamma_{1}}^{2,\sigma}(\overline{D_{+}}) \to C^{0,\sigma}(\overline{D_{+}}) \times C_{0}^{1,\sigma}(\overline{\Gamma_{2}})$$
  defined by $L(u) = (J(u), [\dfrac{\partial u}{\partial \nu_{2}}+\beta\dfrac{\partial u}{\partial \tau_{2}} + du]|_{\Gamma_{2}})$ is a Fredholm map of index zero where $\nu_{2}$ is the outward unit normal vector of $D_{+}$ on $\Gamma_{2}$, $\tau_{2}$ is the unit tangent vector of $\Gamma_{2}$.
\end{thm}
\emph{Proof:}
By Theorem 2.7.6,
$$\tilde{L}(u) = L(u) + ((-\sup|c|-1)u, (\sup|d| + 1)u\big|_{\Gamma_{2}})$$
is one-to-one and onto. \\
Since $((-\sup|c|-1)u, (\sup|d| + 1)u\big|_{\Gamma_{2}})$ is a compact operator, we have
$$\text{index}(\tilde{L}) = \text{index}(L) = 0.$$
\hfill q.e.d.\\

We have the following direct corollary.
\begin{corollary}
  Let $[f_{0}] \in \mathcal M$ and $\Phi: C_{\perp}^{2,\sigma}(\overline{\Gamma_{1}}, N_{+}) \times C_{\Gamma_{1}}^{2,\sigma}(\overline{D_{+}}) \to C^{0,\sigma}(\overline{D_{+}}) \times C_{0}^{1,\sigma}(\overline{\Gamma_{2}})$ be as in Section 3.2.4, then $D_{2}\Phi(f_{0}|_{\Gamma_{1}}, 0)$ is a Fredholm map of index zero. \\
  Let $K = \text{ker}[D_{2}\Phi(f_{0}|_{\Gamma_{1}}, 0)]$, then $[K \times \{0\}] \oplus \text{image}[D_{2}\Phi(f_{0}|_{\Gamma_{1}}, 0)] = C^{0,\sigma}(\overline{D_{+}}) \times C_{0}^{1,\sigma}(\overline{\Gamma_{2}})$.
\end{corollary}
\emph{Proof:} By Theorem 3.2.2,
$$D_{2}\Phi(f_{0}|_{\Gamma_{1}},0)(h) = (\tilde{g}^{ij}h_{ij}+Q(h), [\alpha_{2}\dfrac{\partial h}{\partial \nu_{2}} + \beta_{2}\dfrac{\partial h}{\partial \mu_{2}}+Rh]|_{\Gamma_{2}} )$$
where $Q$ is a first order differential operator, $\tilde{g}^{ij}$ is the inverse of $\tilde{g}_{ij} = <\partial_{i}f_{0}, \partial_{j}f_{0}>$, $\alpha_{2}, \beta_{2}, R \in C^{1,\sigma}(\overline{\Gamma_{2}})$ such that $\alpha > 0$, $\beta(S_{1}) = \beta(S_{2}) = 0$.
Then by Theorem 3.3.1,
$$D_{2}\Phi(f_{0}|_{\Gamma_{1}},0)$$
is a Fredholm map of index zero.\\
Let $K = \text{ker}[D_{2}\Phi(f_{0}|_{\Gamma_{1}},0)]$, then for any $\kappa \in K$ and $u \in C_{\Gamma_{1}}^{2,\sigma}(\overline{D_{+}})$, by the Green's formula in Section 3.2.4, we have\\
$-\int_{D_{+}} u[DH_{f_{0}}(\kappa)]dA_{f_{0}} + \int_{\Gamma_{1}} u[\alpha_{1}\dfrac{\partial \kappa}{\partial \nu_{1}} + \beta_{1}\dfrac{\partial \kappa}{\partial \tau_{1}}]|_{\Gamma_{1}} ds_{f_{0}} + \int_{\Gamma_{2}} u[D\Theta_{f_{0}}(\kappa)]ds_{f_{0}}$\\
$=-\int_{D_{+}} \kappa[DH_{f_{0}}(u)]dA_{f_{0}} + \int_{\Gamma_{1}} \kappa[\alpha_{1}\dfrac{\partial u}{\partial \nu_{1}} + \beta_{1}\dfrac{\partial u}{\partial \tau_{1}}]|_{\Gamma_{1}} ds_{f_{0}} + \int_{\Gamma_{2}} \kappa[D\Theta_{f_{0}}(u)]ds_{f_{0}}$.\\
Which is
$$0 + 0  + 0 = -\int_{D_{+}} \kappa[DH_{f_{0}}(u)]dA_{f_{0}} + 0 + \int_{\Gamma_{2}} \kappa[D\Theta_{f_{0}}(u)]ds_{f_{0}}.$$
Thus there is no $u \in C_{\Gamma_{1}}^{2,\sigma}(\overline{D_{+}})$ such that
$$DH_{f_{0}}(u) = \kappa, D\Theta_{f_{0}}(u) = 0.$$
Which means
$$(\kappa, 0) \notin \text{image}[D_{2}\Phi(f_{0}|_{\Gamma_{1}}, 0)].$$
Since $D_{2}\Phi(f_{0}|_{\Gamma_{1}}, 0)$ is a Fredholm map of index zero, we get the desired conclusion.
\hfill q.e.d.

\subsection{Structure of the space of minimal half disks}
Follow the same idea in White\cite{White1987}, we have the following local structure theorem of $\mathcal M$.
\begin{thm}
  Let $[f_{0}] \in \mathcal{M}$. If $\text{ker}[D_{2}\Phi(f_{0}|_{\Gamma_{1}}, 0)] = \{0\}$, then there exists a $C^{q-1}$ map
  $$F: U \subset C_{\perp}^{2,\sigma}(\overline{\Gamma_{1}}, N_{+}) \to C_{*}^{2,\sigma}(\overline{D_{+}}, N_{+})$$
  defined on some neighborhood $U$ of $f_{0}|_{\Gamma_{1}}$ such that for all $\gamma \in U$,\\
  (1) $[F(\gamma)] \in \mathcal M$.\\
  (2) $F(\gamma)|_{\Gamma_{1}} = \gamma$.\\
  (3) $F(f_{0}|_{\Gamma_{1}}) = f_{0}$\\
  (4) There is a neighborhood $W$ of $f_{0}$ in $C_{*}^{2,\sigma}(\overline{D_{+}}, N_{+})$ such that if $f \in W$ and $[f] \in \mathcal M$, then $f|_{\Gamma_{1}} \in U$ and
  $$[F(f|_{\Gamma_{1}})] = [f].$$
  (5) Furthermore, we can choose $U$ and $W$ above so that they are open in the $C^{2,\sigma^{'}}$ topology for any $0< \sigma^{'} < \sigma$ (we need this because $C_{\perp}^{2,\sigma}(\overline{\Gamma_{1}}, N_{+})$ is only separable with a weaker  $C^{2,\sigma^{'}}$ norm).
\end{thm}
\emph{Proof:} Since $\text{ker}[D_{2}\Phi(f_{0}|_{\Gamma_{1}}, 0)]=\{0\}$, we know $D_{2}\Phi(f_{0}|_{\Gamma_{1}}, 0)$ is an isomorphism. \\
Therefore by the implicit function theorem, there are neighborhoods $U$ of $f_{0}|_{\Gamma_{1}}$ in $C_{\perp}^{2,\sigma}(\overline{\Gamma_{1}}, N_{+})$, $V$ of $0$ in $C_{\Gamma_{1}}^{2,\sigma}(\overline{D_{+}})$ and a $C^{q-1,\beta-\sigma}$ map $$F_{0}: U \to V$$
such that for any $\gamma \in U$, $F_{0}(\gamma)$ is the unique solution in $V$ to $$\Phi(\gamma,F_{0}(\gamma)) = 0.$$
Define $F(\gamma) = \phi(\gamma) + F_{0}(\gamma)\vec{p}$. Then $F$ is what we need.\\
Conclusion (4) comes from the fact that for any $f \in C_{*}^{2,\sigma}(\overline{D_{+}}, N_{+})$ which is near to $f_{0}$, there is a $h \in C_{\Gamma_{1}}^{2,\sigma}(\overline{D_{+}})$ which is near to $0$ and a $\tilde{f} \in C_{*}^{2,\sigma}(\overline{D_{+}}, N_{+})$ such that
$$[f] = [\tilde{f}], \tilde{f} = \phi(f|_{\Gamma_{1}}) + h\vec{p}.$$
To prove (5), for any $0<\sigma^{'}<\sigma$, repeat the previous steps with $C^{2,\sigma^{'}}$ spaces. \\
Then there are neighborhoods $U^{'}$ of $f_{0}|_{\Gamma_{1}}$ in $C_{\perp}^{2,\sigma^{'}}(\overline{\Gamma_{1}}, N_{+})$, $V^{'}$ of $0$ in $C_{\Gamma_{1}}^{2,\sigma^{'}}(\overline{D_{+}})$ and a $C^{q-1,\beta-\sigma^{'}}$ map
$$F_{0}^{'}: U^{'} \to V^{'}$$
such that for any $\gamma \in U^{'}$, $F_{0}^{'}(\gamma)$ is the unique solution in $V^{'}$ to
$$\Phi(\gamma,F_{0}^{'}(\gamma)) = 0.$$
In the following we show that for any $\gamma \in U^{'} \cap C_{\perp}^{2,\sigma}(\overline{\Gamma_{1}}, N_{+})$, $$F_{0}^{'}(\gamma) \in C_{\Gamma_{1}}^{2,\sigma}(\overline{D_{+}}).$$
For any $\gamma \in U^{'} \cap C_{\perp}^{2,\sigma}(\overline{\Gamma_{1}}, N_{+})$, since $$\Phi(\gamma,F_{0}^{'}(\gamma)) = (H_{\phi(\gamma)+F_{0}^{'}(\gamma)\vec{p}}, \Theta_{\phi(\gamma)+F_{0}^{'}(\gamma)\vec{p}}) = 0,$$
we have that $u = F_{0}^{'}(\gamma) \in C_{\Gamma_{1}}^{2,\sigma^{'}}(\overline{D_{+}})$ satisfies the following quasilinear elliptic equation
$$\begin{cases}
  \tilde{g}^{ij}u_{ij}=f  &\text{on $D_{+}$},\\
  u = 0 &\text{on $\Gamma_{1}$}, \\
  -\alpha u_{y} + \beta u_{x} = g &\text{on $\Gamma_{2}$},
\end{cases}$$
where
$$\tilde{g}^{ij} = \tilde{g}^{ij}(x,Du), f = f(x,u,Du) \in C^{1,\sigma^{'}}(\overline{D_{+}}) \subset C^{0,\sigma}(\overline{D_{+}})$$
such that $\tilde{g}^{ij}$ is diagonal at $S_{1}$ and $S_{2}$,
$$\alpha = \alpha(x,\vec{p}, D\vec{p},\phi(\gamma), D\phi(\gamma), u) \in C^{1,\sigma}(\overline{\Gamma_{2}}),$$
$$\beta = \beta(x,\vec{p}, D\vec{p},\phi(\gamma), D\phi(\gamma), u) \in C^{1,\sigma}(\overline{\Gamma_{2}}),$$
$$g = g(x,\vec{p}, D\vec{p},\phi(\gamma), D\phi(\gamma), u) \in C^{1,\sigma}(\overline{\Gamma_{2}})$$
such that $\alpha > 0, \beta(S_{1}) = \beta(S_{2}) = 0$, $g(S_{1}) = g(S_{2}) = 0$. \\
Thus by Theorem 2.7.6, there is a unique $v \in C_{\Gamma_{1}}^{2,\sigma}(\overline{D_{+}})$ which is a solution to
$$\begin{cases}
  \tilde{g}^{ij}v_{ij} - v=f - u  &\text{on $D_{+}$},\\
  u = 0 &\text{on $\Gamma_{1}$}, \\
  -\alpha v_{y} + \beta v_{x} + v = g + u &\text{on $\Gamma_{2}$}.
\end{cases}$$
Let $w = u - v \in C_{\Gamma_{1}}^{2,\sigma^{'}}(\overline{D_{+}})$, then $w$ solves the following equation
$$\begin{cases}
  g^{ij}w_{ij} - w=0  &\text{on $D_{+}$},\\
  w = 0 &\text{on $\Gamma_{1}$}, \\
  -\alpha w_{y} + \beta w_{x} + w = 0 &\text{on $\Gamma_{2}$}.
\end{cases}$$
By the uniqueness of the solutions to the above linear elliptic problem in $C^{2,\sigma^{'}}(\overline{D_{+}})$, we have
$$w = 0,$$
which means
$$u = v \in C_{\Gamma_{1}}^{2,\sigma}(\overline{D_{+}}).$$
Therefore, for any $\gamma \in U^{'} \cap C_{\perp}^{2,\sigma}(\overline{\Gamma_{1}}, N_{+})$,
$$F_{0}^{'}(\gamma) \in C_{\Gamma_{1}}^{2,\sigma}(\overline{D_{+}}).$$
And by the uniqueness of $F_{0}$ on $U$, for any $\gamma \in U^{'} \cap U$,
$$F_{0}^{'}(\gamma) = F_{0}(\gamma).$$
This finishes the proof of (5).
\hfill q.e.d.\\

Next we prove the local structure theorem of $\mathcal M$ when there are some Jacobi fields, i.e, $\text{ker}[D_{2}\Phi(f_{0}|_{\Gamma_{1}}, 0)] \neq \{0\}$.
\begin{thm}
  Let $[f_{0}] \in \mathcal{M}$ and $K = \text{ker}[D_{2}\Phi(f_{0}|_{\Gamma_{1}}, 0)]$, then there exist $C^{q-1}$ maps
  $$F: U \subset C_{\perp}^{2,\sigma}(\overline{\Gamma_{1}}, N_{+}) \times K \to C_{*}^{2,\sigma}(\overline{D_{+}}, N_{+}),$$
  $$g: C_{\perp}^{2,\sigma}(\overline{\Gamma_{1}}, N_{+}) \times K \to K,$$
  defined on some neighborhood $U = U_{1} \times U_{2}$ of $<f_{0}|_{\Gamma_{1}}, 0>$ in $C_{\perp}^{2,\sigma}(\overline{\Gamma_{1}}, N_{+}) \times K$ such that for all $<\gamma, \kappa> \in U$,\\
  (1) $F(\gamma,\kappa)|_{\Gamma_{1}} = \gamma$.\\
  (2) $F(f_{0}|_{\Gamma_{1}},0) = f_{0}$.\\
  (3) $[F(\gamma,\kappa)] \in \mathcal M$ if and only if $g(\gamma,\kappa) = 0$.\\
  (4) $[D_{2}F(f_{0}|_{\Gamma_{1}},0)](\kappa) = \kappa\vec{p}$ for any $\kappa \in K$.\\
  (5) $g^{-1}(0) \cap U$ is a $C^{q-1}$ submanifold of codimension equal to $\dim(K)$. It contains $\{0\} \times K$ in its tangent space at $<f_{0}|_{\Gamma_{1}},0>$.\\
  (6) For every $\epsilon > 0$, there is a neighborhood $W$ of $f_{0}$ in $C_{*}^{2,\sigma}(\overline{D_{+}}, N_{+})$ such that if $f \in W$ and $[f] \in \mathcal M$, then there is a $\kappa \in K$ with $\|\kappa\|_{2,\sigma} < \epsilon$ such that
  $$[F(f|_{\Gamma_{1}}, \kappa)] = [f].$$
  (7) Furthermore, we can choose $U$ and $W$ above so that they are open in the $C^{2,\sigma^{'}}$ topology for any $0<\sigma^{'} < \sigma$ (we need this because $C_{\perp}^{2,\sigma}(\overline{\Gamma_{1}}, N_{+})$ is only separable with a weaker  $C^{2,\sigma^{'}}$ norm).
\end{thm}
\emph{Proof:} Define
$$\tilde{\Phi}: C_{\perp}^{2,\sigma}(\overline{\Gamma_{1}}, N_{+}) \times K \times K_{\Gamma_{1}}^{\perp} \to \text{image}[D_{2}\Phi(f_{0}|_{\Gamma_{1}},0)]$$
by
$$\tilde{\Phi}(\gamma,\kappa,h) = \Pi_{\text{image}[D_{2}\Phi(f_{0}|_{\Gamma_{1}},0)]} \circ \Phi(\gamma,\kappa+h)$$
where $K_{\Gamma_{1}}^{\perp}$ is the orthogonal complement of $K$ in $C_{\Gamma_{1}}^{2,\sigma}(\overline{D_{+}})$ with respect to the $L_{2}$ norm. \\
Thus $D_{3}\tilde{\Phi}(f_{0}|_{\Gamma_{1}},0,0)$ is an isomorphism. By the implicit function theorem, there are neighborhoods $U = U_{1} \times U_{2}$ of $<f_{0}|_{\Gamma_{1}}, 0>$ in $C_{\perp}^{2,\sigma}(\overline{\Gamma_{1}}, N_{+}) \times K$, $V$ of $\{0\}$ in $K_{\Gamma_{1}}^{\perp}$ and a $C^{q-1,\beta-\sigma}$ map
$$F_{0}: U \to V$$
such that for any $(\gamma,\kappa) \in U$, $F_{0}(\gamma,\kappa)$ is the unique solution in $V$ to
$$\tilde{\Phi}(\gamma,,\kappa,F_{0}(\gamma,\kappa)) = 0.$$
Define
$$F(\gamma,\kappa) =\phi(\gamma)+ (\kappa+F_{0}(\gamma,\kappa))\vec{p},$$
$$g(\gamma,\kappa) = \Pi_{K} \circ H_{F(\gamma,\kappa)}.$$
By Corollary 3.3.1,
$$[K \times \{0\}] \oplus \text{image}[D_{2}\Phi(f_{0}|_{\Gamma_{1}},0)]= C^{0,\sigma}(\overline{D_{+}}) \times C_{0}^{1,\sigma}(\overline{\Gamma_{2}}),$$
we know
$$[F(\gamma,\kappa)] \in \mathcal M \text{ if and only if } g(\gamma,\kappa) = 0.$$
This proves (1), (2) and (3).\\
To prove (4), since for any $(\gamma,\kappa) \in U$,
$$\Pi_{\text{image}[D_{2}\Phi(f_{0}|_{\Gamma_{1}},0)]} \circ (H_{\phi(\gamma)+(\kappa+F_{0}(\gamma,\kappa))\vec{p}}, \Theta_{\phi(\gamma)+(\kappa+F_{0}(\gamma,\kappa))\vec{p}}) = 0.$$
Therefore\\
$\text{ \quad $0$} \\
  = \dfrac{d}{dt}\Big|_{t=0} \Pi_{\text{image}[D_{2}\Phi(f_{0}|_{\Gamma_{1}},0)]} \circ (H_{f_{0}+(\kappa t+F_{0}(f_{0}|_{\Gamma_{1}},\kappa t))\vec{p}}, \Theta_{f_{0}+(\kappa t+F_{0}(f_{0}|_{\Gamma_{1}},\kappa t))\vec{p}})\\
  = \Pi_{\text{image}[D_{2}\Phi(f_{0}|_{\Gamma_{1}},0)]} \circ (DH_{f_{0}}(\kappa+D_{2}F_{0}((f_{0}|_{\Gamma_{1}},0)(\kappa))), D\Theta_{f_{0}}(\kappa+D_{2}F_{0}((f_{0}|_{\Gamma_{1}},0)(\kappa))))\\
  = \Pi_{\text{image}[D_{2}\Phi(f_{0}|_{\Gamma_{1}},0)]} \circ (DH_{f_{0}}(D_{2}F_{0}((f_{0}|_{\Gamma_{1}},0)(\kappa))), D\Theta_{f_{0}}(D_{2}F_{0}((f_{0}|_{\Gamma_{1}},0)(\kappa))))\\
  = (DH_{f_{0}}(D_{2}F_{0}((f_{0}|_{\Gamma_{1}},0)(\kappa))), D\Theta_{f_{0}}(D_{2}F_{0}((f_{0}|_{\Gamma_{1}},0)(\kappa))))
$.\\
Thus
$$D_{2}F_{0}((f_{0}|_{\Gamma_{1}},0)(\kappa)) \in K.$$
But by definition,
$$F_{0}(\gamma,\kappa) \in K_{\Gamma_{1}}^{\perp},$$
so we have
$$D_{2}F_{0}((f_{0}|_{\Gamma_{1}},0)(\kappa)) \in K \cap K_{\Gamma_{1}}^{\perp} = \{0\}.$$
Then $$[D_{2}F(f_{0}|_{\Gamma_{1}}, 0)](\kappa) = \kappa \vec{p} + D_{2}F_{0}((f_{0}|_{\Gamma_{1}},0)(\kappa))\vec{p} = \kappa \vec{p}.$$
To prove (5), we see for any $\kappa \in K$,
\begin{align*}
  \dfrac{d}{dt}\Big|_{t=0}g(f_{0}|_{\Gamma_{1}}, \kappa t)
  &= \Pi_{K} \circ DH_{f_{0}}([D_{2}F(f_{0}|_{\Gamma_{1}}, 0)](\kappa)) \\
  &= \Pi_{K} \circ DH_{f_{0}}(\kappa) \\
  &= 0.
\end{align*}
So we have
$$D_{2}g(f_{0}|_{\Gamma_{1}}, 0) = 0.$$
Then we show that $D_{1}g(f_{0}|_{\Gamma_{1}}, 0)$ is of full rank.\\
Let $\kappa \in K$ and $\kappa_{n} = [\alpha_{1}\dfrac{\partial \kappa}{\partial \nu_{1}} + \beta_{1}\dfrac{\partial \kappa}{\partial \tau_{1}}]\Big|_{\Gamma_{1}}$, by the Calder\'{o}n unique continuation theorem, $\dfrac{\partial \kappa}{\partial \nu_{1}}$ can't vanish on any open subset of $\Gamma_{1}$ unless $\kappa = 0$ (if $\dfrac{\partial \kappa}{\partial \nu_{1}} = 0$ on some open subset $\Gamma \subset \Gamma_{1}$, then $\kappa = 0$ on any smooth domain $M \subset \Omega$ such that $\Gamma \subset \partial M$, thus $\kappa = 0$ by considering $M$ tends to $\Omega$). \\
Since $\kappa|_{\Gamma_{1}} = 0$, we know $\kappa_{n}$ cannot vanish on any open subset of $\Gamma_{1}$ unless $\kappa = 0$. \\
Let $h \in C^{2,\sigma}(\overline{D_{+}})$ such that $h|_{\Gamma_{1}}$ is $C^{0}$ close to $\kappa_{n}|_{\Gamma_{1}}$, by the Green's formula in Section 3.2.4, we have
\begin{align*}
 &-\int_{D_{+}} u[DH_{f_{0}}(\kappa)]dA_{f_{0}} + \int_{\Gamma_{1}} u\kappa_{n} ds_{f_{0}} + \int_{\Gamma_{2}} u[D\Theta_{f_{0}}(\kappa)]ds_{f_{0}}\\
=&-\int_{D_{+}} \kappa[DH_{f_{0}}(u)]dA_{f_{0}} + \int_{\Gamma_{1}} \kappa u_{n} ds_{f_{0}} + \int_{\Gamma_{2}} \kappa[D\Theta_{f_{0}}(u)]ds_{f_{0}},\\
 &0 + \int_{\Gamma_{1}} u\kappa_{n} ds_{f_{0}} + 0 \\
=&-\int_{D_{+}} \kappa[DH_{f_{0}}(u)]dA_{f_{0}} + 0 + \int_{\Gamma_{2}} \kappa[D\Theta_{f_{0}}(u)]ds_{f_{0}}.
\end{align*}
Since $u|_{\Gamma_{1}}$ is $C^{0}$ close to $\kappa_{n}$, we know $(DH_{f_{0}}(u), D\Theta_{f_{0}}(u)) \neq (0,0)$. \\
By (1), $F(\gamma,\kappa)|_{\Gamma_{1}} = \gamma$, we have
$$[D_{1}F(f_{0}|_{\Gamma_{1}},0)(u\vec{p}|_{\Gamma_{1}})]|_{\Gamma_{1}} = u\vec{p}|_{\Gamma_{1}}$$
and
\begin{align*}
  \dfrac{d}{dt}\Big|_{t=0} g(f_{0}|_{\Gamma_{1}}+t u \vec{p}|_{\Gamma_{1}}, 0)
  &= \dfrac{d}{dt}\Big|_{t=0} \Pi_{K} \circ H_{F(f_{0}|_{\Gamma_{1}}+t u \vec{p}|_{\Gamma_{1}}, 0)}\\
  &= \dfrac{d}{dt}\Big|_{t=0} H_{F(f_{0}|_{\Gamma_{1}}+t u \vec{p}|_{\Gamma_{1}}, 0)}\\
  &= DH_{f_{0}}(\tilde{u})
\end{align*}
where $\tilde{u}\vec{p} = [D_{1}F(f_{0}|_{\Gamma_{1}},0)(u\vec{p}|_{\Gamma_{1}})]$ and $\tilde{u}|_{\Gamma_{1}} = u|_{\Gamma_{1}}$.\\
Since for any $(\gamma, \kappa) \in U$,
$$\Theta_{F(\gamma, \kappa)} = 0$$
and
$$(DH_{f_{0}}(u), D\Theta_{f_{0}}(u)) \neq (0,0),$$
we have
$$DH_{f_{0}}(\tilde{u}) \neq 0.$$
Thus we have
$$\text{rank}[D_{1}g(f_{0}|_{\Gamma_{1}},0)] = \dim(K).$$
This proves (5).\\
Conclusion (6) and (7) are similar to Theorem 3.4.1.
\hfill q.e.d.\\

Next we have the global structure theorem.
\begin{thm}
  Let $\mathcal M$ be as previous, then \\
  (1) $\mathcal{M}$ is a $C^{q-1}$ Banach manifold modelled on $C_{\perp}^{2,\sigma}(\overline{\Gamma_{1}}, N_{+})$ with a countable family of coordinate charts given in Theorem 3.4.2.\\
  (2) The map $\Pi: \mathcal M \to C_{\perp}^{2,\sigma}(\overline{\Gamma_{1}}, N_{+})$ defined by $\Pi([f]) = f|_{\Gamma_{1}}$ is a $C^{q-1}$ Fredholm map of index zero. \\
  (3) If $D\Pi([f_{0}])$ has $k$-dimensional kernel, then a neighborhood $U$ of $[f_{0}]$ in $\mathcal M$ may be identified with a $C^{q-1}$ codimension-$k$ submanifold of $C_{\perp}^{2,\sigma}(\overline{\Gamma_{1}}, N_{+}) \times \mathbb R^k$ (so that $\Pi$ corresponds to the obvious projection of $C_{\perp}^{2,\sigma}(\overline{\Gamma_{1}}, N_{+}) \times \mathbb R^k$). \\
  (4) $U$ may be chosen small enough so that there is a $k$-dimensional subspace $V$ of $C_{\perp}^{2,\sigma}(\overline{\Gamma_{1}}, N_{+})$ and $V \times \{0\}$ complements the tangent space to $U$ at each point of $U$.
\end{thm}
\emph{Proof:} Let $\omega_{1},...,\omega_{k}$ be real analytic two-forms on a neighborhood of $f_{0}(D_{+})$. Define
$$\Omega: C_{*}^{2,\sigma}(\overline{D_{+}}, N_{+}) \to C_{\perp}^{2,\sigma}(\overline{\Gamma_{1}}, N_{+}) \times \mathbb R^{k}$$
by
$$\Omega(f) = (f|_{\Gamma_{1}}, \int_{D_{+}}f^{*}\omega_{1},...,\int_{D_{+}}f^{*}\omega_{k}).$$
Thus we have
$$\Omega \circ F: C_{\perp}^{2,\sigma}(\overline{\Gamma_{1}}, N_{+}) \times K \to C_{\perp}^{2,\sigma}(\overline{\Gamma_{1}}, N_{+}) \times \mathbb R^{k}$$
and
$$\dfrac{d}{dt}\Big|_{t=0}\int_{D_{+}}[F(f_{0}|_{\Gamma_{1}}, t\kappa)]^{*}\omega = \dfrac{d}{dt}\Big|_{t=0}\int_{D_{+}}(f_{0} + t\kappa)^{*}\omega.$$
We can choose $\omega_{1},...,\omega_{k}$ such that the map $K \to \mathbb R^{k}$ defined by
$$\kappa \to \dfrac{d}{dt}\Big|_{t=0}\left(\int_{D_{+}}(f_{0} + t\kappa)^{*}\omega_{1},..., \int_{D_{+}}(f_{0} + t\kappa)^{*}\omega_{k}  \right)$$
is an isomorphism. \\
Therefore by the implicit function theorem, $\Omega \circ F$ is a diffeomorphism near $(f_{0}|_{\Gamma_{1}},0)$ and we have
$$\Omega \circ F(f_{0}|_{\Gamma_{1}},0) = (f_{0}|_{\Gamma_{1}}, 0,...,0).$$
By Theorem 3.4.2, $g^{-1}(0)$ is a codimension-$k$ submanifold of $C_{\perp}^{2,\sigma}(\overline{\Gamma_{1}}, N_{+}) \times K$ near $(f_{0}|_{\Gamma_{1}},0)$. Then $S = \Omega \circ F(g^{-1}(0))$ is a codimension-$k$ submanifold of $C_{\perp}^{2,\sigma}(\overline{\Gamma_{1}}, N_{+}) \times \mathbb R^{k}$. \\
Thus $\Omega$ is a bijective map from $\mathcal W = \{[f] \in \mathcal M: f \in W\}$ to a neighborhood of $\Omega(f_{0}) \in S$.\\
In order to show that $\{\mathcal W, \Omega\}$ forms an atlas for $\mathcal M$, we need to show that the transition maps are well-defined. \\
Let $\{\mathcal W_{1}, \Omega_{1}\}$ and $\{\mathcal W_{2}, \Omega_{2}\}$ be two coordinate charts. We have
$$\Omega_{1} \circ \Omega_{2}^{-1}|_{S_{2}} = \Omega_{1} \circ F_{2} \circ (\Omega_{2} \circ F_{2})^{-1}|_{S_{2}}.$$
This proves (1).\\
Conclusion (2) and (3) are directly from Theorem 3.4.2.
\hfill q.e.d.

\subsection{Mod-2 degree}
Under suitable assumptions, we can prove that $\Pi$ has a well-defined mod-2 degree.
\begin{thm}
  Assume $q \geq 3$. Let $\Pi: \mathcal M \to C_{\perp}^{2,\sigma}(\overline{\Gamma_{1}}, N_{+})$ be the projection map defined by $\Pi([f]) = f|_{\Gamma_{1}}$ and $g:[0,1] \to C_{\perp}^{2,\sigma}(\overline{\Gamma_{1}}, N_{+})$ be a $C^{1}$ embedding, then there is a $\tilde{g}: [0,1] \to C_{\perp}^{2,\sigma}(\overline{\Gamma_{1}}, N_{+})$ such that $\tilde{g}$ is arbitrarily close to $g$ and $\Pi$ is transversal to $\tilde{g}$.
\end{thm}
\emph{Proof:} By Theorem 3.4.3, $\mathcal M$ has a countable cover. Thus by Theorem 3.1 in Smale\cite{Smale1965}, there is a $\tilde{g}: [0,1] \to C_{\perp}^{2,\sigma}(\overline{\Gamma_{1}}, N_{+})$ such that $\tilde{g}$ is arbitrarily close to $g$ and $\Pi$ is transversal to $\tilde{g}$.
\hfill q.e.d.\\

\begin{thm}
  Assume $q \geq 3$. Let $\mathcal M^{'}$ be an open subset of $\mathcal M$ and $Z$ be a connected open subset of $C_{\perp}^{2,\sigma}(\overline{\Gamma_{1}}, N_{+})$. If $\Pi$ maps $\mathcal M^{'}$ properly into $Z$, then for a generic $\gamma \in Z$,
  $$\#(\Pi^{-1}(\gamma) \cap \mathcal M^{'}) \quad \text{mod 2}$$
  is a constant (1 or 0).
\end{thm}
\emph{Proof:} For any $\gamma_{0}, \gamma_{1} \in Z$, there is a $g:[0,1] \to Z$ which is a $C^{1}$ embedding such that $g(0) = \gamma_{0}$, $g(1) = \gamma_{1}$. \\
Then by Theorem 3.5.1, there is a $\tilde{g}: [0,1] \to Z$ such that $\tilde{g}$ is arbitrarily close to $g$ and $\Pi$ is transversal to $\tilde{g}$. \\
Therefore $\Pi^{-1}[\tilde{g}([0,1])] \cap \mathcal M^{'}$ is a one-dimensional manifold and we can choose $\tilde{g}(0), \tilde{g}(1)$ to be regular values of $\Pi$ such that
$$\partial \{\Pi^{-1}[\tilde{g}([0,1])] \cap \mathcal M^{'}\} = \Pi^{-1}\{\tilde{g}(0), \tilde{g}(1)\}.$$
By the assumption that $\Pi|_{\mathcal M^{'}}$ is proper, we know $\Pi^{-1}[\tilde{g}([0,1])] \cap \mathcal M^{'}$ is compact and
$$\#\Pi^{-1}\{\tilde{g}(0)\} + \#\Pi^{-1}\{\tilde{g}(1)\} = \#\partial \{\Pi^{-1}[\tilde{g}([0,1])] \cap \mathcal M^{'}\}$$
is an even number. Thus we have
$$\#\Pi^{-1}\{\tilde{g}(0)\} = \#\Pi^{-1}\{\tilde{g}(1)\} \text{   (mod $2$)}.$$
\hfill q.e.d.\\

Remark: Given $N = (\mathbb B_{+}^3, g)$, a three-dimensional Riemannian manifold which is diffeomorphic to a three-dimensional half ball with weakly mean convex boundaries consisting of the half sphere $S$ and the disk $\Sigma$ such that $S$ meets $\Sigma$ orthogonally along their common boundary.

If we consider the following spaces: \\
(1) $C_{\perp}^{2,\sigma}(\overline{\Gamma_{1}}, S) = \{\gamma:\Gamma_{1} \to S \, \big| \, \gamma $ is an embedding and $C^{2,\sigma}$ up to the boundary,  $\gamma(\Gamma_{1}) \subset S, \partial [\gamma(\Gamma_{1})] = \overline{\gamma(\Gamma_{1})} \cap \partial S, \gamma(\Gamma_{1}) \perp \Sigma\}$.\\
(2) $C_{*}^{2,\sigma}(\overline{D_{+}}, N) = \{f: D_{+} \to N \, \big|\, f$ is an embedding and $C^{2,\sigma}$ up to the boundaries and corners, $f(\Gamma_{1}) \subset S, f(\Gamma_{2}) \subset \Sigma, f(D_{+}) \subset N$, $f(D_{+})$ is transversal to $\partial N, f|_{\Gamma_{1}} \in C_{\perp}^{2,\sigma}(\overline{\Gamma_{1}}, S)\}$.\\
Two maps $f_{1}, f_{2} \in C_{*}^{2,\sigma}(\overline{D_{+}}, N)$ are said to be equivalent if $f_{1} = f_{2} \circ u$ for some $C^{2,\sigma}$ diffeomorphism $u: D_{+} \to D_{+}$ such that $u(x) = x$ for all $x \in \Gamma_{1}$. We use $[f]$ to denote the equivalence class of $f$.\\
(3) $\mathcal M = \{[f]: f \in C_{*}^{2,\sigma}(\overline{D_{+}}, N), H_{f} = 0, \Theta_{f} = 0\}$.

Then this is the case we introduce in the introduction and all theorems in this section still hold.

\end{document}